\documentclass[12pt]{amsart}

\usepackage{graphicx}
\usepackage{caption}
\usepackage{subcaption}

\usepackage{preamble_W}

\newcommand{\setOfReals}{\mathbb{R}}
\newcommand{\setOfNaturals}{\mathbb{N}}
\newcommand{\setOfNonnegativeIntegers}{{\mathbb{N}_0}}
\newcommand{\setOfPositiveReals}{{\setOfReals_{+}}}

\newcommand{\borel}[1]{\mathcal{B} (#1 )}

\newcommand{\spaceOfMeasures}[1]{\mathcal{M}(#1)}

\newcommand{\spaceOfOccupationMeasures}[2]{\mathcal{M}_{#1}(#2)}

\newcommand{\sys}[1]{\textsc{#1}}

\newcommand{\Bin}[2]{\sys{Binomial}(#1,#2)}

\newcommand{\Exponential}[1]{\sys{Exponential}(#1)}
\newcommand{\Unif}[2]{\sys{Uniform}[#1,#2]}
\newcommand{\Beta}[2]{\sys{Beta}(#1,#2)}

\newcommand{\measureIntegral}[2]{\langle#1, #2 \rangle }

\newcommand{\optionalVariation}[1]{ [#1] }
\newcommand{\predictableVariation}[1]{ \langle#1\rangle }

\newcommand{\absolute}[1]{\lvert #1 \rvert }

\newcommand{\leb}{\mathsf{Leb}}
\newcommand{\Wasserstein}[2]{\mathsf{W}_{#2}\left(#1\right)}
\newcommand{\boundedLipschitz}[1]{\mathsf{d}_{\mathsf{BL}}\left(#1\right)}
\newcommand{\Kantorovich}[1]{\mathsf{d}_{\mathsf{KR}}\left(#1\right)}

\newcommand{\TotalVariation}[1]{\mathsf{d}_{\mathsf{TV}}\left(#1\right)}

\newcommand{\cadlag}{c\`adl\`ag }

\newcommand{\defeq}{\coloneqq}

\newcommand{\indicator}[2]{\mathsf{1}_{#1}\left(#2\right)}

\newcommand{\myOperator}[1]{\mathsf{#1}}

\newcommand{\differential}[1]{\mathrm{d} #1}

\newcommand{\timeDerivative}[1]{\frac{\differential}{\differential t} #1 }

\newcommand{\eqstop}{.}
\newcommand{\eqcomma}{,}

\newcommand{\norm}[1]{\left\lVert#1\right\rVert}

\newcommand{\modulusOfcontinuity}{\mathsf{m}}

\newcommand{\PRM}{\zeta}
\newcommand{\compPRM}{\tilde{\zeta}} 

\newcommand{\E}{\mathsf{E}}
\newcommand{\Eof}[1]{\E\left[#1 \right]}

\newcommand{\prob}{\mathsf{P}}
\newcommand{\probOf}[1]{\prob\left(#1\right)}

\newcommand{\history}[1]{\mathcal{F}_{#1}  }

\newcommand{\ConvInProb}{\xrightarrow[]{ \hspace*{4pt}  \text{         P   } }}

\newcommand{\ConvInDist}{  \overset{\hspace*{4pt}  \mathcal{D}  }{\implies } }

\newcommand{\myExp}[1]{\exp \left( #1 \right)  }

\newcommand{\ie}{\textit{i.e.}}
\newcommand{\eg}{\textit{e.g.}}

\newcommand{\nX}{X^{(n)}}

\newcommand{\nU}{U^{(n)}}

\newcommand{\nS}{S^{(n)}}

\newcommand{\nI}{I^{(n)}}

\newcommand{\nR}{R^{(n)}}

\newcommand{\nK}{K^{(n)}}

\newcommand{\nM}{M^{(n)}}
\newcommand{\nm}{m^{(n)}}

\newcommand{\nmu}{\mu^{(n)}}

\newcommand{\vep}{\varepsilon}

\newcommand{\nerror}{\mathcal{E}^{(n)}}

\def\i{\iota}

\def\d{\mathrm{d}}

\def\nS{S^{(n)}}
\def\nI{I^{(n)}}
\def\nR{R^{(n)}}

\def\nX{X^{(n)}}

\title[Coloured epidemic models]{Coloured epidemic models: Functional Law of Large Numbers and Propagation of Chaos}
\author[K. Dutta]{Kushankur Dutta}
\address{Kushankur Dutta, 
School of Mathematical Sciences, 
University of Nottingham, 
University Park, 
Nottingham NG7 2RD, 
United Kingdom
}
\email{Kushankur.Dutta@nottingham.ac.uk}

\author[O. Izyumtseva]{Olga Izyumtseva}
\address{Olga Izyumtseva, 
School of Mathematical Sciences, 
University of Nottingham, 
University Park, 
Nottingham NG7 2RD, 
United Kingdom
}
\email{Olga.Iziumtseva1@nottingham.ac.uk}

\author[W. R. KhudaBukhsh]{Wasiur R. KhudaBukhsh}
\address{Wasiur R. KhudaBukhsh, 
School of Mathematical Sciences, 
University of Nottingham, 
University Park, 
Nottingham NG7 2RD, 
United Kingdom
}
\email{wasiur.khudabukhsh@nottingham.ac.uk}

\author[G. A. Rempa{\l}a]{Grzegorz A. Rempa{\l}a}
\address{Grzegorz A. Rempa{\l}a, 
Division of Biostatistics, 
College of Public Health, 
The Ohio State University, 
380E Cunz Hall, 
1841 Neil Ave., 
Columbus, OH 43210, 
USA}
\email{rempala.3@osu.edu}

\begin{document}
\maketitle

\begin{abstract}
In this paper, we study a stochastic \acf{SIR} model where the infection and the recovery rates depend on individual covariates for susceptibility and infectiousness of  the infector and the infectee. Such models allow  explicit nonlinearity in the incidence term. They are also important from a practical perspective, as they allow for the incorporation of individual heterogeneity into the epidemic process. Statistical estimates for crucial epidemiological parameters, such as the basic reproduction number, herd immunity threshold, could be vastly different, and even biased, when the population heterogeneity is ignored in the mathematical model. 

\vspace{10pt}

We describe our epidemic model as an \ac{IPS} of \acp{SDE} driven by \acp{PRM}. Our main mathematical contributions are a \acf{FLLN}, which approximates the empirical random measure of the \ac{IPS} by means of a deterministic measure-valued function,  and the propagation of chaos phenomenon, which establishes asymptotic independence of the particles as the population size goes to infinity with an explicit construction of McKean--Vlasov type Kac's ``nonlinear process''. We also briefly mention how the propagation of chaos phenomenon leads to a product-form likelihood function, which forms the basis of the so-called \acf{DSA} method for parameter inference based on sparse data. 
\end{abstract}


\tableofcontents

\section*{Preface}
Sara Del Valle (Los Alamos National Laboratory, USA), Joel Miller (La Trobe University, Australia), Rick Durrett (Duke University, USA), and the third author of this paper, Wasiur R. KhudaBukhsh (University of Nottingham, UK), organised a Banff International Research Station (BIRS) workshop titled ``Preparing for the next pandemic'' from June 12 to June 17, 2022. During the workshop, Gabriela Gomes (University of Strathclyde, UK) presented the  (deterministic) frailty models, which are popular in infectious disease epidemiology because of their ability to capture heterogeneity in susceptibility of individuals, and asked how they could arise from stochastic models. A rigorous derivation of such models from an individual-based stochastic description appears to be absent from the literature. This led to the development of the ``coloured epidemic models'', which include the frailty model as an example. 

The title of the paper is inspired by graph theory where the notion of \emph{coloured graphs} is used to distinguish graphs with heterogeneous vertices with individual features (colours) from ordinary graphs where vertices have no such distinguishing feature beyond their degrees (number of connections/neighbours). Perhaps a more descriptive title would be ``\acl{FLLN} and propagation of chaos for a heterogeneous \ac{SIR} model with individual covariates for susceptibility and infectiousness''. 

\section{Introduction}
\label{sec:introduction}

\subsection{Why does heterogeneity matter?}
\label{sec:heterogeneity}
Accounting for heterogeneity is important for infectious disease epidemiology. Although the standard argument for introducing heterogeneity in mathematical models has been their perceived greater biological realism, there is also a practical necessity. When the underlying heterogeneity is ignored, the estimates for key epidemiological parameters, such as the reproduction numbers, herd immunity threshold, could be drastically different, and even biased, leading to suboptimal public health policies, and potential loss of human lives. In \cite{Britton2020HerdImmunity}, the authors use a heterogeneous model and demonstrate that population heterogeneity can affect disease-induced immunity significantly because the proportion of infected individuals in groups with the highest contact rates is greater than that in groups with low contact rates. 

In \cite{Andersson1998heterogeneity}, the authors consider an epidemic model with varying susceptibilities, and compare the asymptotic final epidemic size with the corresponding size for the same epidemic model on a  homogeneous population based on an earlier work of Ball and colleagues \cite{Ball1993FinalSizeGenStoch,Ball1985detStoch}. Their main finding is that while the largest epidemic arises in heterogeneous populations for less contagious diseases, the homogeneous population produces larger epidemic sizes when the disease is very contagious. Interestingly, the probability of a large outbreak is \emph{always} minimal in the homogeneous setup regardless of the contagiousness of the disease. In \cite{Ball1995finalsize_het_infectivity}, the authors study the final size of an epidemic amongst a closed homogeneously mixing population with different types of infective. In \cite{Neal2007coupling}, Neal considers an epidemic model with variable susceptibilities and infectivities, called the \emph{variable generalised stochastic epidemic}, and shows that the final size of the epidemic is the same for the random graph epidemic model, which is also heterogeneous model. The authors in \cite{Ball2007infector_dependent_severity} consider an epidemic model where infected individuals have different severities of the disease. They provide large population limits and also discuss the effect of vaccination on the final size of the epidemic. Heterogeneity in infectivity is considered in \cite{Forien2021Varying,Pang2022FCLT,Forien2025immunity}. In a recent work \citet{Xue_2025Scaling}, Xue considers individual infection and recovery rates on a complete graph model, and proves large-graph (hydrodynamic) limits and also quantifies the fluctuations in terms of a generalised \ac{OU} process. For a similar model, the question of phase transitions is considered in \cite{Xue2018PhaseTransitions}. 

Although the focus of this paper is on stochastic models, there is also a vast amount of literature on the deterministic side, which corroborates the importance of studying heterogeneity in the context of epidemics. 
In \cite{Montalban2022HerdImmunity}, the authors consider heterogeneity in susceptibility in a deterministic compartmental \ac{SEIR} model and study the impact of heterogeneity on herd immunity threshold. In \cite{Margheri2015Correlation}, the authors study an endemic \ac{SIS} model and test the impact of the  second centred moment of the population susceptibility on the infection prevalence. We refer the interested readers to \cite{Seymour2022BNP,Kotounou2026Global,Margheri2017Vaccine,Rodrigues2009Reinfection,Inaba2001Kermack,Novozhilov2008heterogeneous} and references therein for other recent studies of heterogeneity in the deterministic context. In a recent paper \cite{Izyumtseva2026Sellke}, the authors provide a derivation of the deterministic frailty model (heterogeneity in susceptibility) from an individual-based stochastic model using the Sellke construction \cite{AnderssonBritton2000}. It is shown in \cite{Izyumtseva2026Sellke} that the nonlinearities in the incidence term (power law, logarithmic etc.) in a large class of \ac{ODE}-based \ac{SIR} models are in fact a consequence of individual heterogeneity in the population. 

The present work is closely related to \cite{Izyumtseva2026Sellke}, which considers heterogeneity in susceptibility and uses the Sellke construction to obtain an explicit reduction in terms of cumulative infection pressure and the Laplace transform of the susceptibility distribution. Here, we allow the transmission rate to depend on covariates of both the susceptible and the infectious individuals through a general kernel $\beta(\cdot, \cdot)$, as well as allow covariate-dependent recovery. The frailty model considered in \cite{Izyumtseva2026Sellke} is, therefore, an important special case of the present framework, while the more general formulation naturally leads to a
measure-valued limit and to the propagation-of-chaos result, which additionally facilitates a likelihood-based inference.









\subsection{Our contribution}
In this paper, we consider an individual-based stochastic \ac{SIR} model in terms of a system of \acp{SDE} driven by independent \acp{PRM}. Each individual is endowed with covariates that affect their susceptibility, also called frailty, and infectiousness. We prove an \ac{FLLN} that shows that the empirical random measure of the epidemic process converges weakly to a deterministic, continuous (probability) measure-valued function, which is characterised as a weak solution to a differential equation in \Cref{thm:flln}. We use a tightness-uniqueness argument to prove the \ac{FLLN}. To be more precise, we first show that the collection of empirical random measures (indexed by varying population sizes) is relatively compact in the space of all probability measures on the space of finite (probability) measure-valued \cadlag functions. We use a tightness criterion from \cite{Dawson1993Measure_valued} for this purpose. Next, we show that the limit points of convergent subsequences are unique by showing that an integral equation, the candidate limit point, admits a unique solution. The weak \ac{FLLN} in \Cref{thm:flln} is further strengthened to $L^1(\Omega, \history{}, \prob)$ convergence in \Cref{thm:wasserstein_L1_conv}. We use the bounded Lipschitz distance (also called the Dudley metric; see \Cref{example:DudleyMetric} in Appendix \ref{sec:IPM}) and the fact that uniformly bounded and equicontinuous functions on a compact subset of Euclidean space are totally bounded to show that the $L^1(\Omega, \history{}, \prob)$ distance between the empirical random measure and the limiting deterministic, continuous measure-valued function vanishes as the population size increases to infinity, when the space of finite measures on the space of covariates is equipped with the bounded Lipschitz distance.


Convergence of empirical random measures to a deterministic measure is intimately related to the phenomenon of propagation of chaos \cite{Sznitman1991Topics,Hauray2014KacChaos}, a concept that originated in physics, and was popularised by Kac among analysts and probabilists. Propagation of chaos for an \ac{IPS} implies asymptotic independence of the particles. As such, it has important applications in statistical inference since it allows one to write a product-form likelihood function, which is advantageous from a practical standpoint. Because of its relevance to parameter inference via the \ac{DSA} method, we prove the propagation of chaos phenomenon in the  $L^1(\Omega, \history{}, \prob)$ sense in \Cref{thm:L1_prop_chaos}. We use the coupling method to perform a pathwise analysis of the \ac{IPS} as a system of interacting \acp{SDE} and a system of \ac{iid} \acp{SDE}, akin to the McKean--Vlasov equation in the probability literature. 

Even though we focus on the \ac{SIR} model, our results can be extended to other compartmental models, such as the \ac{SEIR}, \ac{SIS} models. The \ac{SEIR} model is often regarded as more realistic than the \ac{SIR} model since it accounts for the incubation period of the disease. However, as shown in \cite{KhudaBukhsh2024HowTo}, the \ac{SEIR} model can be approximated to any desired level of accuracy by an \ac{SIR} model with time-varying coefficients. In \cite{Izyumtseva2024RandomEffects}, the authors consider an \ac{SIR} model with (dynamic) random effects that affect the infection rates for all individuals. However, we note that the model considered in \cite{Izyumtseva2024RandomEffects} is different from the one we consider in this paper, since it does not allow for individual heterogeneity. It is worth mentioning that our model can be viewed as an extension of the multitype epidemic models \cite{Ball1993FinalSizeGenStoch,Neal2006Multitype}, \cite[Section 3]{Andersson1998heterogeneity} and is also closely related to Kenah's pairwise model \citep{kenah2010contact,kenah2015semiparametric,kenah2013nonparametric,Sharker2024pairwise} based on contact intervals. 


\subsection{Structure of the paper}
The rest of the paper is structured as follows. We conclude this section with a list of notations and notational conventions. In \Cref{sec:stoch_model}, we describe our stochastic epidemic model with covariates in terms of a pure jump Markov process constructed as a solution to a system of \acp{SDE} driven by \acp{PRM} with Lebesgue intensities. We provide a simulation algorithm and a number of specific examples that could be useful for practical purposes. In \Cref{sec:mean_field_limit}, we present a weak \ac{FLLN} for the  empirical (random) measure of the epidemic model. The \ac{FLLN} for the empirical measure is closely related to the notion of propagation of chaos, and its use in the so-called \ac{DSA} method of parameter inference. We discuss this connection and prove a propagation of chaos result in  $L^1(\Omega, \history{}, \prob)$  (\Cref{thm:L1_prop_chaos}) in \Cref{sec:propagation_of_chaos} and \Cref{sec:parameter_inference}. 
A short conclusion is presented in \Cref{sec:conclusion}. 
For the sake of completeness, and to appeal to the wider analysis and applied mathematics communities, we provide additional proofs and extensive background material in the Appendices  \ref{sec:additional}, \ref{sec:IPM}, \ref{sec:limits_of_stoch_pro}, and \ref{sec:prop_chaos_appendix}. A list of acronyms used in this paper is provided in Appendix \ref{sec:acronyms}.


\subsection{Notations}
\label{sec:notations}
The symbol $\setOfReals$ denotes the set of real numbers. 
The sets of natural numbers, non-negative integers, non-negative real numbers are denoted by $\setOfNaturals, \setOfNonnegativeIntegers$, and $\setOfPositiveReals$, respectively. The indicator function (also known as the \emph{characteristic function} in mathematical analysis) of a set $A$ will be denoted by $\indicator{A}{\cdot}$, \ie, $\indicator{A}{x} = 1$ if $x\in A$, and zero otherwise. 
The space of continuous functions from a metric space $E$ to a metric space $F$ will be denoted by $C(E, F)$ with the subset $C_b(E, F)$ containing the bounded, continuous functions when $F$ is a subset of $\setOfReals$. The space $D([0,T], F)$ denotes  the space of \cadlag functions (those that are right continuous and have finite left-hand limits) from the interval $[0,T]$ to $F$, with $0 < T < \infty$. The definition of $D([0, \infty), F)$ follows similarly. In our case, the spaces $E$ and $F$ will be complete, separable metric spaces. The space $C(E, F)$ will be equipped with the supremum norm, whereas the space $D([0,T], F)$ will be equipped with the Skorokhod topology (see \cite[Chapter 3]{Ethier:1986:MPC}, \cite{Billingsley1999Convergence}, or \cite[Chapter 12]{Whitt2002StochLimits}), unless stated otherwise. Sometimes, it will be convenient to use the supremum norm for both spaces. For a \cadlag function $f$, we denote the left-hand limit of $f$ at $t$ by $f(t-)$.   For two real numbers $a$ and $b$, the symbols $a \wedge b$ and  $a \vee  b$ denote, respectively, the minimum and maximum of  $a$ and $b$. 

The  Borel $\sigma$-field on a metric or a topological space $E$ is denoted by $\borel{E}$. The set $\spaceOfMeasures{E}$ will denote the space of finite (non-negative) measures on $E$ equipped with the topology of weak convergence \citep{Billingsley1999Convergence}. For $r>0$, $\spaceOfOccupationMeasures{r}{E} \subset \spaceOfMeasures{E}$ will denote the space of (non-negative) measures $\nu$ on $E$ such that $\nu(E) =r$. Therefore, the set $ \spaceOfOccupationMeasures{1}{E}$ will denote the space of probability measures on $E$. We use $\leb$ to denote the Lebesgue measure on $\setOfReals$.  Given a 
measure $\mu$ on a measurable space $(U, {\mathcal U})$ and a  function $g:U\to \setOfReals$ that is integrable  with respect to $\mu$, we denote the integral of $g$ with respect to the measure $\mu$ by $\measureIntegral{g}{\mu}$. That is, $\measureIntegral{g}{\mu} = \int g \differential{\mu} = \int_{U} g(x) \mu(\differential{x}) $. We will not mention the set $U$ explicitly whenever it is clear from the context. 
    
For a Polish metric space $(E, d)$, the Vaser\v{s}te\u{\i}n (Wasserstein in Western form) distance of order $p$ between two probability measures $\mu, \nu \in \spaceOfOccupationMeasures{1}{E}$ will be denoted by $\Wasserstein{\mu, \nu}{p}$. Several integral probability metrics will be considered in this paper. For $\mu, \nu \in \spaceOfOccupationMeasures{1}{E}$, the notations $\Kantorovich{\mu, \nu}, \TotalVariation{\mu, \nu}, $ and $\boundedLipschitz{\mu, \nu}$ will denote the Kantorovich--Rubinstein distance, the total variation distance, and the bounded Lipschitz distance, respectively. See Appendix \ref{sec:IPM} for the definitions of these distances. \cite{Villani2009optimaltransport} is an excellent resource for optimal transport in general and \cite{Mueller1997IntegralProbMetrics,Zolotarev1983ProbabilityMetrics,Rachev1991ProbabilityMetrics}, about integral probability metrics.

\section{Stochastic model} 
\label{sec:stoch_model}
Our stochastic model is an extension of the standard \ac{CTMC}-based \ac{SIR} model \citep{AnderssonBritton2000}. We allow the infection and the recovery rates to depend on  individual covariates. We will assume that the covariates do not change over time. Further, we will assume the covariates can be described by random variables taking values in a set $\Theta$. For practical purposes, we will assume $\Theta$ is a compact subset of the Euclidean space $\setOfReals^d$ for some $d\in \setOfNaturals$. In order to keep the notations simple, we will not distinguish between covariates that affect infectivity and susceptibility separately, and will assume that the individual covariates affect both. This is not a problem and does not imply lack of generality, since we could always divide the covariates into two (not necessarily disjoint) subsets, and define the infection and recovery rates as functions of the two subsets separately. For example, see \cite{Sharker2024pairwise,kenah2010contact}.

Let $(\Omega, \history{}, \prob)$ be a given, large enough probability space. For instance, the probability space $([0, 1], \borel{[0, 1]}, \leb)$ will be sufficiently large for our purposes. All random variables (also  random elements) will be defined on this space. 
Let us assume that the population size is $n$. Let $\nX_i (t) \defeq (\nS_i(t), \nI_i(t), \nR_i(t))$ denote the state of the $i$-th individual at time $t$, where the stochastic processes $\nS_i(t), \nI_i(t), \nR_i(t)$ are defined as follows: $\nS_i(t)$ is one if the $i$-th individual is susceptible at time $t$, and zero otherwise. Similarly, $\nI_i(t)$ is one if the $i$-th individual is infected at time $t$, and zero otherwise. Finally, $\nR_i(t) = 1- \nS_i(t) - \nI_i(t)$ for all $t\ge 0$. At any point in time, an individual can have only one of the three immunological statuses: susceptible, infected, or recovered. The stochastic processes are superscripted by $n$ to make their dependence on the total population size $n$ explicit. Furthermore, we assume that the $i$-th individual is endowed with covariates $\theta_i \in \Theta$, for $i=1, \ldots, n$. We assume that the covariates are \ac{iid} random variables. They are also assumed independent of the initial condition $(\nX_1(0), \nX_2(0), \cdots, \nX_n(0))$. Once assigned to the individuals at time $t=0$, the covariates do not change over time. Let 
\begin{align*}
    \nX(t) & \defeq (\nX_1(t), \nX_2(t), \ldots, \nX_n(t))\eqcomma \text{ for } t\ge 0\eqstop 
\end{align*}


The dynamics is entirely driven by the infected individuals in the sense that the dynamics stop (\ie, become constant) if the number of infected individuals becomes zero. Conditioned on $(\theta_1, \theta_2, \dots, \theta_n)$,  if the $i$-th individual is infected at time $t$, then the $i$-th individual makes infectious contact with a susceptible individual $j$ with covariates $\theta_j$ at rate $n^{-1}\beta(\theta_j, \theta_i)$, for some known function $\beta: \Theta\times \Theta \to \setOfPositiveReals$. A contact with an infected individual is sufficient to cause infection in a susceptible individual. The infected individual $i$ recovers (or is removed) at rate $\gamma(\theta_i)$, for some known function $\gamma:\Theta \to \setOfPositiveReals$. Once recovered, the individual $i$ plays no further role in the epidemic.  For simplicity, 
we will assume $ \absolute{\beta(u, v)} \le \hat{\beta} \eqcomma \text{ and } \absolute{\gamma(u)} \le \hat{\gamma} \eqcomma $ for some positive constants $\hat{\beta}$ and $\hat{\gamma}$, and for all $u, v \in \Theta$. 


\subsection{\aclp{SDE} driven by \aclp{PRM}}

We will describe the trajectories of the individuals using a system of \acp{SDE} driven by \acp{PRM}. Readers unfamiliar with the theory of stochastic calculus with respect to \acp{PRM} (or \acp{PPM}) are referred to \citet[Chapter 4]{Applebaum_2009Levy}, \citet[Chapter IV]{IkedaWatanabe2014Stochastic}, or \cite{Bremaud2020PointProcessCalculus}. The lectures by Last and Penrose \citep{Last2018Lectures} are also a good introduction to \acp{PPM}. For the general theory of random measures, \citet{Kallenberg2017RandomMeasures} is an excellent resource. 

Let $\PRM_{i}^{(1)}$ and $\PRM_{i}^{(2)}$ for $i=1, 2, \ldots, n,$ be independent \acp{PRM} on $\setOfPositiveReals\times \setOfPositiveReals$ with intensity $\leb\times \leb$, where $\leb$ is the Lebesgue measure on $\setOfPositiveReals$. We equip the probability space $(\Omega, \history{}, \prob)$ with the  filtration $\{\history{t} : t\ge 0\}$ generated by  
\begin{align}
    \history{t} \defeq \sigma ((\nX_i(0), \theta_i), \PRM_{i}^{(1)}((0, s]\times A), \PRM_{i}^{(2)}((0, s]\times B): i \in \{1, \ldots, n\}; A, B \in \borel{\setOfPositiveReals}; s\le t) \eqcomma 
\end{align}
for $t\ge  0$. Note that the initial data $\{(\nX_i(0), \theta_i): i=1, \ldots, n\}$ is measurable with respect to (included in) $\history{0}$. Furthermore, we include all $\prob$-null sets in  $\history{0}$. The filtration $\{\history{t} : {t\ge 0}\}$ is right continuous 
    \begin{align*}
        \history{t+} \defeq \cap_{s>0}\history{t+s} = \history{t}\eqstop 
    \end{align*}
Therefore, the filtered probability space $(\Omega, \history{}, (\history{t})_{t\ge 0}, \prob) $  is complete. In stochastic analysis literature, this is often referred to as the \emph{usual  conditions} or the Dellacherie conditions (see \cite[Definition 2.25]{karatzas1991brownian} or \cite[Definition 1.3]{jacod2003limit}). 


For the $i$-th individual, we assume that the stochastic process $\nX_i$ solves the following \ac{SDE}: 
\begin{align}
    \begin{aligned}
        \nS_i(t) &= \nS_i(0) -  
        \int_{(0, t] \times \setOfPositiveReals} \indicator{[0, \nS_{i}(s-) n^{-1} \sum_{j=1}^{n}\beta(\theta_i, \theta_j)  \nI_j(s-) ]}{v}\PRM_{i}^{(1)}(\differential{s}\times\differential{v})\eqcomma \\
        \nI_i(t) &=  \nI_{i}(0) + 
        \int_{(0, t] \times \setOfPositiveReals} \indicator{[0, \nS_{i}(s-) n^{-1} \sum_{j=1}^{n}\beta(\theta_i, \theta_j)  \nI_j(s-) ]}{v}\PRM_{i}^{(1)}(\differential{s}\times\differential{v})\\
        &{}\quad \quad  -  \int_{(0, t] \times \setOfPositiveReals} \indicator{[0, \gamma(\theta_i) \nI_{i}(s-) ]}{v}\PRM_{i}^{(2)}(\differential{s}\times\differential{v})\eqcomma \\
        \nR_i(t) & = 1- \nS_i(t) - \nI_i(t)\eqcomma \quad \text{ for } t\ge 0\eqstop 
    \end{aligned}
    \label{eq:sde_for_i_th_individual}
\end{align} 
Let $\Delta_3 \defeq \{ e_1, e_2, e_3\}$, where $e_1 = (1, 0, 0), e_2 = (0, 1, 0), e_3 = (0, 0, 1)$ are the standard basis vectors in $\setOfReals^3$. Then, $\nX_i(t) \in \Delta_3$ for all $i=1, \ldots, n$ and $t\ge 0$. In particular, $\nX_i(t) = e_1$, $\nX_i(t) = e_2$, and $\nX_i(t) = e_3$ correspond to the $i$-th individual being in the susceptible, infected, and recovered states at time $t$, respectively. Naturally, the above system of \acp{SDE} can also be written as 
\begin{align}
    \begin{aligned}
        \nX_i(t) &= \nX_i(0) + (e_2-e_1) 
        \int_{(0, t] \times \setOfPositiveReals} \indicator{[0, \nS_{i}(s-) n^{-1} \sum_{j=1}^{n}\beta(\theta_i, \theta_j)  \nI_j(s-) ]}{v}\PRM_{i}^{(1)}(\differential{s}\times\differential{v})\\ 
        &{}\quad \quad  + (e_3-e_2) \int_{(0, t] \times \setOfPositiveReals} \indicator{[0, \gamma(\theta_i) \nI_{i}(s-) ]}{v}\PRM_{i}^{(2)}(\differential{s}\times\differential{v})\eqcomma \quad \text{for } i =1, 2, \ldots, n\eqstop  
    \end{aligned}
    \label{eq:combined_sde_for_i_th_individual}
\end{align}
If the covariates $\theta_1, \theta_2, \ldots, \theta_n$ were deterministic, the solution $(\nX_1, \nX_2, \ldots, \nX_n)$ to the system of \acp{SDE} in \eqref{eq:combined_sde_for_i_th_individual} is a Markov process with respect to the filtration $\{\history{t} : t\ge 0\}$, \ie, $\Eof{f(\nX(t+s)) \mid \history{t}} = \Eof{f(\nX(t+s)) \mid \nX(t)}$ \ac{a.s.} for all bounded measurable functions $f: \Delta_3^n \to \setOfReals$. However, this is not true when the covariates $\theta_1, \theta_2, \ldots, \theta_n$ are random variables. However, the following is true:
\begin{align*}
    \Eof{f(\nX(t+s)) \mid \history{t}} = \Eof{f(\nX(t+s)) \mid \nX(t), (\theta_1, \theta_2, \ldots, \theta_n)}\eqcomma \quad \ac{a.s.}\eqcomma 
\end{align*}
for all bounded measurable functions $f: \Delta_3^n \to \setOfReals$. 
As such, define the (random) operator $\myOperator{A}_n$ on measurable functions $f: \Delta_3^n \to \setOfReals$ as 
\begin{align}
    \begin{aligned}
    \myOperator{A}_n f(x_1, x_2, \ldots, x_n) &\defeq \sum_{i=1}^{n} \sum_{j=1}^{n} \frac{\beta(\theta_i, \theta_j)}{n}x_i^{(1)} x_j^{(2)}\left( f(x_1, x_2, \ldots, x_i - e_1 +e_2, \ldots, x_n) \right. \\
    &\quad \quad \quad \quad \quad 
    \left. - f(x_1, x_2, \ldots, x_n)  \right) \\
    &\quad 
    + \sum_{i =1}^{n} {\gamma(\theta_i)}x_i^{(2)}\left( f(x_1, x_2, \ldots, x_i - e_2 +e_3, \ldots, x_n) \right. \\
    &\quad \quad \quad \quad \quad 
    \left. - f(x_1, x_2, \ldots, x_n)  \right) \eqcomma 
    \end{aligned}
\end{align}
where $x_i \defeq (x_i^{(1)}, x_i^{(2)}, x_i^{(3)}) \in \Delta_3$ for $i=1, \ldots, n$. Note that the operator $\myOperator{A}_n$ is not a Markov generator itself (since it is random). Nevertheless, the stochastic process 
\begin{align}
    \begin{aligned}
        \nK_f(t) & \defeq f(\nX_1(t), \nX_2(t), \ldots, \nX_n(t)) - f(\nX_1(0), \nX_2(0), \ldots, \nX_n(0)) \\
        &\quad \quad 
        - \int_{0}^{t} \myOperator{A}_n f(\nX_1(s), \nX_2(s), \ldots, \nX_n(s)) \differential{s}
    \end{aligned}
    \label{eq:Dynkin_martingale_individual}
\end{align}
is a zero-mean local martingale with respect to the filtration $\{\history{t}: t\ge 0\}$ for each bounded measurable function $f: \Delta_3^n \to \setOfReals$. In fact, in view of our assumptions on the functions $\beta$ and $\gamma$, and the fact that the covariates $\theta_1, \theta_2, \ldots, \theta_n$ are square integrable (since $\Theta$ is assumed compact), one can show that $\nK_f$ is indeed a true martingale and not just a local martingale. To see that, note that for each $t\ge 0$, the collection $\{\nK_f(t \wedge \tau_k) : k \ge 1\}$ is bounded in $L^2(\Omega, \history{}, \prob)$ for any sequence of stopping times $\{\tau_k: k \ge 1\}$. The collection is therefore uniformly integrable. Therefore, by virtue of  \citet[Proposition 1.8 ]{ChungWilliams90}, the stochastic process $\nK_f$ is indeed an $\history{t}$-martingale, and not just a local martingale.

Before we proceed to study the limit theorems for this stochastic model in \Cref{sec:mean_field_limit}, and \Cref{sec:propagation_of_chaos}, we describe a simulation algorithm for practical purposes, and mention some examples. 

\subsection{Simulation algorithm}
\label{sec:simulation_algorithm} 
Since the covariates do not change over time, and the processes $\nK_f$ in \eqref{eq:Dynkin_martingale_individual} are zero-mean martingales, the trajectories of $(\nX_1, \nX_2, \ldots, \nX_n)$ can be simulated simply by adapting the standard Doob--Gillespie algorithm \citep{Anderson:2011:CTM,Wilkinson2018SMS}, or the jump chain--holding times construction of a \ac{CTMC} \citep{Norris:1998:MarkovChains,Darling2008Differential}. For the readers' convenience, a pseudocode is presented in \Cref{alg:Gillespie}. 

\begin{algorithm}[tbh]
\small\topsep=0in\itemsep=0in\parsep=0in
\normalsize
\begin{algorithmic}[1]
  \caption{%
  \small Pseudocode for the Doob--Gillespie algorithm}\label{alg:Gillespie}
  \Require{The functions $\beta, \gamma$; and a final time $T$}
  \State Initialize $(\nX_i(0), \theta_i) \equiv ((\nS_i(0), \nI_i(0), \nR_i(0)), \theta_i)$ for $i=1, \ldots, n$
  \While{$t<T$}
  \State Calculate rates $\lambda_{i, j}^{(1)} =\beta(\theta_i, \theta_j) \nS_i(t) \nI_j(t)/n$  and $\lambda_{i}^{(2)}=\gamma(\theta_i) \nI_i(t)$ for $i=1, \ldots, n$ and $j \ne i$
  \State Draw $E_{i, j}^{(1)}$ from $\Exponential{\lambda_{i, j}^{(1)}}$ for $i=1, \ldots, n$ and $j \ne i$
  \State Draw $E_{i}^{(2)}$ from $\Exponential{\lambda_{i}^{(2)}}$ for $i=1, \ldots, n$
  \State Set   next  transition  time  $\Delta t = \min\{E_{i, j}^{(1)}, E_{i}^{(2)} : i=1, \ldots, n;  j \ne i\}$
  \State Update time $t \leftarrow t + \Delta t$
  \If {$\Delta t = E_{i, j}^{(1)}$ for some $(i, j)$}
  \State Update $(\nX_i(t) ) \leftarrow (\nX_i(t-)- e_1 + e_2 )$
  \ElsIf {$\Delta t = E_{i}^{(2)}$ for some $i$}
  \State Update $(\nX_i(t) ) \leftarrow (\nX_i(t-)- e_2 + e_3 )$
  \EndIf 
  \EndWhile
\end{algorithmic}
\end{algorithm}


\subsection{Examples}
\label{sec:examples} 

\begin{myExample}[Standard \ac{SIR} model]
    The standard \ac{SIR} model can be realised by setting $\beta(u, v) = \beta_{\star}$ for some $\beta_{\star} >0$. This model is now a classical topic included in many textbooks, \eg, \cite{AnderssonBritton2000}. 
\end{myExample}

\begin{myExample}[Markovian multitype \ac{SIR} model]
    Assume $\Theta = \{1, 2, \ldots, K\}$ for some positive integer $K$. Then, our model corresponds to the Markovian multi-type \ac{SIR} model for any $\beta : \Theta^2 \to \setOfPositiveReals$, and $\gamma : \Theta \to \setOfPositiveReals$. For instance, one could choose $\beta(\theta, \theta') = \beta_\star \theta \theta'$ and $\gamma(\theta) = \gamma_\star \theta$.  See \cite{Ball1995finalsize_het_infectivity,Ball1985detStoch,Andersson1998heterogeneity,Neal2006Multitype}. 
\end{myExample}

\begin{myExample}[\ac{SIR} model on random geometric graph]
    We can study \ac{SIR} epidemics on random geometric graph \citep{Penrose2003RandomGeometricGraphs} by choosing the covariates as spatial points of a point process and then setting  
\begin{align*}
    \beta(\theta_i, \theta_j) = \beta_{\star}\indicator{[0, r]}{\norm{\theta_i- \theta_j}}\eqcomma 
\end{align*}
for some $\beta_{\star}>0$ and $r>0$. In this case, an infected individual can infect a susceptible individual if and only if their physical distance is at most some predefined value $r$. 
\end{myExample}

\begin{myExample}[Exponentially decaying infection rate]
    Instead of entirely disallowing infectious contacts between individuals whose physical distance exceeds a predefined value, it might be desirable to let the infection rate decay exponentially with the distance between an infected individual and a susceptible individual. We can achieve this by setting  
    \begin{align*}
    \beta(\theta_i, \theta_j) = \beta_{\star}\myExp{-\norm{\theta_i- \theta_j}}\eqcomma 
\end{align*}
for some $\beta_{\star}>0$. In this case, the underlying graph is a complete graph, but the infection rate decays exponentially with the distance between the two individuals. 
\end{myExample}

\begin{myExample}[Frailty model]
    The frailty models postulate that the infection rate varies by susceptibility. Therefore, we can  consider the frailty model by setting $\beta(\theta_i, \theta_j) = \beta_{\star} \tilde{\beta}(\theta_i) $ for some bounded function $\tilde{\beta} : \Theta \to [0, \infty)$ and $\beta_{\star}>0$. The most common choice is the linear function $\tilde{\beta}(x) =  x$. We refer the readers to \cite{gomes2022individual,Gomes2014,Montalban2022HerdImmunity,Margheri2015Correlation,Izyumtseva2026Sellke} for an overview of this type models. 
    
\end{myExample}

\begin{myExample}[Epidemics on correlated graphon random graph]\label{example:graphon}
    We can also model epidemics on a random graph sampled from a given graphon. Graphons are fundamental objects describing limits of dense (random) graphs \citep[Chapter 7]{Lovasz2012LargeNetworks}. They are useful for applications to mean-field games and epidemic models \cite{Cui2022Chaos}. Assume that the set $\Theta = \Theta_1\times \Theta_2\times \Theta_3$ is a compact subset of $\setOfReals^{d}\times \setOfReals \times 
    \setOfReals$, for some $d\ge 1$. Write $\theta_i = (\theta_i^{(1)}, \theta_i^{(2)}, \theta_i^{(3)})$, where $\theta_i^{(1)} \in \Theta_1\subset \setOfReals^d$ and $ \theta_i^{(2)}, \theta_i^{(3)} \in \setOfReals$, for $i=1, 2, \ldots, n$.  Assume the random variables $\theta_1^{(j)}, \theta_2^{(j)}, \ldots, \theta_n^{(j)}$ are \ac{iid} for each $j=1, 2, 3$. For each $i=1, 2, \ldots, n$, the random variables $\theta_i^{(1)}, \theta_i^{(2)},$ and $\theta_i^{(3)}$ are also assumed independent. Furthermore, assume the random variables are continuous. Let 
    \begin{align*}
        F(x) \defeq \probOf{\theta_1^{(2)} \theta_2^{(2)} \le x }\eqcomma x \in \setOfReals\eqcomma 
    \end{align*}
    denote the \ac{CDF} of the random variable $\theta_1^{(2)}\theta_2^{(2)}$. 
    Define  
    \begin{align*}
        \beta(\theta_i, \theta_j) = \indicator{[0, w(\theta_i^{(3)}, \theta_j^{(3)})]}{F(\theta_i^{(2)}\theta_j^{(2)}) } \tilde{\beta}(\theta_i^{(1)}, \theta_j^{(1)})\eqcomma 
    \end{align*}
    where $w : \Theta_3^2 \to [0, 1]$ is a graphon \citep[Chapter 7]{Lovasz2012LargeNetworks}, and $\tilde{\beta} : \Theta_1^2 \to \setOfPositiveReals$ is a function such that $\tilde{\beta}(u, v) \le \hat{\beta}$. Since $F(\theta_i^{(2)}\theta_j^{(2)})\sim \Unif{0}{1}$, for all $i\ne j$ (by the probability integral transform), the conditional probability that there is an edge between $i$ and $j$, conditioned on $\theta_i^{(3)}, \theta_j^{(3)}$, is $w(\theta_i^{(3)}, \theta_j^{(3)})$. That is, the random graph is precisely a sample from the graphon $w$. Conventionally, one chooses $\Theta_3 = [0, 1]$ so that the graphon is defined on the unit square $[0, 1]^2$. However, note that, in our construction, the edges between two pairs $(i, j)$ and $(i,k)$ are not independent, unlike classical random graphs sampled from graphons. 
\end{myExample}

\begin{myExample}[Epidemics on correlated \ac{ER} random graphs]\label{example:ErdosRenyi}
    \ac{ER} random graphs are the most basic class of random graphs for which many analytic computations are feasible \cite{vanDerHofstad2017RGCNvol1}. 
    Setting $w(\theta, \theta')=p$ for all $\theta, \theta'$, and  for some $p\in [0, 1]$ in \Cref{example:graphon} yields an epidemic model on a correlated \ac{ER} random graph on $n$ vertices with connection probability $p$. As with \Cref{example:graphon}, the edges between  two pairs $(i, j)$ and $(i,k)$ are not independent in this construction, unlike classical \ac{ER} random graphs. 
\end{myExample}


\section{Mean-field limit}
\label{sec:mean_field_limit}
In this section, we will prove an \ac{FLLN} for the system of \acp{SDE} in \eqref{eq:combined_sde_for_i_th_individual}. Let us define the empirical (random) measure 
\begin{align}
    \nmu_t (A\times B) \defeq \frac{1}{n} \sum_{i=1}^{n} \delta_{(\nX_i(t), \theta_i)}(A\times B) \eqcomma 
\end{align}
for $t\ge 0$ and $A\times B \in \borel{\Delta_3 \times \Theta}$, the Borel $\sigma$-field of subsets of $\Delta_3 \times \Theta$. Note that the trajectories of the stochastic process $\{\nmu_t: t\ge 0\}$ lie in $D([0, \infty),  \spaceOfOccupationMeasures{1}{\Delta_3 \times \Theta})$, the space of \cadlag functions on $[0, \infty)$ with values in $\spaceOfOccupationMeasures{1}{\Delta_3 \times \Theta}$, the space of probability measures on the set $\Delta_3 \times \Theta$. To be more precise, the values of the process $\{\nmu_t: t\ge 0\}$  are probability measures on $\Delta_3 \times \Theta$ of the form $n^{-1} \sum_{i=1}^{n} \delta_{(x_i, y_i)}$. Let us denote the subset of such probability measures by $\spaceOfOccupationMeasures{1, \text{point}}{\Delta_3 \times \Theta}$.  In fact, the stochastic process $\nmu$ is a Markov process. In order to describe the generator of this Markov process, let us define the map $\lambda_{\phi} : \left(\Delta_3\times \Theta\right)\times \spaceOfOccupationMeasures{1}{\Delta_3 \times \Theta} \to \setOfReals$ as follows 
\begin{align}
    \begin{aligned}
        \lambda_{\phi}((x, \theta), \nu) &\defeq x^{(1)} \measureIntegral{ \beta(\theta, \cdot)  }{\nu(\{e_2\}\times \cdot )}\left(\phi \left(\nu -\frac{1}{n} \delta_{(x, \theta)} + \frac{1}{n}\delta_{(x - e_1 + e_2, \theta)} \right) 
             - \phi \left(\nu \right) \right) \\
             &{}\quad 
             + \gamma(\theta)x^{(2)} \left( \phi \left(\nu -\frac{1}{n} \delta_{(x, \theta)} + \frac{1}{n}\delta_{(x - e_2 + e_3, \theta)} \right) 
             - \phi \left(\nu \right) \right)\eqcomma
    \end{aligned}
             \label{eq:defn_lambda_phi}
\end{align}
for each bounded, measurable function $\phi : \spaceOfOccupationMeasures{1}{\Delta_3 \times \Theta} \to \setOfReals$, where the first integral $\measureIntegral{\beta(\theta, \cdot)  }{\nu (\{e_2\}\times \cdot )}$ is defined as follows:
\begin{align*}
    \measureIntegral{\beta(\theta, \cdot)  }{\nu (\{e_2\}\times \cdot )} = \frac{1}{n} \sum_{j=1}^{n} \beta(\theta, \theta_j) \indicator{\{e_2\}}{x_j} =  \frac{1}{n} \sum_{j=1}^{n} \beta(\theta, \theta_j) x_j^{(2)}\eqcomma 
\end{align*}
for $\nu = \frac{1}{n} \sum_{j=1}^{n} \delta_{(x_j, \theta_j)} \in \spaceOfOccupationMeasures{1}{\Delta_3\times \Theta}$. When $\nu$ is not atomic, the integral should be interpreted as 
\begin{align}
    \measureIntegral{\beta(\theta, \cdot)  }{\nu (\{e_2\}\times \cdot )} = \int_{\Theta} \beta(\theta, p)\nu(\{e_2\}\times \differential{p})\eqstop 
    \label{eq:integral_interpretation}
\end{align}

The next lemma shows that the measure-valued process $\{\nmu_t : t\ge 0\}$ is indeed a Markov process, and describes its generator.


\begin{myLemma}
    Assume $\{\nX_1(0), \nX_2(0), \ldots, \nX_n(0)\}$ is an exchangeable collection of random variables. Furthermore, assume $\theta_1, \theta_2, \ldots, \theta_n$ are \ac{iid}. 
    Then, the stochastic process $\{\nmu_t: t\ge 0\}$ is an $\spaceOfOccupationMeasures{1, \text{point}}{\Delta_3 \times \Theta}$-valued continuous-time jump Markov process with (infinitesimal) generator $\myOperator{L}_n$ acting on a bounded measurable real-valued function $\phi$ of the measure $\nu = \frac{1}{n} \sum_{j=1}^{n}\delta_{(x_j, y_j)} \in \spaceOfOccupationMeasures{1, \text{point}}{\Delta_3 \times \Theta}$ as follows
    \begin{align}
        \myOperator{L}_n \phi \left(\nu\right) &= n \measureIntegral{ \lambda_{\phi}(\cdot, \nu) }{\nu}\eqcomma 
        \label{eq:generator_of_measure_valued_process}
    \end{align}
    where the function $\lambda_{\phi}$ is defined in \eqref{eq:defn_lambda_phi}. 
\label{lemma:generator_of_measure_valued_process}
\end{myLemma} 
\begin{proof}[Proof of \Cref{lemma:generator_of_measure_valued_process}]
    It is straightforward that the process $\{ \nmu_t : t\ge 0\}$ is a Markov process. For example, see \citet[Proposition 2.3.3]{Dawson1993Measure_valued}. To calculate its (infinitesimal) generator \citep{Ethier:1986:MPC,Bobrowski2020Generators}, note that the process $\{\nmu_t : t\ge 0\}$ jumps whenever there is a jump in the system of \acp{SDE} in \eqref{eq:combined_sde_for_i_th_individual}. The trajectories of the process $\{\nmu_t : t\ge 0\}$ can be described as solutions to the following \ac{SDE} (written in the integral form) 
    \begin{align*}
        \nmu_t (\cdot) &= \nmu_0 (\cdot) + \sum_{i=1}^{n} 
        \int_{(0, t] \times \setOfPositiveReals} \left(-\frac{1}{n}\delta_{(\nX_i(s-), \theta_i)}(\cdot)+ \frac{1}{n}\delta_{(\nX_i(s-)-e_1+e_2, \theta_i)}(\cdot)\right) \\
        &{} \quad \quad \quad \quad \times 
         \indicator{[0, \nS_{i}(s-) n^{-1}\sum_{j=1}^{n}\beta(\theta_i, \theta_j)  \nI_j(s-) ]}{v}\PRM_{i}^{(1)}(\differential{s}\times\differential{v})\\ 
        &{}\quad \quad  + \sum_{i=1}^{n} \int_{(0, t] \times \setOfPositiveReals} \left( -\frac{1}{n}\delta_{(\nX_i(s-), \theta_i)}(\cdot)+ \frac{1}{n}\delta_{(\nX_i(s-)-e_2+e_3, \theta_i)}(\cdot)\right) \\
        &\quad \quad \quad \times 
        \indicator{[0, \gamma(\theta_i) \nI_{i}(s-) ]}{v}\PRM_{i}^{(2)}(\differential{s}\times\differential{v})\eqcomma 
    \end{align*}
    with the initial data $\nmu_0(\cdot) = \frac{1}{n} \sum_{i=1}^{n} \delta_{(\nX_i(0), \theta_i)}(\cdot)$. By an application of the It\^o's formula for \acp{SDE} driven by \acp{PRM} (\eg, see \cite[Lemma~4.4.5]{Applebaum_2009Levy}, \cite[Theorem~5.1]{IkedaWatanabe2014Stochastic}) for a (real-valued)   bounded, measurable function $\phi$ on the space $\spaceOfOccupationMeasures{1, \text{point}}{\Delta_3 \times \Theta}$, we have 
    \begin{align}
        \phi(\nmu_t)  &= \phi(\nmu_0)+  \sum_{i=1}^{n} 
        \int_{(0, t] \times \setOfPositiveReals} \left(\phi\left(\nmu_{s-} -\frac{1}{n}\delta_{(\nX_i(s-), \theta_i)}+ \frac{1}{n}\delta_{(\nX_i(s-)-e_1+e_2, \theta_i)}\right) \right. \nonumber\\
        &\quad \quad \quad \left. - \phi\left(\nmu_{s-}\right) \right) 
         \indicator{[0, \nS_{i}(s-) n^{-1}\sum_{j=1}^{n} \beta(\theta_i, \theta_j) \nI_j(s-) ]}{v}\PRM_{i}^{(1)}(\differential{s}\times\differential{v}) \nonumber \\ 
        &{}\quad \quad  + \sum_{i=1}^{n} \int_{(0, t] \times \setOfPositiveReals} \left(\phi\left(\nmu_{s-} -\frac{1}{n}\delta_{(\nX_i(s-), \theta_i)}+ \frac{1}{n}\delta_{(\nX_i(s-)-e_2+e_3, \theta_i)}\right) 
        - \phi\left(\nmu_{s-}\right)\right) \nonumber \\
         &\quad \quad \quad \times 
        \indicator{[0, \gamma(\theta_i) \nI_{i}(s-) ]}{v}\PRM_{i}^{(2)}(\differential{s}\times\differential{v}) \nonumber  \\
        &= \phi(\nmu_0) + \int_{(0, t]} \sum_{i=1}^{n} \left( \frac{1}{n} \sum_{j=1}^{n} \beta(\theta_i, \theta_j) \nS_{i}(s) \nI_j(s) \right) \nonumber \\
        &{}\quad \quad \times   \left(\phi\left(\nmu_{s} -\frac{1}{n}\delta_{(\nX_i(s), \theta_i)}+ \frac{1}{n}\delta_{(\nX_i(s)-e_1+e_2, \theta_i)}\right) 
        - \phi\left(\nmu_{s}\right) \right)\differential{s} \nonumber  \\
        &{}\quad \quad 
        + \int_{(0, t]} \sum_{i=1}^{n} \gamma(\theta_i) \nI_{i}(s) \left(\phi\left(\nmu_{s} -\frac{1}{n}\delta_{(\nX_i(s), \theta_i)}+ \frac{1}{n}\delta_{(\nX_i(s)-e_2+e_3, \theta_i)}\right) \right. \nonumber \\
        &{}\quad \quad \quad \quad  \left.
        - \phi\left(\nmu_{s}\right)\right) \differential{s} + \nM_{\phi}(t) \nonumber 
        \\
        & = \phi(\nmu_0) +  \int_{(0, t]} \myOperator{L}_n \phi (\nmu_s) \differential{s} + \nM_{\phi}(t)\eqcomma \label{eq:ito_martingale_mu}
    \end{align}
    where the operator $\myOperator{L}_n$ is given by 
    \begin{align}
        \myOperator{L}_n \phi \left(\nu\right) &= \sum_{i=1}^{n} \left( \frac{1}{n} \sum_{j=1}^{n} \beta(y_i, y_j) x_i^{(1)} x_j^{(2)}  \right) 
        \left( \phi \left(\nu -\frac{1}{n} \delta_{(x_i, y_i)} + \frac{1}{n}\delta_{(x_i - e_1 + e_2, y_i)} \right) 
             - \phi \left(\nu \right) \right) \nonumber \\
        &{}\quad + \sum_{i=1}^{n} \gamma(y_i)x_i^{(2)} 
        \left( \phi \left(\nu -\frac{1}{n} \delta_{(x_i, y_i)} + \frac{1}{n}\delta_{(x_i - e_2 + e_3, y_i)} \right) 
             - \phi \left(\nu \right) \right) \nonumber \\
        &= n \measureIntegral{\lambda_{\phi}(\cdot, \nu)}{\nu}\eqcomma \nonumber 
    \end{align}
    and the stochastic process $\{\nM_{\phi}(t) : t\ge 0\}$ is a square-integrable zero-mean $\history{t}$-martingale given by 
    \begin{align*}
        \nM_{\phi}(t) &=  \sum_{i=1}^{n} 
        \int_{(0, t] \times \setOfPositiveReals} \left(\phi\left(\nmu_{s-} -\frac{1}{n}\delta_{(\nX_i(s-), \theta_i)}+ \frac{1}{n}\delta_{(\nX_i(s-)-e_1+e_2, \theta_i)}\right) \right.\\
        &\quad \quad \quad \left. - \phi\left(\nmu_{s-}\right) \right) 
         \indicator{[0, \nS_{i}(s-) n^{-1} \sum_{j=1}^{n} \beta(\theta_i, \theta_j)  \nI_j(s-) ]}{v} \compPRM_{i}^{(1)}(\differential{s}\times\differential{v})\\ 
        &{}\quad \quad  + \sum_{i=1}^{n} \int_{(0, t] \times \setOfPositiveReals} \left(\phi\left(\nmu_{s-} -\frac{1}{n}\delta_{(\nX_i(s-), \theta_i)}+ \frac{1}{n}\delta_{(\nX_i(s-)-e_2+e_3, \theta_i)}\right) 
        - \phi\left(\nmu_{s-}\right)\right) \\
         &\quad \quad \quad \times 
        \indicator{[0, \gamma(\theta_i) \nI_{i}(s-) ]}{v}\compPRM_{i}^{(2)}(\differential{s}\times\differential{v}) \eqcomma 
    \end{align*}
    with $\compPRM_{i}^{(1)}$ and $\compPRM_{i}^{(2)}$ denoting the compensated \acp{PRM} associated with $\PRM_{i}^{(1)}$ and $\PRM_{i}^{(2)}$ respectively, for $i\in \{1, \ldots, n\}$. Taking expectation on both sides of \eqref{eq:ito_martingale_mu}, and then taking derivative with respect to $t$ and evaluating at $t=0$ yields the generator $\myOperator{L}_n $ in \eqref{eq:generator_of_measure_valued_process} of the $\spaceOfOccupationMeasures{1, \text{point}}{\Delta_3 \times \Theta}$-valued Markov process $\{\nmu_t : t\ge 0\}$.

\end{proof}

In order to prove a weak convergence of the sequence $\{\nmu : n\ge 1\}$, we will use a tightness-uniqueness argument. Following this strategy, we will first show that the collection $\{\nmu : n\ge 1\}$ is relatively compact as $D([0, T], \spaceOfOccupationMeasures{1}{\Delta_3 \times \Theta})$-valued random variables. To be more precise, we show the sequence of probability laws of  the collection $\{\nmu : n\ge 1\}$ is relatively compact in the space of probability measures on $D([0, T], \spaceOfOccupationMeasures{1}{\Delta_3 \times \Theta})$ equipped with the topology of weak convergence. Once relative compactness is established, we can extract a weakly convergent subsequence. We will then show the limit points of the convergent subsequences are unique. Therefore, (weak) convergence holds along the entire sequence. 
 
The standard introductory textbook on this line of arguments is \cite{Billingsley1999Convergence}, which provides necessary tools to deal with real-valued stochastic processes. For more general stochastic processes that take values in metric spaces, or Banach spaces, we refer to \cite[Chapter 3]{Ethier:1986:MPC}. Other excellent resources are \cite{jacod2003limit,Whitt2002StochLimits,kallenberg2021foundations}. 

\subsection{Tightness of the empirical measure}
\label{sec:flln_tightness}

Next, we show that the sequence of stochastic processes $\{\nmu : n\ge 1\}$ is relatively compact as $D([0, T], \spaceOfOccupationMeasures{1}{\Delta_3 \times \Theta})$-valued random variables. 
To that end, let us first define some useful notations. Define the function $\Lambda_f : \left(\Delta_3\times \Theta \right)\times \spaceOfOccupationMeasures{1}{\Delta_3\times \Theta} \to \setOfReals$ as follows
\begin{align}
    \begin{aligned}
        \Lambda_{f}((x, \theta), \nu) &\defeq x^{(1)} \measureIntegral{ \beta(\theta, \cdot)  }{\nu(\{e_2\}\times \cdot )}\left(f \left({x - e_1 + e_2, \theta} \right) 
             - f \left({x, \theta} \right) \right) \\
             &{}\quad 
             + \gamma(\theta)x^{(2)} \left( f \left({x - e_2 + e_3, \theta} \right) 
             - f \left(x, \theta \right) \right)\eqcomma
    \end{aligned}
             \label{eq:defn_Lambda_f}
\end{align}
where the first integral $\measureIntegral{\beta(\theta, \cdot)   }{\nu ( \{e_2\} \times \cdot )}$ should be interpreted as \eqref{eq:integral_interpretation}.

\begin{myLemma}
    Assume $\{\nX_1(0), \nX_2(0), \ldots, \nX_n(0)\}$ is an \ac{iid} collection of $\Delta_3$-valued random variables. Furthermore, assume $\theta_1, \theta_2, \ldots, \theta_n$ are \ac{iid} random variables taking values in a compact subset $\Theta$ of the Euclidean space. The sequence of probability laws of  $\{\nmu : n\ge 1\}$ is tight in $\spaceOfOccupationMeasures{1}{ D([0, T], \spaceOfOccupationMeasures{1}{\Delta_3 \times \Theta})}$, the space of probability measures on $D([0, T], \spaceOfOccupationMeasures{1}{\Delta_3 \times \Theta})$, equipped with the topology of weak convergence.

    \label{lemma:tightness-empirical_measure}
\end{myLemma}
\begin{proof}[Proof of \Cref{lemma:tightness-empirical_measure}]
    Let $f$ be an element of $C(\Delta_3\times \Theta, \setOfReals)$. Our aim is to first show that the sequence of real-valued stochastic processes $\{\measureIntegral{f}{\nmu}: n \ge 1\}$ is tight in $D([0, \infty), \setOfReals)$. Verification of the compact containment condition (see \Cref{defn:compact_containment_condition} in Appendix~\ref{sec:limits_of_stoch_pro}) is straightforward since 
    \begin{align*}
        \measureIntegral{f}{\nmu_t} = \frac{1}{n} \sum_{i=1}^{n} f(\nX_i(t), \theta_i)
    \end{align*}
    and $f$ is bounded (by virtue of being defined on a compact domain). Indeed, for every $\varepsilon>0$, we can find a positive real number $K_\varepsilon$ (\eg, any number strictly larger than $\norm{f}_{\infty}$) such that 
    \begin{align}
        \inf_{n} \probOf{\absolute{\measureIntegral{f}{\nmu_t} } < K_\varepsilon \, \forall t \in [0, T] } \ge  1- \varepsilon\eqstop
        \label{eq:compact_containment_f_mu} 
    \end{align}
    
    To verify the modulus of continuity condition (see Appendix~\ref{sec:limits_of_stoch_pro} for background), note that since $\Delta_3$ is a finite set, and $\Theta$ is assumed to be compact, the map $\phi : \spaceOfOccupationMeasures{1}{\Delta_3\times \Theta } \to \setOfReals$ defined by 
    \begin{align*}
        \phi (\nu) \defeq \measureIntegral{f}{\nu}  
    \end{align*}
    is a bounded, measurable map. Therefore, from \eqref{eq:ito_martingale_mu} in the proof of \Cref{lemma:generator_of_measure_valued_process}, the stochastic process  
    \begin{align*}
        \nM_{\phi}(t) &= \phi(\nmu_t) -\phi(\nmu_0) - \int_{0}^{t}\myOperator{L}_n \phi(\nmu_s) \differential{s}
    \end{align*}
    is a zero-mean square integrable $\history{t}$-martingale. Now, see that 
    \begin{align*}
        \int_{0}^{t}\myOperator{L}_n \phi(\nmu_s) \differential{s} &= \int_{0}^{t} n \measureIntegral{\lambda_{\phi}(\cdot, \nmu_s)}{\nmu_s} \differential{s}\\
        & = \int_{0}^{t} n \int_{\Theta}\left(\phi(\nmu_s - \frac{1}{n}\delta_{(e_1, \theta) } + \frac{1}{n} \delta_{(e_2, \theta)}) - \phi(\nmu_s)  \right)\int_{\Theta} \beta(\theta, \theta') \\
        &{}\quad \quad \quad 
        \times \nmu_s(\{e_2\}\times \differential{\theta'}) \nmu_s(\{e_1\}\times \differential{\theta})\differential{s}\\
        &{}\quad + \int_{0}^{t} n \int_{\Theta}\left(\phi(\nmu_s - \frac{1}{n}\delta_{(e_2, \theta) } + \frac{1}{n} \delta_{(e_3, \theta)}) - \phi(\nmu_s)  \right)\gamma(\theta) \\
        &{}\quad \quad \quad 
        \times \nmu_s(\{e_2\}\times \differential{\theta}) \differential{s}\\
        & = \int_{0}^{t} n \int_{\Theta}\left(\measureIntegral{f}{\nmu_s - \frac{1}{n}\delta_{(e_1, \theta) } + \frac{1}{n} \delta_{(e_2, \theta)}} -\measureIntegral{f}{\nmu_s} \right)\int_{\Theta} \beta(\theta, \theta') \\
        &{}\quad \quad \quad 
        \times \nmu_s(\{e_2\}\times \differential{\theta'}) \nmu_s(\{e_1\}\times \differential{\theta})\differential{s}\\
        &{}\quad + \int_{0}^{t} n \int_{\Theta}\left(\measureIntegral{f}{\nmu_s - \frac{1}{n}\delta_{(e_2, \theta) } + \frac{1}{n} \delta_{(e_3, \theta)}} -\measureIntegral{f}{\nmu_s} \right)\gamma(\theta) \\
        &{}\quad \quad \quad 
        \times \nmu_s(\{e_2\}\times \differential{\theta}) \differential{s}\\
        & = \int_{0}^{t} n \int_{\Theta}\frac{1}{n}\left(f(e_2, \theta) - f(e_1, \theta) \right)\int_{\Theta} \beta(\theta, \theta') 
        \nmu_s(\{e_2\}\times \differential{\theta'}) \nmu_s(\{e_1\}\times \differential{\theta})\differential{s}\\
        &{}\quad + \int_{0}^{t} n \int_{\Theta} \frac{1}{n} \left(f(e_3, \theta) - f(e_2, \theta) \right)\gamma(\theta) 
        \nmu_s(\{e_2\}\times \differential{\theta}) \differential{s}\\
        &= \int_{0}^{t} \measureIntegral{\Lambda_f(\cdot, \nmu_s) }{\nmu_s} \differential{s}\eqcomma 
    \end{align*}
    where the map $\Lambda_f$ is defined in \eqref{eq:defn_Lambda_f}. Therefore, for each $f \in C(\Delta_3\times \Theta, \setOfReals)$, the stochastic process 
    \begin{align*}
        \nM_{\phi}(t) &= \measureIntegral{f}{\nmu_t} -\measureIntegral{f}{\nmu_0} - \int_{0}^{t} \measureIntegral{\Lambda_f(\cdot, \nmu_s) }{\nmu_s} \differential{s}
    \end{align*}
    is a zero-mean square integrable $\history{t}$-martingale, with $\phi(\nu) = \measureIntegral{f}{\nu}$. The predictable quadratic variation of the martingale $\nM_{\phi}$ is given by 
    \begin{align*}
        \predictableVariation{\nM_{\phi}}(t) &=  \sum_{i=1}^{n} \sum_{j=1}^{n}
        \int_{0}^{t} \left(\phi\left(\nmu_{s} -\frac{1}{n}\delta_{(\nX_i(s), \theta_i)}+ \frac{1}{n}\delta_{(\nX_i(s)-e_1+e_2, \theta_i)}\right) 
        - \phi\left(\nmu_{s}\right) \right)^2 \\
        &{} \quad \quad \quad \quad \times 
         n^{-1}\beta(\theta_i, \theta_j) \nS_{i}(s) \nI_j(s) \differential{s}\\ 
        &{}\quad \quad  + \sum_{i=1}^{n} \int_{0}^{t} \left(\phi\left(\nmu_{s} -\frac{1}{n}\delta_{(\nX_i(s-), \theta_i)}+ \frac{1}{n}\delta_{(\nX_i(s)-e_2+e_3, \theta_i)}\right) 
        - \phi\left(\nmu_{s}\right)\right)^2 \\
         &\quad \quad \quad \quad \times 
        \gamma(\theta_i) \nI_{i}(s)\differential{s} \\
        & = \int_{0}^{t} \Gamma_n(\phi, \phi)(\nmu_s)\differential{s}\eqcomma 
    \end{align*}
    where the operator $\Gamma_n$ is known as the ``carr\'e du champ'' operator \citep{Bakry2014AnalysisGeometry,Ledoux2000Geometry}, and is given by 
    \begin{align*}
        \Gamma_n(\phi, \psi) \defeq \myOperator{L}_n (\phi \psi) - \phi \myOperator{L}_n (\psi) -  \psi \myOperator{L}_n (\phi)
    \end{align*}
    for bounded, measurable functions $\phi, \psi : \spaceOfOccupationMeasures{1}{\Delta_3\times \Theta} \to \setOfReals$. The operator $\Gamma_n$
    measures how far the operator $\myOperator{L}_n$ is from being a derivation (see \Cref{defn:derivation} in Appendix \ref{sec:additional}). With our choice $\phi(\nu) = \measureIntegral{f}{\nu}$, we further simplify the  expression for the quadratic variation to get
    \begin{align*}
        \predictableVariation{\nM_{\phi}}(t) &= n \int_{0}^{t} \int_{\Theta} \left(\phi\left(\nmu_{s} -\frac{1}{n}\delta_{(e_1, \theta)}+ \frac{1}{n}\delta_{(e_2, \theta)}\right) 
        - \phi\left(\nmu_{s}\right) \right)^2 \\
        &{} \quad\quad\quad   \times \int_{\Theta} \beta(\theta, \theta') \nmu_s(\{e_2\}\times \differential{\theta'}) \nmu_s(\{e_1\}\times \differential{\theta}) \differential{s} \\
        &{}\quad + n \int_{0}^{t} \int_{\Theta} \left(\phi\left(\nmu_{s} -\frac{1}{n}\delta_{(e_2, \theta)}+ \frac{1}{n}\delta_{(e_3, \theta)}\right) 
        - \phi\left(\nmu_{s}\right)\right)^2 \gamma(\theta) \nmu_s(\{e_2\}\times \differential{\theta}) \differential{s}\\
        & = \frac{1}{n} \int_{0}^{t} \int_{\Theta} \left(f(e_2, \theta) - f(e_1, \theta)\right)^2 
        \int_{\Theta} \beta(\theta, \theta') \nmu_s(\{e_2\}\times \differential{\theta'}) \nmu_s(\{e_1\}\times \differential{\theta}) \differential{s} \\
        &{}\quad + \frac{1}{n} \int_{0}^{t} \int_{\Theta} \left(f(e_3, \theta) - f(e_2, \theta)\right)^2 \gamma(\theta) \nmu_s(\{e_2\}\times \differential{\theta}) \differential{s}\\
        &{}\leq \frac{C(t)\norm{f}_{\infty}^2 }{n} \to 0\eqcomma 
    \end{align*}     
    both in probability and in $L^2(\Omega, \history{}, \prob)$
    as $n\to \infty$, 
    for some positive constant $C(t)$ since $\beta(u, v) \le \hat{\beta} $, and $ \gamma(u) \le \hat{\gamma} $ for some positive constants $\hat{\beta}$ and $\hat{\gamma}$, and the random variables $\theta_1, \theta_2, \ldots, \theta_n$ are assumed to be \ac{iid} taking values in a compact Euclidean space. Therefore,
    \begin{align}
        \sup_{t \le T} \predictableVariation{\nM_{\phi}}(t) \to 0\eqcomma 
        \label{eq:conv_of_quad_var}
    \end{align}
    in probability and in $L^2(\Omega, \history{}, \prob)$
    as $n \to \infty$. Now, for any positive real numbers $\varepsilon, \eta$, the  Lenglart--Rebolledo inequality (see \cite[Lemma 3.7]{Whitt2007MCLT}, and \cite[Remark 4.17]{karatzas1991brownian}) gives us 
        \begin{align}
            \probOf{\sup_{t\leq T} \absolute{\nM_{\phi}(t)} >\eta}\leq \varepsilon +\probOf{ \predictableVariation{\nM_{\phi} }(t) > \varepsilon\eta^2}\eqstop 
            \label{eq:Lenglart-Rebolledo}
        \end{align}
    By virtue of the convergence of the quadratic variation in \eqref{eq:conv_of_quad_var}, letting $n\to \infty$ in \eqref{eq:Lenglart-Rebolledo} yields 
    \begin{align}
        \sup_{t\leq T} \absolute{\nM_{\phi}(t)} \ConvInProb 0 \eqcomma 
        \label{eq:martingale_phi_conv_zero}
    \end{align}
    as $n\to \infty$. Therefore, for $0\le s < t$, 
    \begin{align*}
        \absolute{\measureIntegral{f}{\nmu_t} - \measureIntegral{f}{\nmu_s} }&= \absolute{ (\nM_{\phi}(t) - \nM_{\phi}(s) ) + \int_s^t \measureIntegral{\Lambda_f(\cdot, \nmu_u) }{\nmu_u} \differential{u} }\\
        &{} \le \absolute{\nM_{\phi}(t) - \nM_{\phi}(s) } + C(t, \norm{f}_{\infty}) \absolute{t-s} \eqcomma
    \end{align*}
    for some constant $C(t, \norm{f}_{\infty})$. Therefore, by virtue of the convergence in \eqref{eq:conv_of_quad_var} and \eqref{eq:martingale_phi_conv_zero},  for any $\varepsilon >0$, and $\eta>0$ there exists a $0<\delta<1$ and $n_0$ such that 
    \begin{align*}
        \sup_{n \ge n_0} \probOf{\modulusOfcontinuity(\measureIntegral{f}{\nmu}, T, \delta) > \eta} \le \varepsilon \eqcomma 
    \end{align*}
    where the modulus of continuity $\modulusOfcontinuity$ is defined in  \eqref{eq:modu_cont_defn} of \Cref{defn:modulus_continuity_C} in Appendix \ref{sec:limits_of_stoch_pro}. 
    The above along with the compact containment condition in \eqref{eq:compact_containment_f_mu}, in fact,  shows that the sequence $\{\measureIntegral{f}{\nmu}: n\ge 1\}$ is $C$-tight in $D([0, T], \setOfReals)$ (see \Cref{def:c-tight-in-d} in Appendix \ref{sec:limits_of_stoch_pro}).


    Since the sequence of real-valued stochastic processes $\{\measureIntegral{f}{\nmu}: n \ge 1\}$ is tight in $D([0, T], \setOfReals)$, for each $f \in C(\Delta_3\times \Theta)$, we can apply \citet[Theorem 3.7.1]{Dawson1993Measure_valued} (restated as \Cref{thm:tightness_criterion_measure_valued} in Appendix \ref{sec:limits_of_stoch_pro} for the sake of completeness) to conclude that the sequence of probability laws of $\{\nmu : n \ge 1\}$ is tight in $\spaceOfOccupationMeasures{1}{ D([0, T], \spaceOfOccupationMeasures{1}{\Delta_3 \times \Theta})}$, the space of probability measures on $D([0, T], \spaceOfOccupationMeasures{1}{\Delta_3 \times \Theta})$, equipped with the topology of weak convergence. This completes the proof. 
\end{proof}

\begin{myRemark}
    \label{rem:compactness}
    The assumption that the set $\Theta$ is a compact subset of a Euclidean space can be removed in several ways. For example, \citet[Theorem 3.7.1 (b)]{Dawson1993Measure_valued} can be applied to show relative compactness of the collection $\{\nmu : n\ge 1\}$ in the space of \cadlag functions with values in the space of tempered measures on $\setOfReals^d$ \citep[Section 3.1.5]{Dawson1993Measure_valued}. The other approach would be to first show relative compactness in the vague topology by establishing the relative compactness of evaluations of continuous functions with compact support as real-valued stochastic processes and then upgrade it to the topology of weak convergence. This is precisely the strategy adopted in \cite{Fournier2004microscopic}. In order to keep our technical arguments as simple as possible, we do not pursue this line of argument here. The assumption of a compact covariate space $\Theta$ is sufficient for most practical purposes.  
\end{myRemark}

By virtue of \Cref{lemma:tightness-empirical_measure}, we can extract a convergent subsequence $\{\mu^{(n_k)} : k\ge 0\}$. If the limit points of such convergent subsequences are unique, we can conclude the whole sequence $\{\nmu : n\ge 1\}$ converges to this limit point. To this end, let us now proceed to identify the limit point(s) and verify their uniqueness.

\subsection{Uniqueness of solution to a measure-valued differential equation}
\label{sec:flln_uniqueness}

Consider the following differential equation written in the weak form: for all bounded, continuous functions $f: \Delta_3\times \Theta \to \setOfReals$, 
\begin{align}
    \measureIntegral{f}{\mu_t} = \measureIntegral{f}{\mu_0} + \int_{0}^{t} \measureIntegral{\Lambda_{f}(\cdot, \mu_s) }{\mu_s} \differential{s}\eqcomma 
    \label{eq:limiting_ode}
\end{align}
where the function $\Lambda_f : \left(\Delta_3\times \Theta \right)\times \spaceOfOccupationMeasures{1}{\Delta_3\times \Theta} \to \setOfReals$ is defined in \eqref{eq:defn_Lambda_f}. It is not difficult to see that there is a unique solution $g_f$ defined as 
\begin{align*}
    g_f(t) \defeq \measureIntegral{f}{\mu_t}\eqcomma \text{ for } t\ge 0\eqcomma 
\end{align*}
satisfying the integral equation \eqref{eq:limiting_ode}, for each fixed $f\in C_b( \Delta_3\times \Theta , \setOfReals)$. We want to prove that there is a unique (probability) measure-valued solution to \eqref{eq:limiting_ode}.

\begin{myLemma}
    Let $\Theta \subseteq \setOfReals^d$, for some $d\ge 1$. 
    There is a unique measure-valued solution to the differential equation \eqref{eq:limiting_ode} satisfying $\measureIntegral{1}{\mu_t} =1 $ for all $t\le T < \infty$. 
    \label{lemma:flln_uniqueness}
\end{myLemma}
\begin{proof}[Proof of \Cref{lemma:flln_uniqueness}]
    Let $f \in \mathbb{F}_{\mathsf{TV}} \defeq \{ g: \Delta_3\times \Theta \to \setOfReals : \norm{g}_{\infty} \le 1\} $. Let us assume $\mu$ and $\nu$ are two different (measure-valued) solutions to \eqref{eq:limiting_ode} with $\mu_0= \nu_0$ (equality of measures is understood in the sense that integrals of all bounded, continuous functions with respect to the two measures are equal, \ie, $\measureIntegral{g}{\mu_0} = \measureIntegral{g}{\nu_0}$ for all bounded, continuous functions $g$), and $\measureIntegral{1}{\mu_t} = \measureIntegral{1}{\nu_t} = 1 $ for all $t\le T$. Then, we have 
    \begin{align*}
        \measureIntegral{f}{\mu_t - \nu_t} & =  \int_0^t \left(\measureIntegral{\Lambda_{f}(\cdot, \mu_s) }{\mu_s} - \measureIntegral{\Lambda_{f}(\cdot, \nu_s) }{\nu_s} \right)\differential{s} + \measureIntegral{f}{\mu_0 - \nu_0} \\
        & = \int_0^t \left( \int_{\Theta} \left(f(e_2, \theta)-f(e_1, \theta)  \right)  \int_{\Theta} \beta(\theta, \theta') \mu_s(\{e_2\}\times  \differential{\theta'}) \mu_s(\{e_1\} \times \differential{\theta})  \right) \differential{s} \\
        &{}\quad - \int_0^t \left( \int_{\Theta} \left(f(e_2, \theta)-f(e_1, \theta)  \right)  \int_{\Theta} \beta(\theta, \theta') \nu_s(\{e_2\}\times  \differential{\theta'}) \nu_s(\{e_1\} \times \differential{\theta})  \right) \differential{s}\\
        & \quad + \int_0^t \int_{\Theta}\gamma(\theta) \left(f(e_3, \theta) - f(e_2, \theta)\right) \mu_s(\{e_2\}\times \differential{\theta}) \differential{s}\\
        &{}\quad  
        - \int_0^t \int_{\Theta}\gamma(\theta) \left(f(e_3, \theta) - f(e_2, \theta)\right) \nu_s(\{e_2\}\times \differential{\theta}) \differential{s} + \measureIntegral{f}{\mu_0 - \nu_0} \eqstop 
    \end{align*}

    Since $ \absolute{\beta(u, v)} \le \hat{\beta} \eqcomma \text{ and } \absolute{\gamma(u)} \le \hat{\gamma} \eqcomma $ for some positive constants $\hat{\beta}$ and $\hat{\gamma}$, we have 
    \begin{align}
        \absolute{\measureIntegral{f}{\mu_t - \nu_t}} & \le  \absolute{\measureIntegral{f}{\mu_0 - \nu_0}} + \left( 2 \hat{\beta} + \hat{\gamma}  \right) 2 \norm{f}_{\infty} \int_0^t \TotalVariation{\mu_s , \nu_s} \differential{s} \eqstop 
        \label{eq:measure_eq_pre_Gronwall}
    \end{align}
    Taking supremum over all functions $f \in \mathbb{F}_{\mathsf{TV}}$, and using Gr\"onwall's inequality gives us 
    \begin{align*}
        \sup_{t\le T}\TotalVariation{\mu_t, \nu_t} &\le  \TotalVariation{\mu_0, \nu_0} \myExp{2 \left( 2 \hat{\beta} + \hat{\gamma}  \right) T}\eqstop 
    \end{align*} 
    Since $\TotalVariation{\mu_0, \nu_0} = 0$ by assumption, this gives us the desired uniqueness of (measure-valued) solution to \eqref{eq:limiting_ode} admitting $\measureIntegral{1}{\mu_t} = 1 $ for all $t\le T$. 
\end{proof}

We are now ready to present our \ac{FLLN} for the empirical random measures. 

\subsection{\acl{FLLN}}
\begin{myTheorem}
\label{thm:flln}
Assume $\{\nX_1(0), \nX_2(0), \ldots, \nX_n(0)\}$ is an \ac{iid} collection of $\Delta_3$-valued random variables. Furthermore, assume $\theta_1, \theta_2, \ldots, \theta_n$ are \ac{iid} random variables taking values in $\Theta$, a compact subset of the Euclidean space. Additionally we assume that the rate functions $\beta: \Theta\times \Theta \to \setOfPositiveReals$ and $\gamma: \Theta \to \setOfPositiveReals$ are continuous. Then, 
\begin{align*}
    \nmu \ConvInDist \mu
\end{align*}
in $D([0, T], \spaceOfOccupationMeasures{1}{\Delta_3\times \Theta})$
as $n\to \infty$ for some deterministic continuous (probability) measure-valued function $\mu \in C([0, T], \spaceOfOccupationMeasures{1}{\Delta_3\times\Theta})$ satisfying the differential equation: 
for all bounded, continuous functions $f: \Delta_3\times \Theta \to \setOfReals$, 
\begin{align*}
    \measureIntegral{f}{\mu_t} = \measureIntegral{f}{\mu_0} + \int_{0}^{t} \measureIntegral{\Lambda_{f}(\cdot, \mu_s) }{\mu_s} \differential{s}\eqcomma 
\end{align*}
where the function $\Lambda_f : \left(\Delta_3\times \Theta \right)\times \spaceOfOccupationMeasures{1}{\Delta_3\times \Theta} \to \setOfReals$ is defined in \eqref{eq:defn_Lambda_f}. 

\end{myTheorem}
\begin{proof}[Proof of \Cref{thm:flln}]
    Under our assumption about the initial data, note that 
    \begin{align*}
    \probOf{\{\omega \in \Omega : \nmu_0(\omega) \ConvInDist \mu_0(\omega)\} } = 1
    \end{align*}
    for some probability measure $\mu_0$ on $\Delta_3\times \Theta$ by the standard  \ac{SLLN} (\eg, see \citet[Theorem 3]{Varadarajan1958LLN}, and also \citet{parthasarathy2005probability}). Next, by virtue of \Cref{lemma:tightness-empirical_measure}, we know the sequence $\{\nmu : n\ge 1\}$ is relatively compact. Therefore, we 
    can extract a convergent subsequence of $\{\nmu : n\ge 1\}$, which we continue to denote as $\{\nmu : n\ge 1\}$. If all convergent subsequences have the same limit point, we can conclude that the convergence holds along the entire sequence. To this end, we first proceed to identify the limit point. 

    Note that, for each bounded, continuous $f : \Delta_3\times \Theta \to \setOfReals$, and for each $t>0$,
    \begin{align}
        \sup_{t \le T} \absolute{\measureIntegral{f}{\nmu_t} - \measureIntegral{f}{\nmu_{t-}}} &{} \le \frac{1}{n} \sum_{i=1}^{n} \absolute{f(\nX_i(t), \theta_i) - f(\nX_i(t-), \theta_i) } 
        \le \frac{2\norm{f}_{\infty}}{n} \eqcomma 
        \label{eq:f_mu_modulus_of_continuity}
    \end{align}
    since only one of the components $\nX_1, \nX_2, \ldots, \nX_n$ can jump at a time (since the points  of the independent \acp{PRM} are distinct almost surely). Therefore, 
    \begin{align}
        \lim_{n\to \infty}\Eof{ \sup_{t \le T} \absolute{\measureIntegral{f}{\nmu_t} - \measureIntegral{f}{\nmu_{t-}}}} =0\eqcomma \nonumber
    \end{align}
    which, in turn, implies the limit points of $\{ \measureIntegral{f}{\nmu} : n \ge 1\}$ are continuous functions by virtue of \citet[Theorem 13.4]{Billingsley1999Convergence}. Taking supremum over $f \in \mathbb{F}_{\mathsf{TV}}$ in \eqref{eq:f_mu_modulus_of_continuity} and then taking expectation and then, the limit as $n\to \infty$, we can conclude that $$\lim_{n\to\infty}\Eof{\TotalVariation{\nmu_t, \nmu_{t-}}} =0.$$ 
    Thus,  the limit points of $\{\nmu : n\ge 1\}$ also lie in $C([0, T], \spaceOfOccupationMeasures{1}{\Delta_3\times \Theta})$ almost surely. This proves the assertion about the continuity of the limit point $\mu$  in the statement of \Cref{thm:flln}. 

    Next, we note that for each $f \in C(\Delta_3\times \Theta, \setOfReals)$ with $\phi(\nu) = \measureIntegral{f}{\nu}$, the stochastic process 
    \begin{align*}
        \nM_{\phi}(t) &= \measureIntegral{f}{\nmu_t} -\measureIntegral{f}{\nmu_0} - \int_{0}^{t} \measureIntegral{\Lambda_f(\cdot, \nmu_s) }{\nmu_s} \differential{s}
    \end{align*}
    is a zero-mean square integrable $\history{t}$-martingale. From the proof of \Cref{lemma:tightness-empirical_measure}, we know that 
    \begin{align*}
        \sup_{t\le T} \absolute{\nM_{\phi}(t)} \to 0
    \end{align*}
    in probability 
    as $n\to \infty$. In fact, the above convergence holds in $L^2(\Omega, \history{}, \prob)$ as well. (See \Cref{lem:martingale_L2_f}.) Therefore, any limit point $\mu$ of a convergent subsequence of $\{\nmu : n\ge 1\}$ must satisfy 
    \begin{align*}
        \measureIntegral{f}{\mu_t} = \measureIntegral{f}{\mu_0} + \int_{0}^{t} \measureIntegral{\Lambda_f(\cdot, \mu_s) }{\mu_s} \differential{s}\eqcomma 
    \end{align*}
    almost surely. 
    Since the above integral equation admits a unique measure-valued solution satisfying $\measureIntegral{1}{\mu_t} = 1$ for all $t\le T$ by virtue of \Cref{lemma:flln_uniqueness}, the proof of \Cref{thm:flln} would be complete if we show  $\measureIntegral{1}{\nmu_t} = 1 $ for all $t\le T$, almost surely. However, this is straightforward. 
\end{proof}

\begin{myRemark}
    Since the limit point $\mu$ is a deterministic function, the convergence in \Cref{thm:flln} holds not only in the weak sense but also in probability. 
\end{myRemark}

\begin{myRemark}\label{rem:D_infinity_conv}
    If $f \in D([0, \infty), E)$ with some metric space $E$, then the restriction of $f$ to $[0, T]$ is an element of $D([0, T], E)$. In that sense, the function $f$ itself can be thought of as elements of $D([0, T], E)$. Therefore, by virtue of \Cref{lemma:M_infinity_conv}, the convergence in \Cref{thm:flln} also holds in $D([0, \infty), \spaceOfOccupationMeasures{1}{\Delta_3\times \Theta})$. See also the discussion in Section 16 of \cite{Billingsley1999Convergence} (in particular, Theorems 16.2 and 16.7).
\end{myRemark}

There is a huge amount of literature on the topic of convergence of empirical random measures for \acp{IPS}, which are described in terms of an interacting system of \acp{SDE} driven by Wiener processes.  The work of \citet{Oelschlager1984LLN} is particularly relevant for us since it considers the case of jump processes as well. It appears that the framework of \citet{Oelschlager1984LLN} could be adapted to our setup. However, we expect that the verification of the provided sufficient conditions in \citet{Oelschlager1984LLN} would be no less work than proving the result directly. As mentioned in \Cref{rem:compactness}, the compactness of $\Theta$ is only used to prove the tightness, and could be removed. We will prove a stronger result in \Cref{thm:wasserstein_L1_conv} in the next section. 
Moreover, we will present a way to describe the limiting behaviour of a particle with a given covariate. We will employ an analogue of what is known as the ``nonlinear'' equation or McKean--Vlasov equation in the literature on propagation of chaos \citep{Sznitman1991Topics}.  

\section{Propagation of chaos}
\label{sec:propagation_of_chaos}
We will use a coupling argument to prove the propagation of chaos phenomenon in the strong notion of $L^1(\Omega, \history{}, \prob)$-convergence. To be more precise, we will create an \ac{iid} collection of particles that share the covariates and the driving noise terms (the individual \acp{PRM} $\PRM_{i}^{(1)}, \PRM_{i}^{(2)}$ for $i=1,\ldots, n$) with the \ac{IPS} $(\nX_1, \nX_2, \ldots, \nX_n)$ solving the system of \acp{SDE} in \eqref{eq:combined_sde_for_i_th_individual}. We will show that the process $(\nX_1, \nX_2, \ldots, \nX_k)$ converges to the system $(X_1, X_2, \ldots, X_k)$ in $L^1(\Omega, \history{}, \prob)$ sense as $n\to \infty$. To this end, let us first create the  system of \ac{iid} particles. 

\subsection{A  system of \acl{iid} particles}
\label{sec:iid_ips}
Let $\mu$ be the continuous (with respect to the time parameter) probability measure-valued solution to the integral equation \eqref{eq:limiting_ode} in \Cref{thm:flln}. 
Consider the following  system of \ac{iid} particles solving the following \acp{SDE}: 
\begin{align}
    \begin{aligned}
        S_i(t) &= S_i(0) -  
        \int_{(0, t] \times \setOfPositiveReals} \indicator{[0,  S_{i}(s-) \measureIntegral{\beta(\theta_i, \cdot)}{\mu_{s}(\{e_2\}\times \cdot)} ]}{v}\PRM_{i}^{(1)}(\differential{s}\times\differential{v})\eqcomma \\
        I_i(t) &=  I_{i}(0) + 
        \int_{(0, t] \times \setOfPositiveReals} \indicator{[0, S_{i}(s-) \measureIntegral{\beta(\theta_i, \cdot)}{\mu_{s}(\{e_2\}\times \cdot)} ]}{v}\PRM_{i}^{(1)}(\differential{s}\times\differential{v})\\
        &{}\quad \quad  -  \int_{(0, t] \times \setOfPositiveReals} \indicator{[0, \gamma(\theta_i) I_{i}(s-) ]}{v}\PRM_{i}^{(2)}(\differential{s}\times\differential{v})\eqcomma \\
        R_i(t) & = 1- S_i(t) - I_i(t)\eqcomma \text{ for } i = 1, 2, \ldots, n\eqstop 
    \end{aligned}
    \label{eq:prop_iid_ith_sde}
\end{align} 
Here, we have reused the initial data including the  covariates $\theta_1, \theta_2, \ldots, \theta_n$, and the driving noise terms $\PRM_i^{(1)}$ and $\PRM_i^{(2)}$ for $i=1, \ldots, n$, from the original system of particles in \Cref{sec:stoch_model}. Let $X_i (t) \defeq (S_i(t), I_i(t), R_i(t))$ for $i=1, \ldots, n$ and $t\ge 0$. Then, $(X_1, X_2, \ldots, X_n)$ solves the following system of \acp{SDE}: 
\begin{align}
    \begin{aligned}
        X_i(t) &= X_i(0) + (e_2-e_1) 
        \int_{(0, t] \times \setOfPositiveReals} \indicator{[0, S_{i}(s-) \measureIntegral{\beta(\theta_i, \cdot)}{\mu_{s}(\{e_2\}\times \cdot)} ]}{v}\PRM_{i}^{(1)}(\differential{s}\times\differential{v})\\ 
        &{}\quad \quad  + (e_3-e_2) \int_{(0, t] \times \setOfPositiveReals} \indicator{[0, \gamma(\theta_i) I_{i}(s-) ]}{v}\PRM_{i}^{(2)}(\differential{s}\times\differential{v})\eqcomma \quad \text{for } i =1, 2, \ldots, n\eqstop  
    \end{aligned}
    \label{eq:prop_combined_iid_ith_sde}
\end{align}
In the next section, we will show that the sequence of processes $(\nX_1, \nX_2, \ldots, \nX_k)$ converges to $(X_1, X_2, \ldots, X_k)$ in $L^1(\Omega, \history{}, \prob)$, for any finite $k\ge 1$. However, before proving this result, we will analyse a useful martingale, which will be used in the proof of the propagation of chaos. 

\subsection{A useful martingale and convergence in $L^1(\Omega, \history{}, \prob)$}
\label{sec:propagation_of_chaos_strong_L1}

For bounded, measurable $f : \Delta_3\times \Theta \to \setOfReals$, consider the zero-mean square-integrable martingale 
    \begin{align*}
        \nm_{f}(t) &= \measureIntegral{f}{\nmu_t} -\measureIntegral{f}{\nmu_0} - \int_{0}^{t} \measureIntegral{\Lambda_f(\cdot, \nmu_s) }{\nmu_s} \differential{s}\eqstop 
    \end{align*}
The following lemma is immediate. 
\begin{myLemma}
For each bounded, measurable $f : \Delta_3\times\Theta \to \setOfReals$, 
    the sequence $\{\sup_{s\le \cdot}\absolute{\nm_f(s)} :n \ge 1\}$ converges to zero in $L^p(\Omega, \history{}, \prob)$  as $n\to \infty$ for any $p\ge 1$. That is, 
    \begin{align*}
        \lim_{n\to \infty} \Eof{\left(\sup_{s\le t}\absolute{\nm_f(s)}\right)^p} = 0\eqcomma \quad \text{for all } t \in [0, \infty)\eqcomma p\ge 1\eqstop 
    \end{align*}
    \label{lem:martingale_L2_f}
\end{myLemma}
\begin{proof}[Proof of \Cref{lem:martingale_L2_f}]
    The (optional) quadratic variation of the martingale $\nm_f$ is given by 
    \begin{align*}
        \optionalVariation{\nm_f}(t) = \sum_{s \le t } \left(\nm_f(s) - \nm_f(s-)\right)^2 \leq \frac{Cn\norm{f}_{\infty}^2 }{n^2} \le \frac{C}{n}\eqcomma 
    \end{align*}
    for some constant $C$, 
    since there are at most $2n$ jumps possible in the time interval $[0, t]$ (at most one infection and one recovery event for each particle), and each jump size is bounded by $2 \norm{f}_{\infty}/n$ by \eqref{eq:f_mu_modulus_of_continuity}. Therefore, by the \ac{BDG} inequality \citep[p. 534]{Budhiraja2019Analysis}, we have 
    \begin{align*}
        \Eof{\left(\sup_{s\le t}\absolute{\nm_f(s)}\right)^p} \le C_p {\Eof{\optionalVariation{\nm_f}(t)^{\frac{p}{2}}}} \le C_p {\left(\frac{C}{n}\right)^{\frac{p}{2}}} \to 0\eqcomma 
    \end{align*}
    for constants $C, C_p$, and for all $t\in [0, \infty)$. This completes the proof. 
\end{proof}

The following theorem proves convergence of the empirical measure in bounded Lipschitz metric (see Appendix \ref{sec:IPM}) in $L^1(\Omega, \history{}, \prob)$ sense. In the proof of \Cref{thm:wasserstein_L1_conv}, we will show that the martingale $\nm_f$ converges to zero uniformly in $f$ as well. This is possible since the set $\mathbb{F}_{\mathsf{BL}} \defeq \{f \in C(\Delta_3\times\Theta, \setOfReals): \norm{f}_{\infty} + \norm{f}_{\mathsf{Lip}} \le 1\}$ is totally bounded when $\Theta$ is compact, in the light of the Arzel\'a--Ascoli theorem \citep[Chapter VI, Theorem 3.8]{Conway2018FunctionalAnalysis}.

\begin{myTheorem}
    Assume $\{\nX_1(0), \nX_2(0), \ldots, \nX_n(0)\}$ is an \ac{iid} collection of $\Delta_3$-valued random variables. Furthermore, assume $\theta_1, \theta_2, \ldots, \theta_n$ are \ac{iid} random variables taking values in $\Theta$, a compact subset of the Euclidean space. Moreover, assume that $\beta : \Theta \times \Theta \to\setOfPositiveReals$ and $\gamma: \Theta \to \setOfPositiveReals$ are Lipschitz functions with $\norm{\beta}_{\mathsf{Lip}} \le 1$ and $\norm{\gamma}_{\mathsf{Lip}} \le 1$.
    Then, the following  convergence holds in  $L^1(\Omega, \history{}, \prob)$: 
    \begin{align*}
        \lim_{n\to \infty} \Eof{ \sup_{t\le T} \boundedLipschitz{\nmu_t, \mu_t}} = 0\eqcomma \quad \text{for all } 0< T < \infty\eqstop 
    \end{align*}
    \label{thm:wasserstein_L1_conv}
\end{myTheorem}
\begin{proof}[Proof of \Cref{thm:wasserstein_L1_conv}]
        Following similar calculations as in the proof of \Cref{lemma:flln_uniqueness} and \Cref{thm:flln}, for each bounded $f \in \mathbb{F}_{\mathsf{BL}} \subset \mathbb{F}_{\mathsf{TV}}$,  we have 
    \begin{align*}
        \absolute{\measureIntegral{f}{\nmu_t- \mu_t}} &\le \absolute{\measureIntegral{f}{\nmu_0 - \mu_0} } + \absolute{ \nm_{f}(t)} + B \int_0^t \boundedLipschitz{\nmu_s, \mu_s} \differential{s}\eqcomma  
    \end{align*}
    where the constant $B$ is given by 
      $  B \defeq  \left( 2 \hat{\beta} + \hat{\gamma}  \right) 2 $.
    Taking supremum over all $f\in \mathbb{F}_{\mathsf{BL}}$, we get 
    \begin{align}
        \boundedLipschitz{\nmu_t, \mu_t} & \le \boundedLipschitz{\nmu_0, \mu_0} + \sup_{f \in \mathbb{F}_{\mathsf{BL}}}\absolute{\nm_f(t)} + B \int_0^t \boundedLipschitz{\nmu_s, \mu_s} \differential{s}\eqstop 
        \label{eq:wasserstein_1_pre_Gronwall}
    \end{align}

    Since $\Theta$ is compact, by the Arzel\'a--Ascoli theorem (see \cite[Chapter VI, Theorem 3.8, p. 175]{Conway2018FunctionalAnalysis}), the set $\mathbb{F}_{\mathsf{BL}}$ is totally bounded. Therefore, for each $\vep >0$, we can find a finite open cover of size $N(\vep)\in \setOfNonnegativeIntegers$ and each of radius less than $\vep$. That is, there exists $f_1, f_2, \ldots, f_{N(\vep)} \in \mathbb{F}_{\mathsf{BL}}$ such that for any $f \in \mathbb{F}_{\mathsf{BL}}$, 
    \begin{align*}
        \norm{f-f_i}_{\infty} \le \vep\eqcomma 
    \end{align*}
    for some $i=1, 2, \ldots, N(\vep)$. The number $N(\vep)$ is called the covering number (see \cite[Definition 1.4.1]{Talagrand2021Upper}, or \cite[Section 10.1, pp. 72-73]{Lifshits2012Lectures}). Using the triangle inequality $$\absolute{\nm_f(t)} \le \absolute{\nm_{f_i}(t)} + \absolute{\nm_{f}(t) - \nm_{f_i}(t)},$$ and then taking supremum over $f\in \mathbb{F}_{\mathsf{BL}}$, we get the following decomposition
    \begin{align*}
        \sup_{f \in \mathbb{F}_{\mathsf{BL}}}\absolute{\nm_f(t)} &\le \max_{i = 1, 2, \ldots, N(\vep)} \absolute{\nm_{f_i}(t)} + \sup_{f, g \in \mathbb{F}_{\mathsf{BL}}: \norm{f-g}_{\infty} < \vep}\absolute{\nm_{f}(t) - \nm_{g}(t)}\\
        &\le \sum_{i=1}^{N(\vep)} \absolute{\nm_{f_i}(t)} + \vep C_t \eqcomma 
    \end{align*}
    where the term $C_t$ can be chosen as 
    $
        C_t \defeq 2 \left(1+ 2 (\hat{\beta} + \hat{\gamma}) t \right)\eqcomma 
    $
    since 
    \begin{align*}
        \absolute{\nm_f(t) - \nm_g(t)} &\le \absolute{\measureIntegral{f-g}{\nmu_t} } + \absolute{\measureIntegral{f-g}{\nmu_0} } \\
        &\quad 
        + \absolute{\int_0^t \left( \measureIntegral{\Lambda_f(\cdot, \nmu_s) }{\nmu_s} - \measureIntegral{\Lambda_g(\cdot, \nmu_s) }{\nmu_s}\right)\differential{s} } \\
        &\le 2 \norm{f-g}_{\infty} + 2 \norm{f-g}_{\infty} (\hat{\beta} + \hat{\gamma}) t        
        \eqstop 
    \end{align*}
    Now, it follows from \eqref{eq:wasserstein_1_pre_Gronwall} that 
    \begin{align*}
        \sup_{t \le T} \boundedLipschitz{\nmu_t, \mu_t} &\le  \boundedLipschitz{\nmu_0, \mu_0} + \sup_{t \le T} \sup_{f \in \mathbb{F}_{\mathsf{BL}}}\absolute{\nm_f(t)} + B \int_0^T  \sup_{r \le s}\boundedLipschitz{\nmu_r, \mu_r} \differential{s}\eqstop 
    \end{align*}
    Applying Gr\"onwall's inequality, we get 
    \begin{align}
        \sup_{t \le T} \boundedLipschitz{\nmu_t, \mu_t} &\le  \left(\boundedLipschitz{\nmu_0, \mu_0} + \sup_{t \le T} \sup_{f \in \mathbb{F}_{\mathsf{BL}}}\absolute{\nm_f(t)}\right) \myExp{BT}\eqstop 
        \label{eq:wasserstein_1_Gronwall_bound}
    \end{align}
    Now we have
    \begin{align*}
        \Eof{\sup_{t \le T} \sup_{f \in \mathbb{F}_{\mathsf{BL}}}\absolute{\nm_f(t)}} &{}\le \sum_{i=1}^{N(\vep)} \Eof{\sup_{t \le T} \absolute{\nm_{f_i}(t)}} + \vep C_T \eqstop 
    \end{align*}
    Taking limit as $n \to \infty$, we get 
    \begin{align*}
        \limsup_{n\to \infty} \Eof{\sup_{t \le T} \sup_{f \in \mathbb{F}_{\mathsf{BL}}}\absolute{\nm_f(t)}} \le \vep C_T \eqcomma 
    \end{align*}
    since each $\Eof{\sup_{t \le T} \absolute{\nm_{f_i}(t)}} \to 0$ as $n \to \infty$ by \Cref{lem:martingale_L2_f}.
    Since $\vep$ is arbitrary, taking limit as $\vep \to 0$, we get 
    \begin{align*}
        \lim_{n\to \infty}\Eof{\sup_{t \le T} \sup_{f \in \mathbb{F}_{\mathsf{BL}}}\absolute{\nm_f(t)}} &{}=0 \eqcomma 
    \end{align*}
     Now, by virtue of the \ac{SLLN}, we have 
    $\lim_{n\to \infty}\Eof{\boundedLipschitz{\nmu_0, \mu_0}}= 0$. Therefore, the proof of \Cref{thm:wasserstein_L1_conv} is now complete when we take expectation and then limit as $n\to \infty$ in \eqref{eq:wasserstein_1_Gronwall_bound}.
    
\end{proof}

Our proof is inspired by the Dudley's bound employed in the study of extremes of Gaussian processes \citep[Theorem 1.4.2]{Talagrand2021Upper}. See also \cite[Section 10]{Lifshits2012Lectures}. We are now ready to state and prove the main result on propagation of chaos. 




\begin{myTheorem}
    Let $\Theta$ be a compact subset of the Euclidean space. Assume that  $(\nX_1, \nX_2, \ldots, \nX_n)$ and $(X_1, X_2, \ldots, X_n)$ are solutions to the \acp{SDE} \eqref{eq:combined_sde_for_i_th_individual} and \eqref{eq:prop_combined_iid_ith_sde} respectively. Here, we additionally assume that $\beta : \Theta \times \Theta \to \setOfPositiveReals$ is a Lipschitz function with $\norm{\beta}_{\mathsf{Lip}} \le 1$.
    Then, for any fixed positive integer $k$, the sequence of stochastic processes $(\nX_1, \nX_2, \ldots, \nX_k)$ converges to $(X_1, X_2, \ldots, X_k)$ in the $L^1(\Omega, \history{}, \prob)$ sense, \ie, for all $T\ge 0$, we have
    \begin{align*}
        \lim_{n\to \infty} \Eof{ \sup_{t\le T} \norm{(\nX_1(t), \nX_2(t), \ldots, \nX_k(t)) - (X_1(t), X_2(t), \ldots, X_k(t)) } } = 0\eqstop 
    \end{align*}
    \label{thm:L1_prop_chaos}
\end{myTheorem}
\begin{proof}[Proof of \Cref{thm:L1_prop_chaos}]
    
For each $i$, note that
\begin{align*}
    \sup_{t\le T} \absolute{\nS_i(t) - S_i(t)} &
    \le \int_{(0, T]\times \setOfPositiveReals}\left|\indicator{[0, \nS_{i}(s-) \measureIntegral{\beta(\theta_i, \cdot)}{\nmu_{s-}(
    \{e_2\}\times \cdot)} ]}{v}  \right.\\
    &{}\quad \quad \quad \left.
    - \indicator{[0, S_{i}(s-) \measureIntegral{\beta(\theta_i, \cdot)}{\mu_{s}(
    \{e_2\}\times \cdot)} ]}{v}\right|\PRM_{i}^{(1)}(\differential{s}\times\differential{v}) \eqstop 
\end{align*}
Therefore, we have 
\begin{align*}
    \Eof{\sup_{t\le T} \absolute{\nS_i(t) - S_i(t)}} &\le \int_{0}^{T} \E\left[\left|\nS_{i}(s-) \measureIntegral{\beta(\theta_i, \cdot)}{\nmu_{s-}(
    \{e_2\}\times \cdot)} \right. \right. \\
    &\quad \quad \left. \left. - S_{i}(s-) \measureIntegral{\beta(\theta_i, \cdot)}{\mu_{s}(
    \{e_2\}\times \cdot)} \right| \right] \differential{s}\\
    &\le \nerror_{S, i}(T) + \int_{0}^{T} \left(\sup_{r\le s} h(r)\right) \Eof{ \sup_{r\le s} \absolute{\nS_i(r) - S_i(r)}} \differential{s}\eqcomma 
\end{align*}
from which Gr\"onwall's inequality gives 
\begin{align}
    \Eof{\sup_{t\le T} \absolute{\nS_i(t) - S_i(t)}} &\le \nerror_{S, i}(T) \myExp{  \int_{0}^{T} \left(\sup_{r\le s} h(r)\right)  \differential{s}}\eqcomma
    \label{eq:S_i_L1_Gronwall} 
\end{align}
where $h(r) \defeq \sup_{\theta \in \Theta }\measureIntegral{\beta(\theta, \cdot)}{\mu_{r}(
    \{e_2\}\times \cdot)} < \hat{\beta} < \infty$.
    The error term $\nerror_{S, i}(T)$ is given by 
\begin{align*}
    \nerror_{S, i}(T) & \defeq \int_{0}^{T} \Eof{\nS_{i}(s)  \absolute{\measureIntegral{\beta(\theta_i, \cdot)}{\nmu_{s}(
    \{e_2\}\times \cdot)} - \measureIntegral{\beta(\theta_i, \cdot)}{\mu_{s}(
    \{e_2\}\times \cdot)}  }} \differential{s}\\
    &{} \le \int_{0}^{T} \Eof{ \absolute{\measureIntegral{\beta(\theta_i, \cdot)}{\nmu_{s}(
    \{e_2\}\times \cdot)} - \measureIntegral{\beta(\theta_i, \cdot)}{\mu_{s}(
    \{e_2\}\times \cdot)}  }} \differential{s} \\
    &{} \le 2 ( 1 + \hat{\beta} ) \int_{0}^{T} \Eof{ \boundedLipschitz{\nmu_{s},  \mu_{s}} }\differential{s} \\
    &{} \le 2 ( 1 + \hat{\beta} ) \int_{0}^{T} \Eof{ \sup_{r \le s} \boundedLipschitz{\nmu_{r}, \mu_{r}} }\differential{s}
\end{align*}
by virtue of \Cref{thm:wasserstein_L1_conv} as $n$ goes to infinity,  $\nerror_{S, i}(T)$ goes to $0$ . Therefore, from \eqref{eq:S_i_L1_Gronwall}, we have 
\begin{align}
    \lim_{n\to \infty}\Eof{\sup_{t\le T} \absolute{\nS_i(t) - S_i(t)}} =0 \eqstop 
    \label{eq:S_i_L1_conv_zero}
\end{align} 
Let us now consider the terms corresponding to the infectious status 
\begin{align*}
    \sup_{t\le T}\absolute{\nI_i(t) - I_i(t)} & \le \int_{(0, T]\times \setOfPositiveReals} \absolute{\indicator{[0, \nS_{i}(s-) \measureIntegral{\beta(\theta_i, \cdot)}{\nmu_{s-}(
    \{e_2\}\times \cdot)} ]}{v} - \indicator{[0, S_{i}(s-) \measureIntegral{\beta(\theta_i, \cdot)}{\mu_{s}(
    \{e_2\}\times \cdot)} ]}{v}} \\
    &\quad\quad  \times \PRM_{i}^{(1)}(\differential{s}\times\differential{v}) \\
    &{}\quad 
    + \int_{(0, T]\times \setOfPositiveReals} \absolute{\indicator{[0, \nI_{i}(s-) \gamma(\theta_i)]}{v} - \indicator{[0, I_{i}(s-) \gamma(\theta_i) ]}{v}}\PRM_{i}^{(2)}(\differential{s}\times\differential{v})\eqstop 
\end{align*}
Therefore, we have  
\begin{align*}
    \Eof{ \sup_{t\le T}\absolute{\nI_i(t) - I_i(t)} } & \le \nerror_{S, i}(T) + \int_{0}^{T} \left(\sup_{r\le s} h(r)\right) \Eof{ \sup_{r\le s} \absolute{\nS_i(r) - S_i(r)}} \differential{s}\\
    &{}\quad 
    + \hat{\gamma} \int_{0}^{T} \Eof{\sup_{r\le s} \absolute{\nI_i(r) - I_i(r)}  }\differential{s}\eqstop 
\end{align*}
By an application of the Gr\"onwall's inequality, we have 
\begin{align*}
   \Eof{ \sup_{t\le T}\absolute{\nI_i(t) - I_i(t)} } & \le \left(\nerror_{S, i}(T) + \int_{0}^{T} \left(\sup_{r\le s} h(r)\right) \Eof{ \sup_{r\le s} \absolute{\nS_i(r) - S_i(r)}} \differential{s} \right) \times \myExp{ \hat{\gamma} T }
\end{align*}
Again, by virtue of \Cref{thm:wasserstein_L1_conv}, and \eqref{eq:S_i_L1_conv_zero}, the right-hand side converges to zero as $n$ tends to infinity. Therefore, 
\begin{align}
    \lim_{n\to \infty}\Eof{\sup_{t\le T} \absolute{\nI_i(t) - I_i(t)}} =0 \eqcomma  
    \label{eq:I_i_L1_conv_zero}
\end{align}
and 
\begin{align}
    \lim_{n\to \infty}\Eof{\sup_{t\le T} \absolute{\nR_i(t) - R_i(t)}} =0 \eqcomma  
    \label{eq:R_i_L1_conv_zero}
\end{align}
since $\nR_i(t) = 1- \nS_i(t) - \nI_i(t)$, and $R_i(t) = 1 -  S_i(t) - I_i(t)$ for all $t\ge 0$. Combining \eqref{eq:S_i_L1_conv_zero}, \eqref{eq:I_i_L1_conv_zero}, and \eqref{eq:R_i_L1_conv_zero}, we have 
\begin{align}
    \lim_{n\to \infty}\Eof{\sup_{t\le T} \norm{\nX_i(t) - X_i(t)}} =0 \eqstop  
    \label{eq:X_i_L1_conv_zero}
\end{align}
Since $i$ is arbitrarily chosen, we have for each fixed $k$, 
\begin{align*}
    \lim_{n\to \infty}\Eof{\sum_{i=1}^{k}\sup_{t\le T} \norm{\nX_i(t) - X_i(t)}} =0 \eqcomma 
\end{align*}
which completes the proof. 

\end{proof}

\section{Examples and applications}
\label{sec:applications}
In this section, we discuss some potential applications of our results either to epidemic models with specific choices of infection and recovery rates $\beta$ and $\gamma$, or to statistical inference. 
    For a bounded, continuous function $f: \Delta_3\times \Theta \to \setOfReals$, define the real-valued function $g_f(t) \defeq \measureIntegral{f}{\mu_t}$, where $\mu$ is the solution of the differential equation in \eqref{eq:limiting_ode}. Then $g_f$ is a solution of the differential equation 
    \begin{align}
        \timeDerivative{g_f}(t) &= \measureIntegral{\Lambda_f(\cdot, \mu_t)}{\mu_t} \\
        & = \int_{\Theta} \left(f(e_2, \theta)-f(e_1, \theta)  \right)  \int_{\Theta} \beta(\theta, \theta') \mu_t(\{e_2\}\times  \differential{\theta'}) \mu_t(\{e_1\} \times \differential{\theta}) \\
        & \quad + \int_{\Theta}\gamma(\theta) \left(f(e_3, \theta) - f(e_2, \theta)\right) \mu_t(\{e_2\}\times \differential{\theta}) \eqcomma 
        \label{eq:g_f_ODE}
    \end{align}
    with initial data $g_f(0) = \measureIntegral{f}{\mu_0}$. A particularly useful choice is $f(x, \theta) = \indicator{\{e_1\}}{x} \indicator{\Theta}{\theta}$ for which the function $g_f$ gives us the proportion of susceptible individuals.  Similarly, we can define $f(x, \theta) = \indicator{\{e_2\}}{x} \indicator{\Theta}{\theta}$ for which the function $g_f$ gives us the proportion of infected individuals. Therefore, we get the following system of \acp{ODE}
    \begin{align}
        \begin{aligned}
            \timeDerivative{S}(t) &= -  \int_{\Theta}  \int_{\Theta} \beta(\theta, \theta') \mu_t(\{e_2\}\times  \differential{\theta'}) \mu_t(\{e_1\} \times \differential{\theta})  \eqcomma \\
        \timeDerivative{I}(t) &= \int_{\Theta} \int_{\Theta} \beta(\theta, \theta') \mu_t(\{e_2\}\times  \differential{\theta'}) \mu_t(\{e_1\} \times \differential{\theta}) \\
        & \quad - \int_{\Theta}\gamma(\theta) \mu_t(\{e_2\}\times \differential{\theta}) \eqcomma \\
        \timeDerivative{R}(t) &= \int_{\Theta}\gamma(\theta) \mu_t(\{e_2\}\times \differential{\theta}) \eqcomma 
        \end{aligned}
        \label{eq:limiting_SIR_ode}
    \end{align}
    with initial data $S(0) = \mu_0(\{e_1\}\times \Theta)$, $I(0) = \mu_0(\{e_2\}\times \Theta)$, and $R(0) = \mu_0(\{e_3\}\times \Theta)$, where $S(t) \defeq \measureIntegral{\indicator{\{e_1\}}{\cdot} \indicator{\Theta}{ \cdot}}{\mu_t}, I(t) \defeq \measureIntegral{\indicator{\{e_2\}}{\cdot} \indicator{\Theta}{ \cdot}}{\mu_t}, R(t) \defeq \measureIntegral{\indicator{\{e_3\}}{\cdot} \indicator{\Theta}{ \cdot}}{\mu_t}$ denote the total proportion of susceptible, infected and recovered individuals in the population at time $t\ge 0$.

     Since $\timeDerivative{S(t) + I(t) + R(t)} = 0$, we have $S(t) + I(t) + R(t) = S(0) + I(0) + R(0)$ for all $t\ge 0$. That is, $\mu_t(\Delta_3\times \Theta) = \mu_0(\Delta_3\times \Theta) = 1$ for all $t\ge 0$. Moreover, choosing any bounded, continuous function $f: \Delta_3\times \Theta \to \setOfReals$ that does not depend on $x$, \ie, $f(x, \theta) = f(\theta)$, we have $\measureIntegral{f}{\mu_t} = \measureIntegral{f}{\mu_0}$ for all $t\ge 0$. 




\subsection{Examples of the limiting differential equations}
\label{sec:mean_field_limit_examples}

\subsubsection{Standard \ac{SIR} model}

    Consider the case where $\beta$ and $\gamma$ are constants. Then, the system of \acp{ODE} in \eqref{eq:limiting_SIR_ode} for the proportions of individuals in the three compartments reduces to the following system of autonomous \acp{ODE}
    \begin{align}
        \begin{aligned}
            \timeDerivative{S}(t) &= -  \beta_{\star} S(t) I(t) \eqcomma \\
        \timeDerivative{I}(t) &= \beta_{\star} S(t) I(t) - \gamma_{\star} I(t) \eqcomma \\
        \timeDerivative{R}(t) &= \gamma_{\star} I(t) \eqcomma 
        \end{aligned}
        \label{eq:standardSIR_ODE}
    \end{align}
    with initial data $S(0) = \mu_0(\{e_1\}\times \Theta)$, $I(0) = \mu_0(\{e_2\}\times \Theta)$, and $R(0) = \mu_0(\{e_3\}\times \Theta)$. The system of \acp{ODE} in \eqref{eq:standardSIR_ODE} is a special case of the integro-differential renewal equations presented by Kermack and McKendrick in  \cite{Kermack1927Contribution}. 

\subsubsection{Standard \ac{SIR} model with frailty}
\label{sec:frailty}
    Consider the frailty model $\beta(\theta, \theta') = \beta_{\star} \theta$ and $\gamma(\theta) = \gamma_{\star}$, with $\Theta$ being a compact interval of $\setOfPositiveReals$. Then, the system of \acp{ODE} in \eqref{eq:limiting_SIR_ode} reduces to the following system of  \acp{ODE}
    \begin{align}
        \begin{aligned}
            \timeDerivative{S}(t) &= - \beta_{\star} I(t)  \int_{\Theta}  \theta  \mu_t(\{e_1\} \times \differential{\theta})   \eqcomma \\
        \timeDerivative{I}(t) &= \beta_{\star} I(t)  \int_{\Theta}  \theta  \mu_t(\{e_1\} \times \differential{\theta})  - \gamma_{\star} I(t) \eqcomma \\
        \timeDerivative{R}(t) &= \gamma_{\star} I(t) \eqcomma 
        \end{aligned}
    \end{align}
    with initial data $S(0) = \mu_0(\{e_1\}\times \Theta)$, $I(0) = \mu_0(\{e_2\}\times \Theta)$, and $R(0) = \mu_0(\{e_3\}\times \Theta)$. The aggregate equations are not closed in $S(t), I(t)$ and $R(t)$, because the incidence depends on the evolving first moment of susceptibility among the susceptible population. In order to get more insights into the dynamics of the frailty model, we define the \ac{MGF} of the measure $\mu_t(\{e_1\}\times \cdot)$ as  $$ m_S(z,t) \defeq  \int_{\Theta} e^{z\theta} \mu_t(\{e_1\}\times \differential \theta) \eqcomma  $$ for $z\in \setOfReals$. Taking $f_z(x,\theta) = \indicator{\{e_1\}}{x}\myExp{z\theta} $ in \eqref{eq:g_f_ODE}, we obtain the following \ac{PDE} 
    \begin{align}
        \begin{aligned}
            \frac{\partial}{\partial t} m_S(z,t) + \beta_{\star} I(t) \frac{\partial }{\partial z}m_S(z,t) & = 0\eqcomma\\
            \text{with } m_S(z,0)= \int_{\Theta} e^{z\theta} \mu_0(\{e_1\}\times \differential{\theta} )  &= S(0) \int_{\Theta} e^{z\theta} \mu_0(\Delta_3\times \differential{\theta} ) = S(0)  \mathsf{M}(z)\eqcomma 
        \end{aligned}
         \label{eq:SIR_frailty_PDE}
    \end{align}   
     since $\mu_0(\{e_1\}\times \differential{ \theta}) = \mu_0(\{e_1\}\times \Theta) \mu_0(\Delta_3\times \differential{ \theta})$ due to the independence of the allocation of the initial immunological statuses and the colours, and where $ \mathsf{M}(z) \defeq  \int_{\Theta} \myExp{z\theta} \mu_0(\Delta_3\times \differential \theta) $ is the \ac{MGF} of the colours (covariates). Let 
    \begin{align*}
        \timeDerivative{A(t)} = \beta_{\star} I(t)\eqcomma \text{ with } A(0) = 0\eqcomma 
    \end{align*}
    so that $A(t) = \beta_{\star}\int_0^t I(s)\differential{s}$ is the cumulative infection pressure. The \ac{IVP} in \eqref{eq:SIR_frailty_PDE} can be solved using standard methods \citep{Evans2010PDEs} to get
    \begin{align*}
        m_S(z,t) = S(0) \mathsf{M}(z-A(t))\eqcomma 
    \end{align*}
    whence using $ S(t) = m_S(0,t) $ we obtain 
    \begin{align}
    \begin{aligned}
    S(t) &= S(0) \mathsf{M}(-A(t))\eqcomma \\
    I(t) &= I(0) + S(0)\left( 1-\mathsf{M}(-A(t)) \right) - \frac{\gamma_\star}{\beta_{\star}}A(t)\eqcomma \\
    R(t) &= R(0) + \frac{\gamma_\star}{\beta_{\star}}A(t)\eqstop 
    \end{aligned}
    \label{eq:SIR_frailty_solutions}
    \end{align}
    \Cref{fig:sir-comparison_frailty} compares the stochastic trajectories generated using \Cref{alg:Gillespie} with the limiting deterministic equations.

\begin{myRemark}
    Suppose the measure $\mu_0(\{e_1\}\times \differential{\theta})$ (the distribution of $\theta$) admits a density (Radon--Nikodym derivative) $$h(\theta) \defeq \frac{ \mu_0(\{e_1\}\times \differential{\theta})}{\leb(\differential{\theta})}$$ with respect to the Lebesgue measure. Then, note that the \ac{MGF} of the measure $\mu_t(\{e_1\}\times \differential{\theta})$ has the following representation
    \begin{align*}
        m_S(z,t) = S(0) \mathsf{M}(z-A(t)) 
        = S(0) \int_{\Theta} e^{ ( z -A(t) ) \theta} \mu_0(\Delta_3 \times \differential \theta) 
         = \int_{\Theta} e^{ z \theta } \Big( S(0) \ e^{ -A(t) \theta} \ h(\theta) \Big) \differential \theta.
    \end{align*}
    Therefore, by the uniqueness of the Laplace transform, the measure $\mu_t( \{e_1\} \times \differential{\theta} )$ admits a density $S(t, \theta)$ with respect to the Lebesgue measure, which is given by 
    \begin{align*}
        S(t, \theta) \defeq \frac{\mu_t( \{e_1\} \times \differential{\theta} ) }{\leb(\differential{\theta})} = S(0) e^{ -A(t) \theta}  h(\theta)\eqstop 
    \end{align*}
    Differentiating with respect to $t$, one obtains 
    \begin{align}
        \timeDerivative{S(t, \theta)} = -\beta_\star I(t) \theta S(t, \theta)\eqcomma \text{ for } t\ge 0\eqstop 
        \label{eq:frailty_S_density_ODE}
    \end{align}
    Differential equations such as \eqref{eq:frailty_S_density_ODE} are popular in mathematical biology literature. See \cite{gomes2022individual,Gomes2014,Margheri2015Correlation,Margheri2017Vaccine}. Rigorous derivations of such differential equations for the densities directly from the stochastic model are scarce. In the susceptibility-only case, the measure-valued limit recovers the scalar Laplace-transform reduction obtained from the Sellke construction in \cite{Izyumtseva2026Sellke}. The Sellke representation gives a particularly explicit description in this special case, while the present framework extends to transmission kernels depending on both susceptibility and infectiousness.

\end{myRemark}

    \begin{minipage}{\textwidth}
    \captionsetup{type=figure}
    \centering

    \begin{subfigure}[tbh]{0.48\linewidth}
        \centering
        \includegraphics[
            width=\linewidth,
            height=1.98\textheight,
            keepaspectratio
        ]{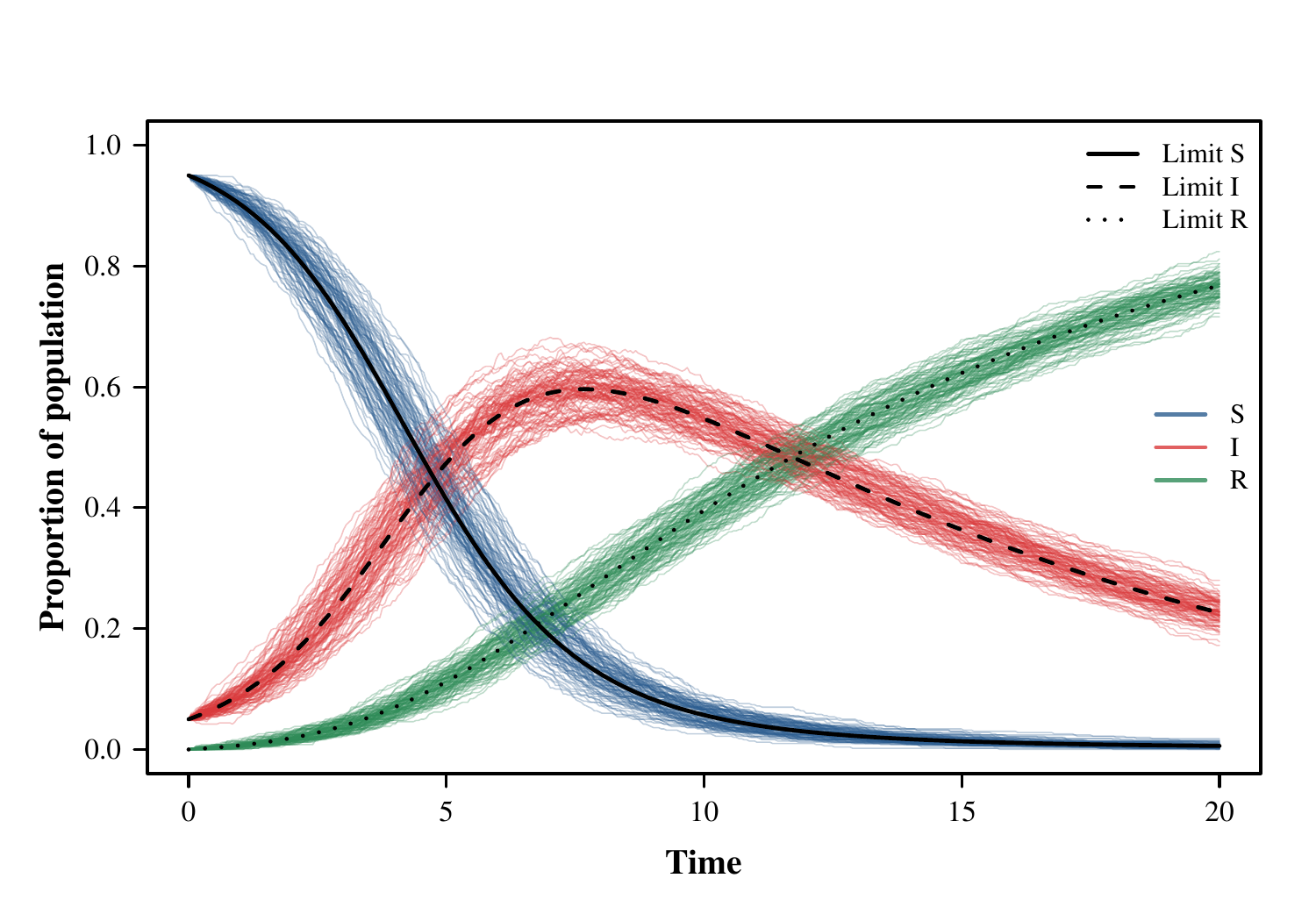}
        \caption{ $ n= 500$, $\theta$
        $\sim$ \Beta{5}{1}, $100$ trajectories.}
    \end{subfigure}
    \hfill
    \begin{subfigure}[tbh]{0.48\linewidth}
        \centering
        \includegraphics[
            width=\linewidth,
            height=1.98\textheight,
            keepaspectratio
        ]{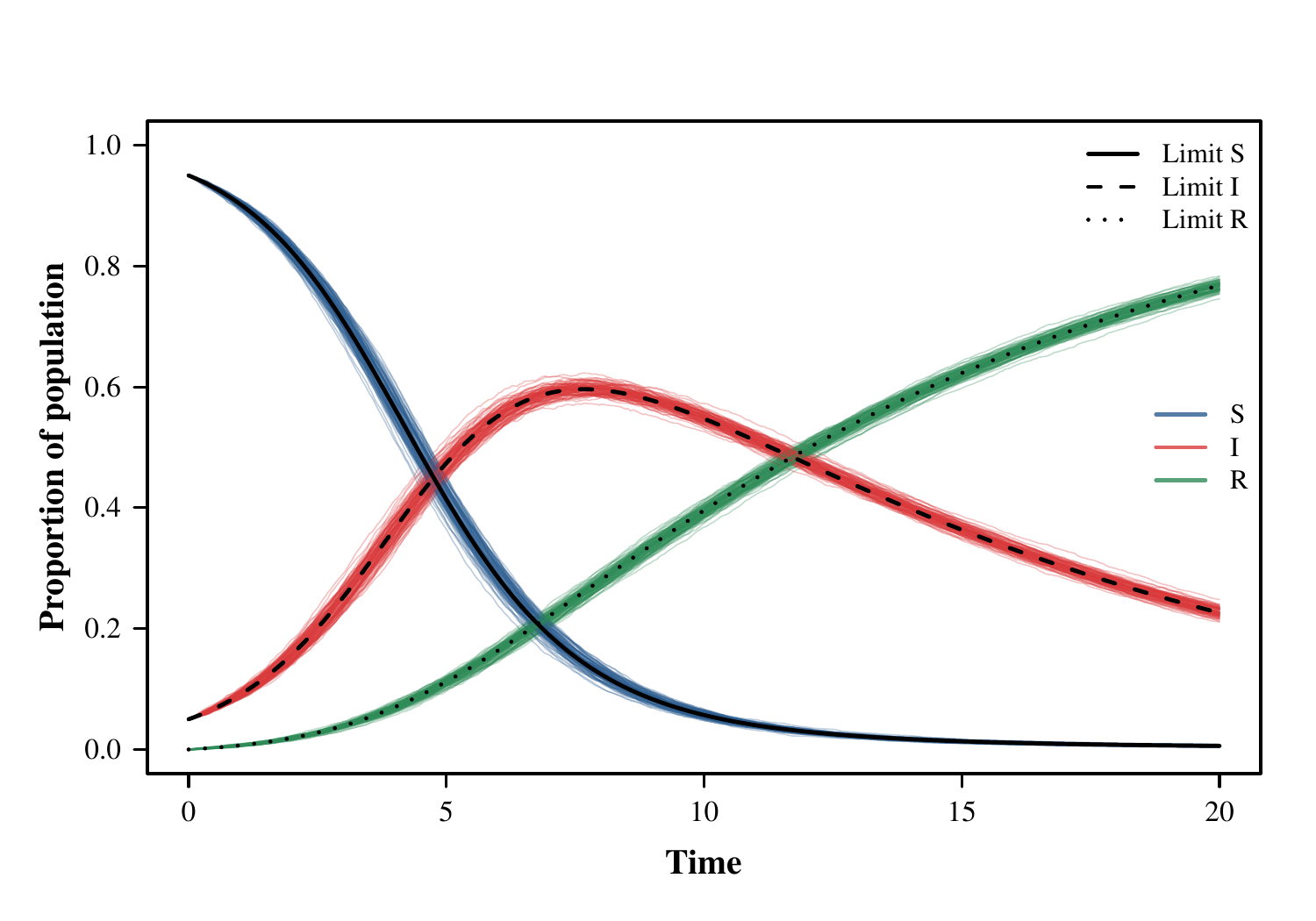}
        \caption{$ n= 5000$, $\theta$
        $\sim$ \Beta{5}{1}, $100$ trajectories.}
    \end{subfigure}


    \begin{subfigure}[tbh]{0.48\linewidth}
        \centering
        \includegraphics[
            width=\linewidth,
            height=1.98\textheight,
            keepaspectratio
        ]{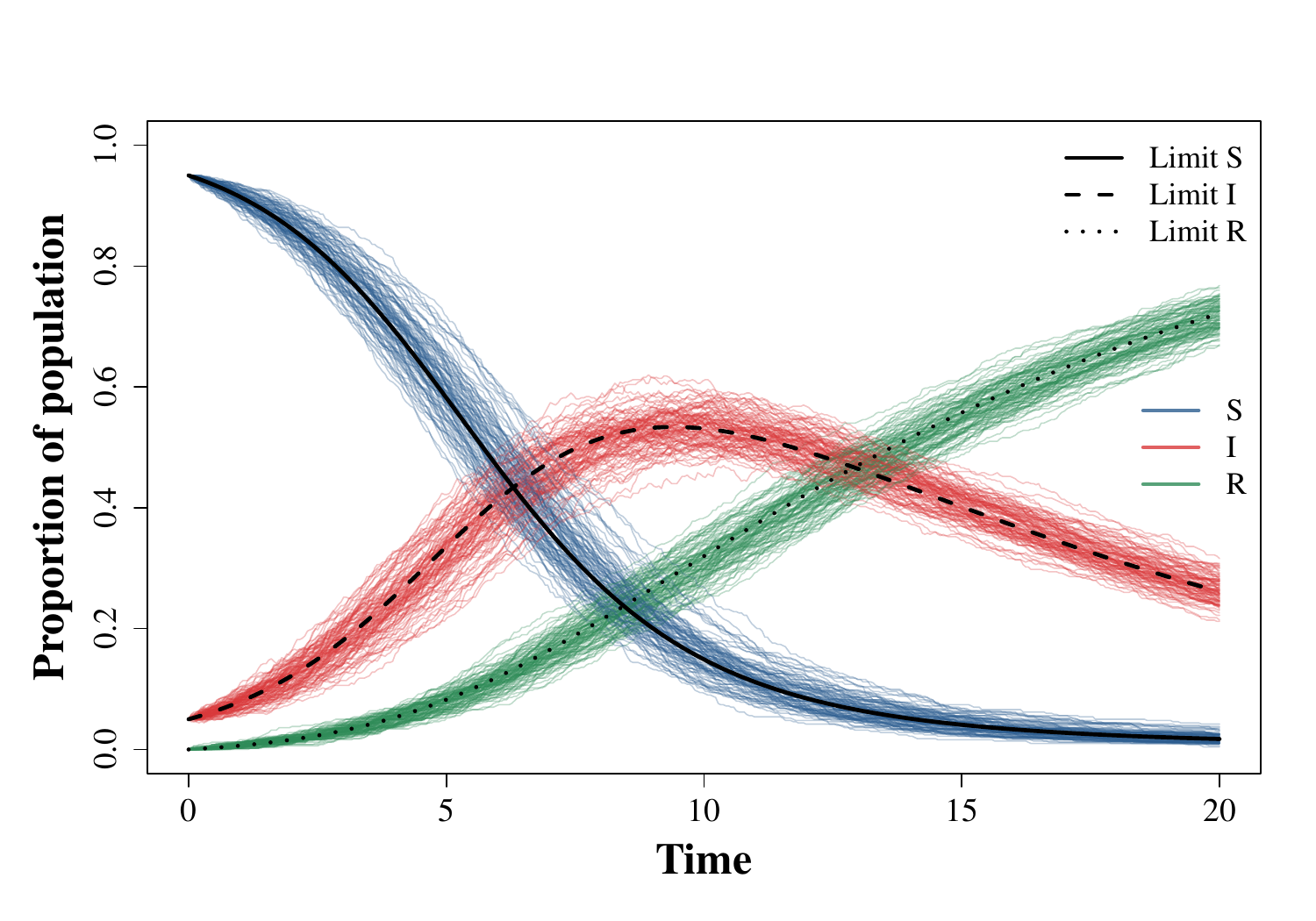}
        \caption{$ n= 500$, $\theta$
        $\sim$ \Unif{0}{1}, $100$ trajectories.}
    \end{subfigure}
    \hfill
    \begin{subfigure}[tbh]{0.48\linewidth}
        \centering
        \includegraphics[
            width=\linewidth,
            height=1.98\textheight,
            keepaspectratio
        ]{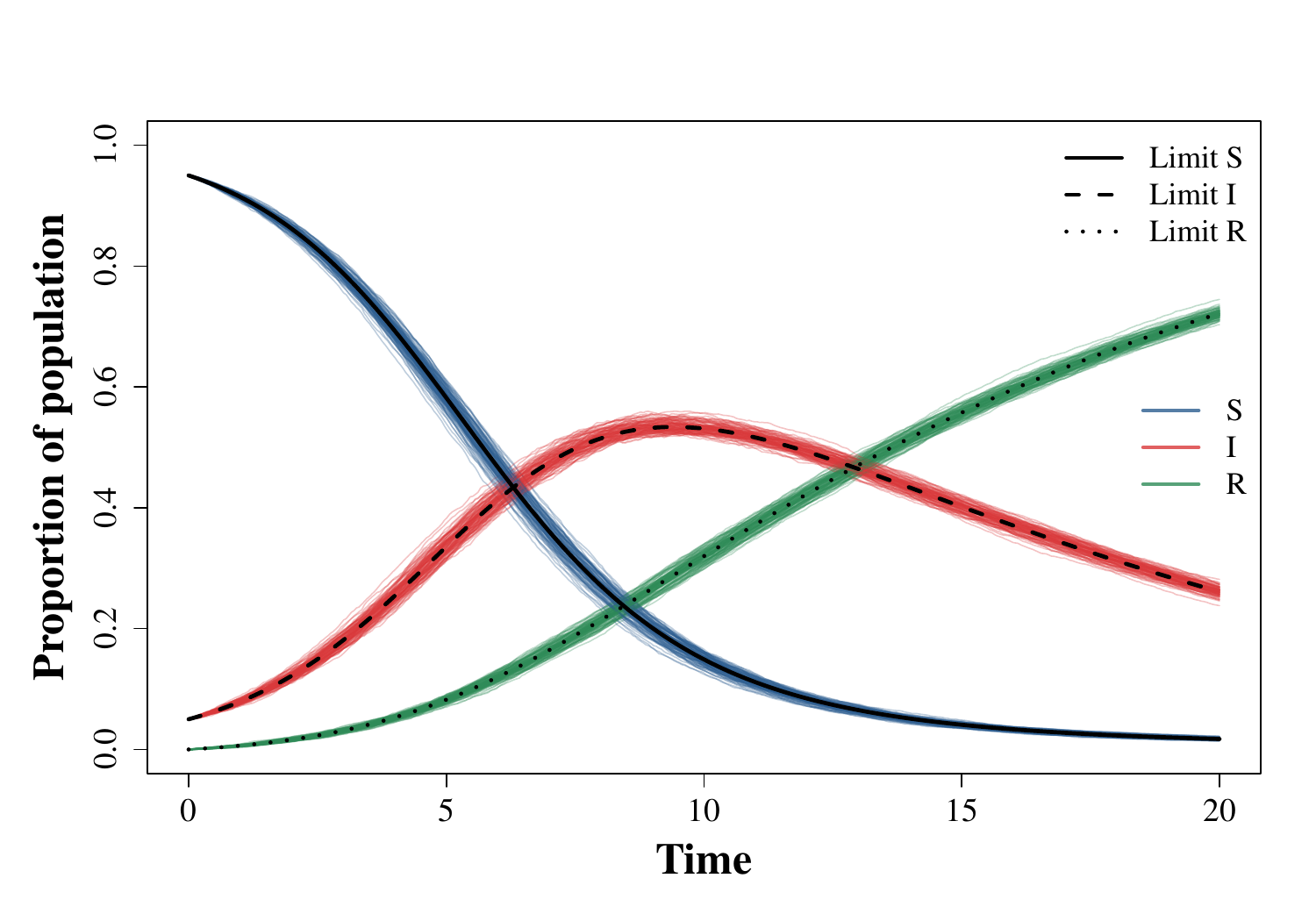}
        \caption{$ n= 5000$, $\theta$
        $\sim$ \Unif{0}{1}, $100$ trajectories.}
    \end{subfigure}

    \begin{subfigure}[tbh]{0.48\linewidth}
        \centering
        \includegraphics[
            width=\linewidth,
            height=1.98\textheight,
            keepaspectratio
        ]{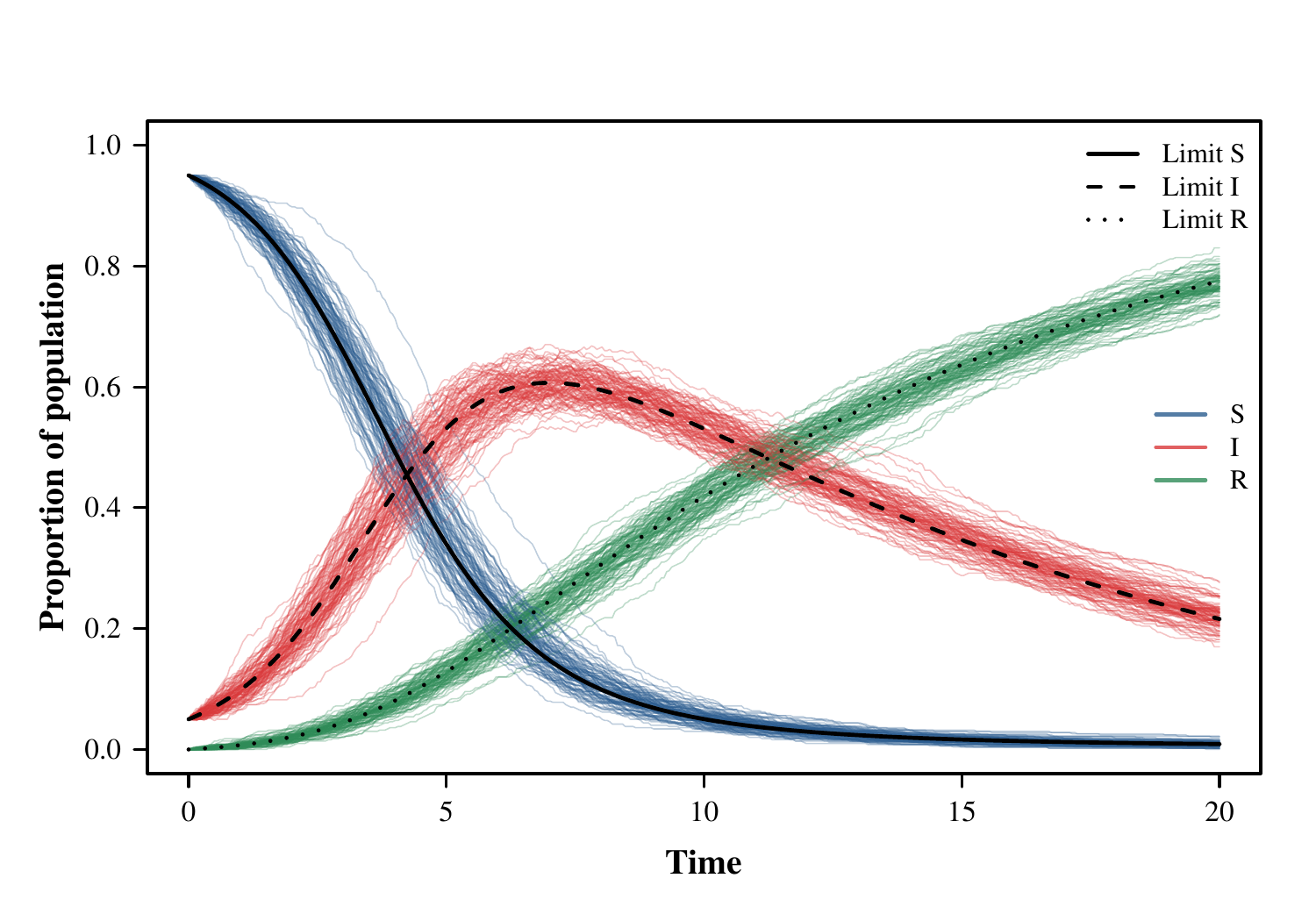}
        \caption{ $ n= 500$, $\theta$
        $\sim$ \Bin{3}{0.6}, $100$ trajectories.}
    \end{subfigure}
    \hfill
    \begin{subfigure}[tbh]{0.48\linewidth}
        \centering
        \includegraphics[
            width=\linewidth,
            height=1.98\textheight,
            keepaspectratio
        ]{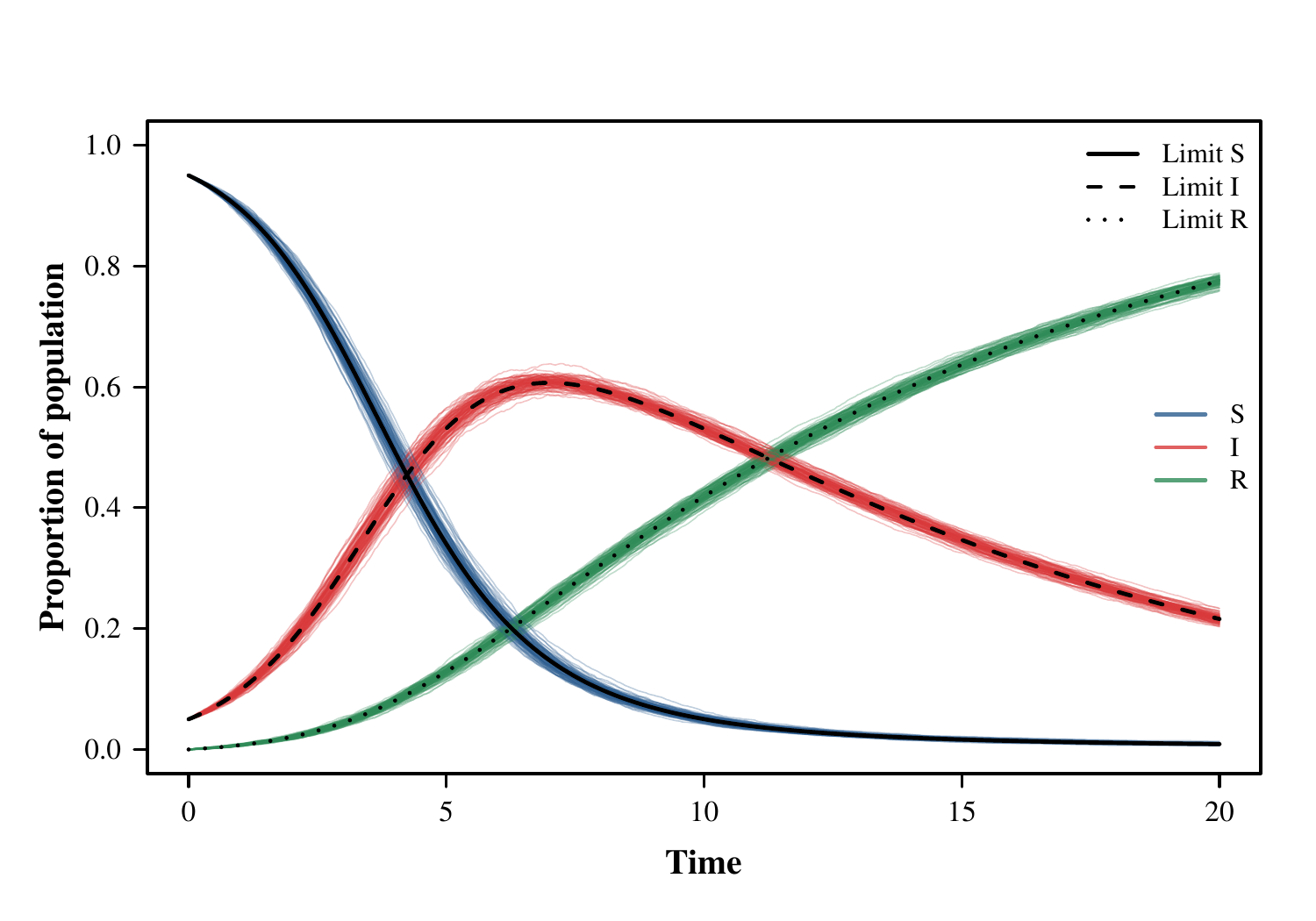}
        \caption{$ n= 5000$, $\theta$
        $\sim$ \Bin{3}{0.6}, $100$ trajectories.}
    \end{subfigure}

    \caption{Comparison of stochastic trajectories of the proportions of susceptible, infected, and recovered population over time obtained via the Doob--Gillespie algorithm (see \Cref{alg:Gillespie}) against the limiting deterministic equations in black dotted lines for three different choices of distribution of the covariates. Here, the infection rate is $\beta ( \theta_i, \theta_j ) = \beta_{\star} \theta_i $ and recovery rate is $\gamma(\theta) = \gamma_{\star}$, under the frailty model of \Cref{sec:frailty}. The figures on the right panel demonstrate that the accuracy of the \ac{FLLN} approximation improves with increasing $n$.
    }
    \label{fig:sir-comparison_frailty}
\end{minipage}

\subsubsection{A variation of the multitype model}
\label{sec:multitype}

    Consider $\beta(\theta, \theta') = \beta_{\star} \theta \theta'$, and $\gamma(\theta) = \gamma_{\star}$, with $\Theta$ being a compact subset of $\setOfPositiveReals$. As before, let us define the \acp{MGF} of the measures $\mu_t(\{e_1\}\times \differential{\theta})$ and $\mu_t(\{e_2\}\times \differential{\theta})$ as follows
    \begin{align*}
        m_S(z, t) \defeq \int_{\Theta} e^{z\theta} \mu_t(\{e_1\}\times \differential \theta)\eqcomma \text{ and } m_I(z, t) \defeq \int_{\Theta} e^{z\theta} \mu_t(\{e_2\}\times \differential \theta)\eqcomma 
    \end{align*}
    for $t\ge 0$ and $z\in \setOfReals$. Let us also define the corresponding first moments of the measures $\mu_t(\{e_1\}\times \differential{\theta})$ and $\mu_t(\{e_2\}\times \differential{\theta})$:
    \begin{align*}
        m_S^{(1)}(t) \defeq \left. \frac{\partial m_S(z, t) }{\partial z }\right|_{z=0} = \int_{\Theta} \theta\, \mu_t(\{e_1\}\times \differential \theta)\eqcomma \\
         m_I^{(1)}(t) \defeq \left. \frac{\partial m_I(z, t) }{\partial z }\right|_{z=0} = \int_{\Theta} \theta\, \mu_t(\{e_2\}\times \differential \theta)\eqstop 
    \end{align*}
    Then, setting $ f(x,\theta) = \mathbf{1}_{\{e_1\}}(x)e^{z\theta}$ for the susceptible population, and $ f(x,\theta) = \mathbf{1}_{\{e_2\}}(x)e^{z\theta}$ for the infected population in \eqref{eq:g_f_ODE}, we get the following \acp{PDE}
    \begin{align}
    \begin{aligned}
        \frac{\partial }{\partial t}m_S(z, t) + \beta_{\star} m_I^{(1)}(t) \frac{\partial }{\partial z} m_S(z, t) &=0\eqcomma \\
         \frac{\partial }{\partial t} m_I(z, t) - \beta_{\star} m_I^{(1)}(t) \frac{\partial}{\partial z}  m_S(z, t) & =- \gamma_\star m_I(z, t)\eqcomma \\
         m_S(z,0) = S(0)  \mathsf{M}(z)\eqcomma\quad  m_I(z,0) &= I(0)  \mathsf{M}(z)\eqcomma 
    \end{aligned}
         \label{eq:multitype_limiting_PDEs}
    \end{align}
    where, as before,  $ \mathsf{M}(z) \defeq  \int_{\Theta} \myExp{z\theta} \mu_0(\Delta_3\times \differential \theta) $ is the \ac{MGF} of the colours (covariates).
    From \eqref{eq:multitype_limiting_PDEs}, we can describe the limiting proportions of the susceptible, the infected, and the recovered individuals as 
    \begin{align}
    \begin{aligned}
        \timeDerivative{S(t)} &=-\beta_{\star} m_S^{(1)}(t)m_I^{(1)}(t)\eqcomma \\
        \timeDerivative{I(t)} &=\beta_{\star} m_S^{(1)}(t)m_I^{(1)}(t)-\gamma_\star I(t)\eqcomma 
        \\
        \timeDerivative{R(t)} &=\gamma_\star I(t).
    \end{aligned}
    \end{align}
    In \Cref{fig:sir-comparison_product}, we  compare the stochastic trajectories generated via \Cref{alg:Gillespie} with the deterministic limiting equations. As we can see from \Cref{fig:sir-comparison_product}, the limiting deterministic equations approximate the stochastic dynamics well when the population size is large.  

\begin{minipage}{\textwidth}
    \captionsetup{type=figure}
    \centering

    \begin{subfigure}[tbh]{0.48\linewidth}
        \centering
        \includegraphics[
            width=\linewidth,
            height=1.98\textheight,
            keepaspectratio
        ]{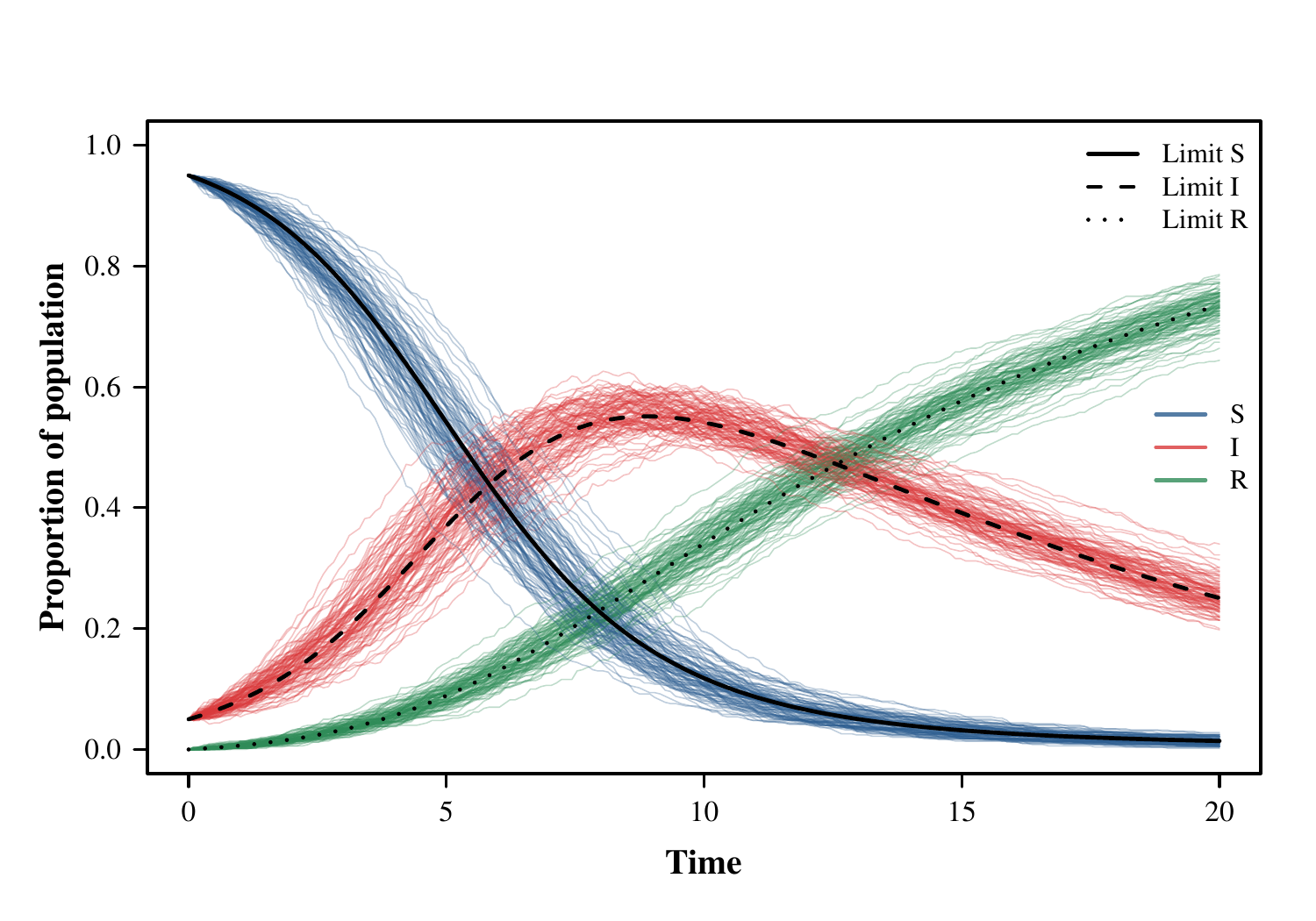}
        \caption{ $ n= 500$, $\theta$
        $\sim$ \Beta{5}{1}, $100$ trajectories.}
    \end{subfigure}
    \hfill
    \begin{subfigure}[tbh]{0.48\linewidth}
        \centering
        \includegraphics[
            width=\linewidth,
            height=1.98\textheight,
            keepaspectratio
        ]{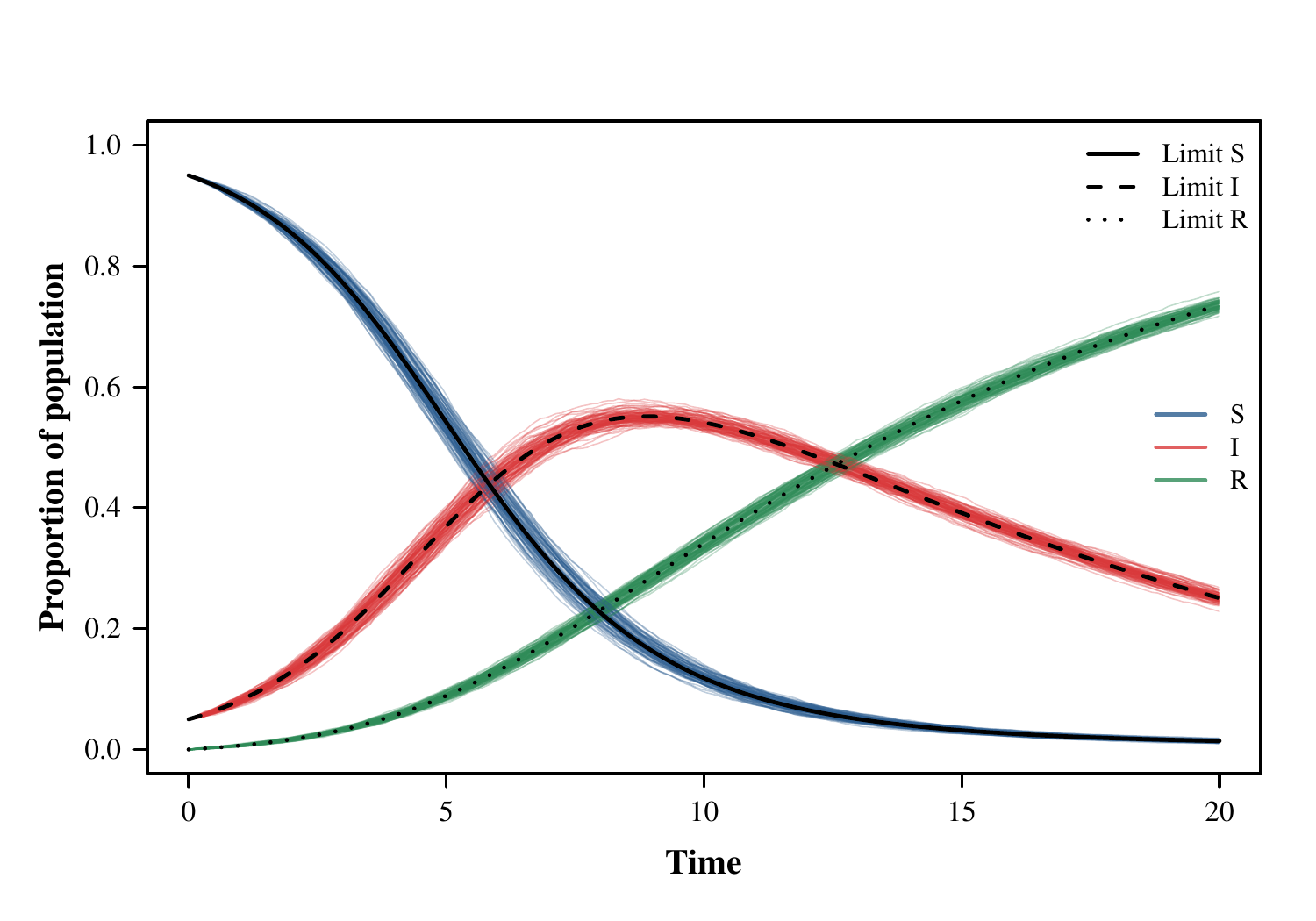}
        \caption{$ n= 5000$, $\theta$
        $\sim$ \Beta{5}{1}, $100$ trajectories.}
    \end{subfigure}


    \begin{subfigure}[tbh]{0.48\linewidth}
        \centering
        \includegraphics[
            width=\linewidth,
            height=1.98\textheight,
            keepaspectratio
        ]{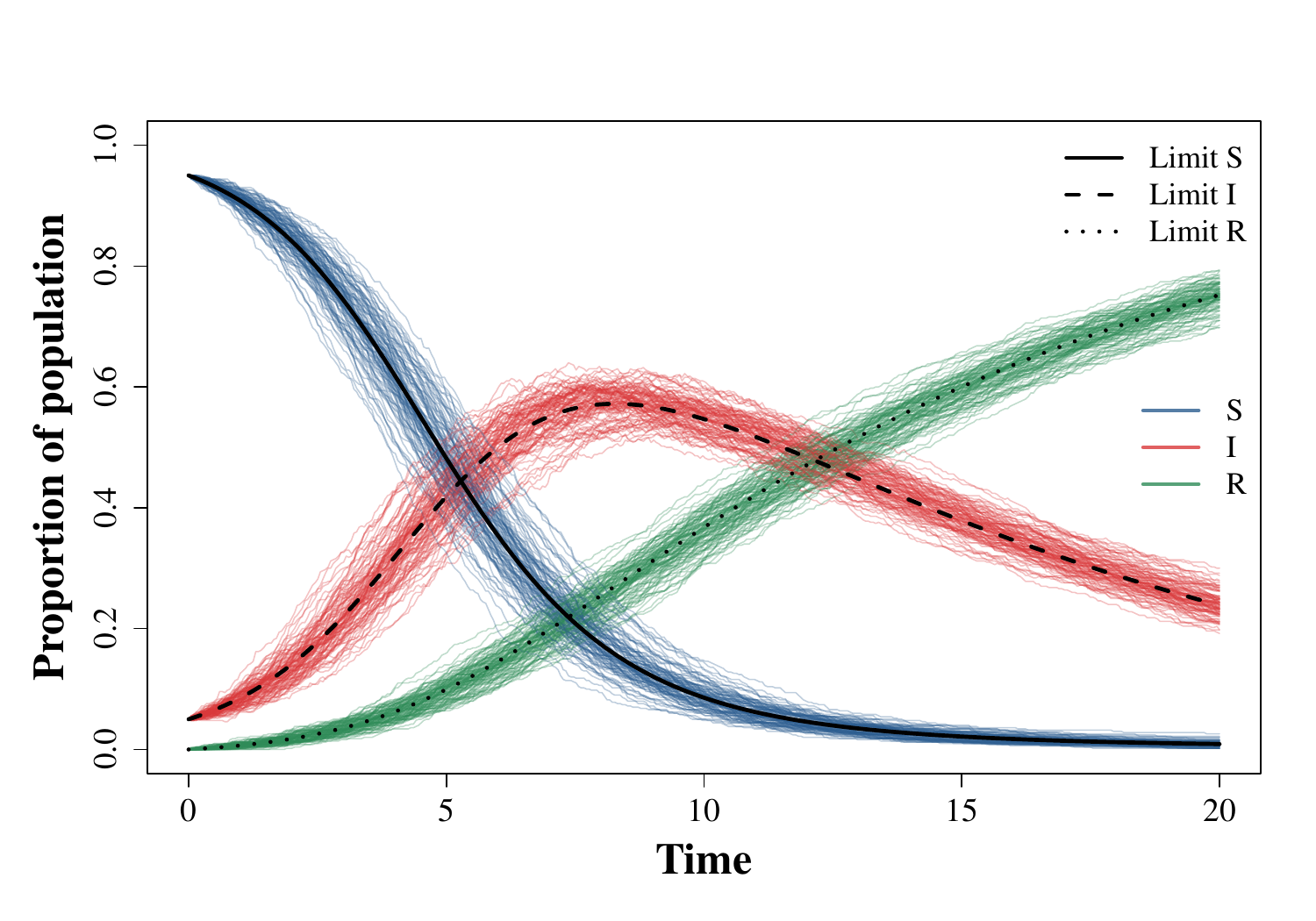}
        \caption{$ n= 500$, $\theta$
        $\sim$ \Unif{0}{1}, $100$ trajectories.}
    \end{subfigure}
    \hfill
    \begin{subfigure}[tbh]{0.48\linewidth}
        \centering
        \includegraphics[
            width=\linewidth,
            height=1.98\textheight,
            keepaspectratio
        ]{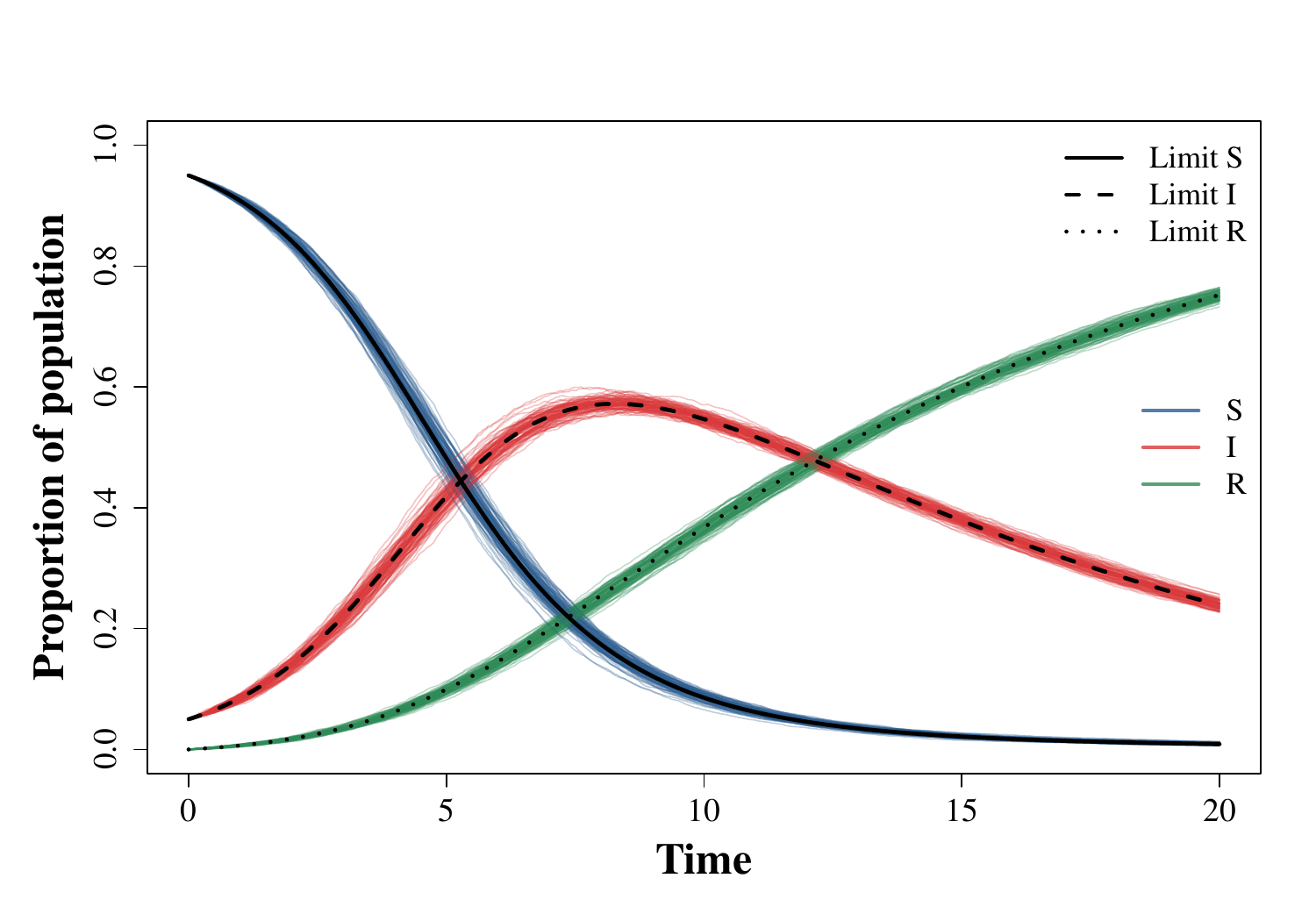}
        \caption{$ n= 5000$, $\theta$
        $\sim$ \Unif{0}{1}, $100$ trajectories.}
    \end{subfigure}

    \begin{subfigure}[tbh]{0.48\linewidth}
        \centering
        \includegraphics[
            width=\linewidth,
            height=1.98\textheight,
            keepaspectratio
        ]{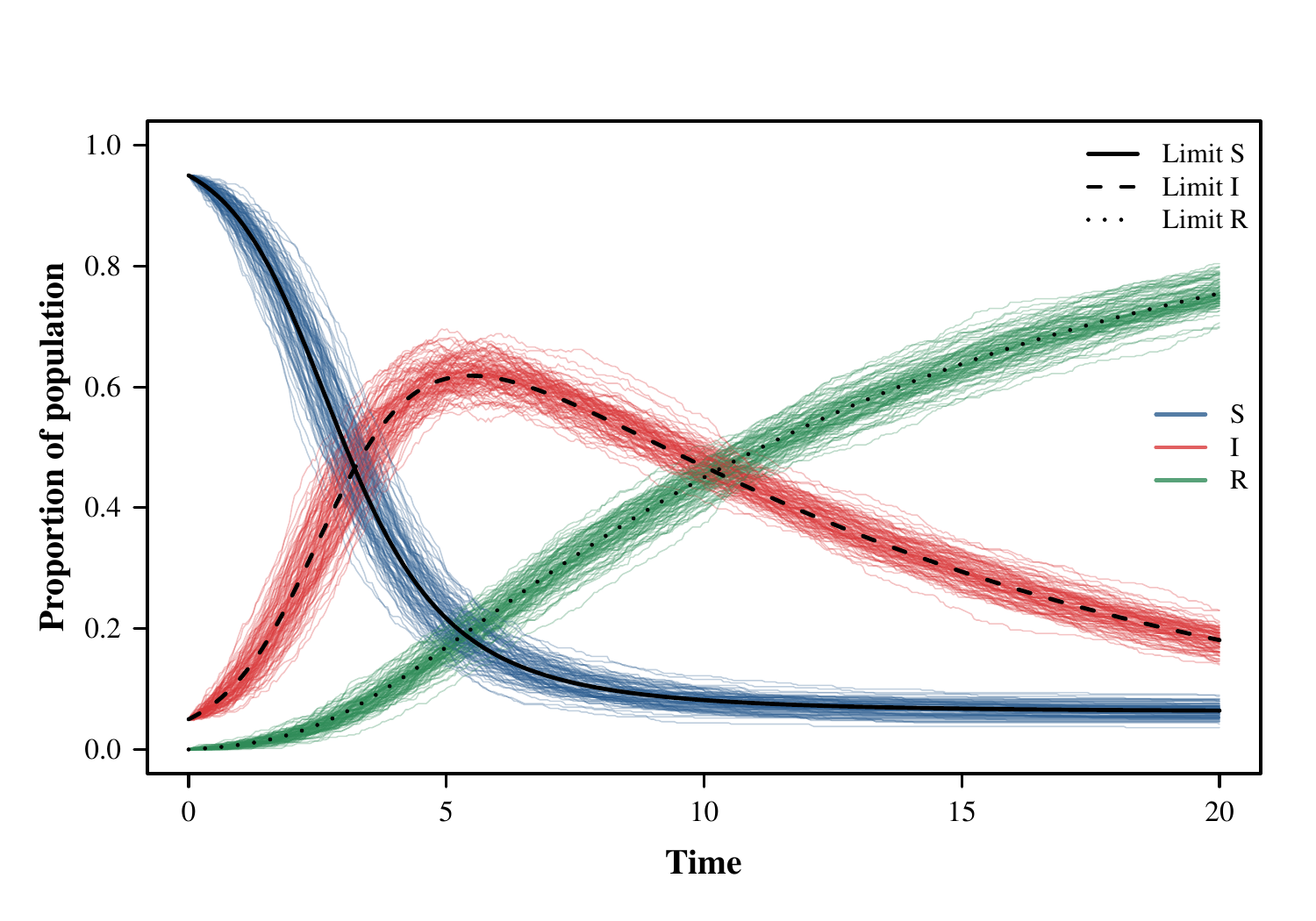}
        \caption{ $ n= 500$, $\theta$
        $\sim$ shifted \Bin{3}{0.6}, $100$ trajectories.}
    \end{subfigure}
    \hfill
    \begin{subfigure}[tbh]{0.48\linewidth}
        \centering
        \includegraphics[
            width=\linewidth,
            height=1.98\textheight,
            keepaspectratio
        ]{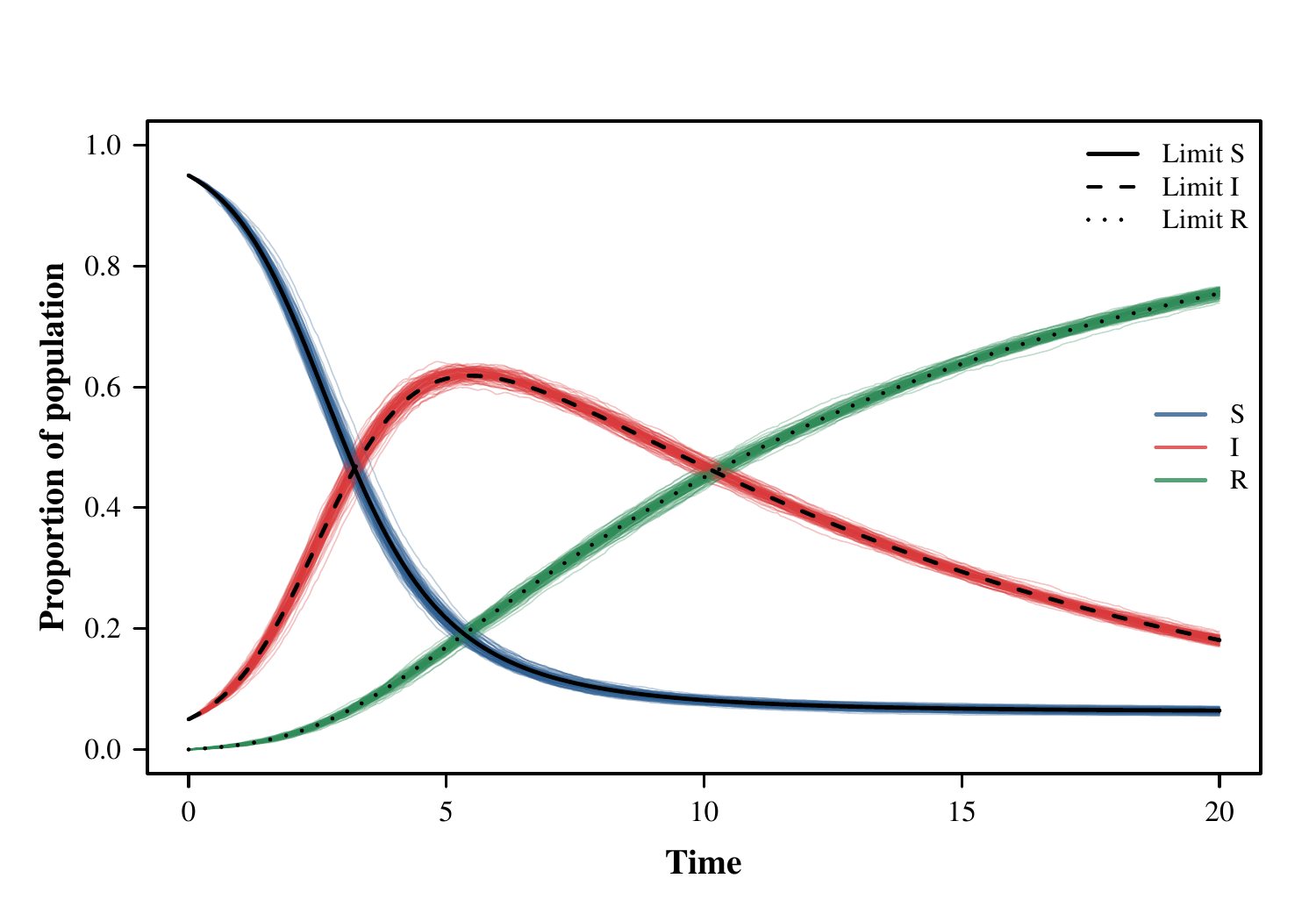}
        \caption{$ n= 5000$, $\theta$
        $\sim$ shifted \Bin{3}{0.6}, $100$ trajectories.}
    \end{subfigure}

    \caption{Comparison of stochastic trajectories of proportions of susceptible, infected, and recovered individuals over time obtained via the  Doob--Gillespie’s algorithm (see \Cref{alg:Gillespie}) with the deterministic limiting equations in black dotted lines for different choices of  distribution of the covariates. Here, the infection rate is $\beta( \theta_i \theta_j ) = \beta_{\star} \theta_i \theta_j $ and recovery rate is $\gamma(\theta) = \gamma_{\star}$, in line with the multitype model in \Cref{sec:multitype}. As we can see by comparing the right panel with the left one, the \ac{FLLN} provides a more accurate approximation when the population size $n$ is large. }
    \label{fig:sir-comparison_product}
\end{minipage}

\subsection{An application of propagation of chaos to  statistics}
\label{sec:parameter_inference}
One of the practical challenges of statistical inference in the context of infectious disease epidemiology is that population-level data are often not available. Instead, we have access to a random sample of  individual-level data, which are possibly partially observed, and/or truncated and censored. Typical data are the infection and/or recovery/death times of a small sample of individuals who either report to a hospital, or are found through contact-tracing efforts.  The \ac{DSA} method \citep{KhudaBukhsh2024HowTo,KhudaBukhsh2020DSA,khuda_bukhsh_2022_projecting,Rempala2023handbook} allows one to perform inference on the parameters of the model based on a random sample of individual-level information on infection and/or recovery times without requiring population-level data. In essence, the method relies on asymptotic independence of the particles, which is a consequence of the propagation of chaos phenomenon. Strictly speaking, a true random sample of infection and/or recovery times is hard to obtain. Nevertheless, the likelihood applies when the observed individuals can reasonably be treated as a random sample, or when the ascertainment mechanism is
conditioned on or modelled separately. 
Since parameter inference is not the main focus of the paper, we  discuss this important application only briefly below. We refer the readers to \cite[Section 4]{Ganguly2026Bharat} for an extensive discussion on the use of propagation of chaos to justify the \ac{DSA} method in the context of \ac{MM} enzyme kinetic reactions. 

The following functional analytic fact about the continuity of jump times of \cadlag functions is straightforward to verify, but we state it for the sake of completeness. For a \cadlag function $f \in D([0, \infty), E)$ taking values in a metric space $E$, let $j_k$ denote its $k$-th jump time. That is, 
\begin{align*}
    j_1(f) &\defeq \inf\{t >0 : f(t) \ne f(t-)\},\\ 
    j_k(f) &\defeq \inf\{t > j_{k-1}(f) : f(t) \ne f(t-)\}\eqcomma \quad \text{for } k\ge 2. 
\end{align*}
Here, we adopt the convention that the infimum of the empty set is infinity. In our model, each particle (individual) can have at most two jumps, one at the time of infection and one at the time of recovery. The following  is immediate.


\begin{myCorollary}
    For each fixed $k\in \setOfNaturals$, we have 
    \begin{align*}
        ((j_1(\nX_i), j_2(\nX_i)): i =1 , 2, \ldots, k) \ConvInDist ((j_1(X_i), j_2(X_i)): i =1 , 2, \ldots, k)
    \end{align*}
    as $n\to \infty$, where $(\nX_1, \nX_2, \ldots, \nX_n)$ solves the system of \acp{SDE} in \eqref{eq:combined_sde_for_i_th_individual}, and $X_1, X_2, \ldots, X_k$ are \ac{iid} satisfying the \acp{SDE} in \eqref{eq:prop_combined_iid_ith_sde}. 
    \label{cor:jump_time_convergence}
\end{myCorollary}
\begin{proof}[Proof of \Cref{cor:jump_time_convergence}]
    Note that $L^1(\Omega, \history{}, \prob)$-convergence  implies weak convergence. Following \Cref{lemma:M_infinity_conv} and the discussion in \Cref{rem:D_infinity_conv}, \Cref{thm:L1_prop_chaos} shows that the sequence of stochastic processes $(\nX_1, \nX_2, \ldots, \nX_k)$, in fact,  converges to $(X_1, X_2, \ldots, X_k)$ in $D([0, \infty), \Delta_3^k)$ by virtue of \cite[Theorems 16.2 and 16.7]{Billingsley1999Convergence}.  Therefore, 
    the proof of \Cref{cor:jump_time_convergence} follows by the continuous mapping theorem \citep[Theorem 5.27]{kallenberg2021foundations}, in the light of \Cref{thm:L1_prop_chaos} and \Cref{lemma:jump_time_continuity_cadlag} in \Cref{sec:additional}. 
\end{proof}

\Cref{cor:jump_time_convergence} implies that the infection and the recovery times of the individuals, which may be infinity, in the original stochastic model of \Cref{sec:stoch_model} converge to those of the \ac{iid} system of particles (individuals) in \Cref{sec:iid_ips}. That is, for any fixed sample size $k$, the joint law of the sampled individual histories converges
to the corresponding product law. Thus, the product likelihood based on the limiting
nonlinear process is asymptotically exact as $n\to \infty$. In this precise sense, a likelihood function based on the limiting \ac{iid} system is \emph{asymptotically exact}. It follows from \eqref{eq:prop_combined_iid_ith_sde}  that, conditioned on the initial data, and the covariate, the distribution of the jump times can be derived as follows:
\begin{align*}
    \probOf{j_1(X_1) > t \mid S_1(0)=1, \theta_1 } &= \myExp{- \int_{0}^{t} \measureIntegral{\beta(\theta_1, \cdot)}{\mu_r(\{e_2\}\times \cdot)} \differential{r}}\eqcomma \\
    \probOf{j_1(X_1) >t \mid I_1(0)=1, \theta_1 } & = \myExp{- \gamma(\theta_1)t }\eqcomma \\
    \probOf{j_2(X_1) > t+s \mid S_1(0)=1,  j_1(X_1) = t, \theta_1} & = \myExp{- \gamma(\theta_1)s }\eqcomma 
\end{align*} 
for $t>0, s>0$. Taking derivative with respect to the time parameter gives us the corresponding densities. When the initial data and the covariate information are not available, the distributions can be calculated using the law of total probability:
\begin{align}
    \begin{aligned}
        F_1(t) \defeq \probOf{j_1(X_1) > t} &= \int_{\Theta} \probOf{j_1(X_1) > t \mid S_1(0)=1, \theta } \mu_0(\{e_1\}\times \differential{\theta}) \\
    &{}\quad 
    + \int_{\Theta}\myExp{- \gamma(\theta)t} \mu_0(\{e_2\}\times \differential{\theta}) + \mu_0(\{e_3\}\times\Theta)\\
    & = \int_{\Theta} \myExp{- \int_{0}^{t} \measureIntegral{\beta(\theta, \cdot)}{\mu_r(\{e_2\}\times \cdot)} \differential{r}} \mu_0(\{e_1\}\times \differential{\theta}) \\
    &{}\quad 
    + \int_{\Theta}\myExp{- \gamma(\theta)t} \mu_0(\{e_2\}\times \differential{\theta}) + \mu_0(\{e_3\}\times\Theta)\eqcomma \\
    &{} = F_{1, S}(t) + F_{1, I}(t) + \mu_0(\{e_3\}\times\Theta)\eqcomma \\
    F_2(t) \defeq \probOf{j_2(X_1) > t} &= \int_{0}^{t} \int_{\Theta} \probOf{j_2(X_1)>t \mid S_1(0)=1, j_1(X_1) =s, \theta} \\
    &{}\quad 
    \times \measureIntegral{\beta(\theta, \cdot)}{\mu_s(\{e_2\}\times \cdot)} \myExp{- \int_{0}^{s} \measureIntegral{\beta(\theta, \cdot)}{\mu_r(\{e_2\}\times \cdot)} \differential{r}} \\
    &\quad 
    \times  \mu_0(\{e_1\}\times \differential{\theta}) \differential{s} + \mu_0(\{e_2, e_3\}\times \Theta) \\
    &{} = F_{2, S}(t) + \mu_0(\{e_2, e_3\}\times \Theta)\eqcomma 
    \end{aligned}
    \label{eq:dsa_survival_functions}
\end{align}
for all $t\ge 0$, where 
\begin{align*}
    F_{1, S}(t) &\defeq \int_{\Theta} \myExp{- \int_{0}^{t} \measureIntegral{\beta(\theta, \cdot)}{\mu_r(\{e_2\}\times \cdot)} \differential{r}} \mu_0(\{e_1\}\times \differential{\theta})\eqcomma \\
    F_{1, I}(t) &\defeq \int_{\Theta}\myExp{- \gamma(\theta)t} \mu_0(\{e_2\}\times \differential{\theta})\eqcomma \\
    F_{2, S}(t) & \defeq F_2(t) - \mu_0(\{e_2, e_3\}\times \Theta)\eqstop 
\end{align*}
The functions $F_1$ and $F_2$, called survival functions in mathematical statistics, completely characterise the probability laws of the first and second jump times of a randomly chosen individual in a large population. 

Let $\Psi$ denote the set of unknown parameters of the system. 
Suppose we observe a random sample of infection times $t_1, t_2, \ldots, t_k >0$ in a time interval $[0, T]$.  
Then, the contribution of the infection times to the likelihood function is given by 
\begin{align}
    \ell_{I}(\Psi \mid t_1, t_2, \ldots, t_k) & = \frac{1  }{ (S(0)- F_{1, S}(T))^k } \prod_{j=1}^k \left(- \timeDerivative{F_{1, S}(t)}\bigg|_{t=t_j} \right) \eqstop 
    \label{eq:likelihood_infection}
\end{align}
If recovery times are available, their contribution to the likelihood can be calculated similarly using the functions $F_{1, I}$ and $F_{2, S}$ conditioning on whether the individual was infected or susceptible initially. 

\begin{myRemark}
    Note that the \ac{DSA} likelihood functions such as $\ell_{I}$ in \eqref{eq:likelihood_infection} have been used in the mathematical biology and applied statistics literature before. The propagation of chaos phenomenon provides a rigorous analytic and probabilistic justification to the likelihood function \citep{Ganguly2026Bharat}. 
\end{myRemark}

\section{Conclusion}
\label{sec:conclusion}
In this paper, we consider a stochastic \ac{SIR} model in terms of a system of \acp{SDE} driven by independent \acp{PRM} where each individual is endowed with covariates that affect their susceptibility, and infectiousness. We prove an \ac{FLLN} that shows that the empirical random measure of the epidemic process converges to a deterministic, continuous (probability) measure-valued function, which is characterised as a weak solution to a differential equation in \Cref{thm:flln}. Under a slightly stricter assumption on the infection and recovery rates, we also prove an \ac{FLLN} in a strong sense in \Cref{thm:wasserstein_L1_conv}. Furthermore, we establish the propagation of chaos phenomenon in \Cref{thm:L1_prop_chaos}, which establishes asymptotic independence between particles, and thereby allows for a likelihood-based statistical inference. Extensions of our model to incorporate non-Markovian dynamics along the lines of \cite{DiLauro2022nonMarkovian} will be pursued in the future.

\appendix

\section{Additional proofs}
\label{sec:additional}

\begin{myDefinition}[Derivation]
    Given an algebra $\mathsf{X}$, and a field or a ring $\mathbb{K}$, a derivation is linear map (with scalars from $\mathbb{K}$) $\mathsf{D}: \mathsf{X} \to \mathsf{X}$ that satisfies the Leibnitz's rule:
    \begin{align*}
        \mathsf{D}(x y) = \mathsf{D}(x)y + x\mathsf{D}(y)\eqcomma 
    \end{align*}
    for all $x, y \in \mathsf{X}$. 
    \label{defn:derivation}
\end{myDefinition}

\begin{myLemma}
    Let $\{(M_t, d_t) : t\ge 0\}$ be a parametric family of metric spaces. Let 
    \begin{align*}
        M_\infty &\defeq \bigcap_{t\ge 0} M_t\eqcomma \\
        d_{\infty}(x, y) &\defeq \int_0^\infty e^{-t} \min\{1, d_t(x, y)\} \differential{t}\eqcomma \text{ for } x, y \in M_\infty\eqstop 
    \end{align*}
    Then, $(M_\infty, d_{\infty})$ is a (possibly empty) metric space. Furthermore, let $y$ and $\{x_n : n\ge 1\} \subset M_{\infty}$ be an element and a sequence in $M_{\infty}$ such that $x_n \to y$ in $(M_t, d_t)$ as $n\to \infty$, for $\leb$-almost all $t\ge 0$. Then, $x_n\to y$ in $(M_{\infty}, d_{\infty})$ as $n\to \infty$.  
    \label{lemma:M_infinity_conv}
\end{myLemma}
\begin{proof}[Proof of \Cref{lemma:M_infinity_conv}]
    It is straightforward that $(M_\infty, d_{\infty})$ is a metric space. The second assertion follows by Lebesgue \ac{DCT}. 
\end{proof}


\begin{myLemma}
    Let $f_n$ and $f$ be elements of $D([0, \infty), E)$,  where  $(E, d)$ is a metric space. Suppose that the jump sizes of $f_n$ and $f$ are uniformly bounded below by a positive number, \ie, there exists a positive number $\delta_{\star}$ such that $d(f_n(t), f_n(t-) ) > \delta_{\star}$ and $d({f(t), f(t-)}) > \delta_{\star}$ if $t$ is a point of (jump) discontinuity. If $f_n \to f$ in $(D([0, T], E), \norm{\cdot}_{\infty})$ for every $T>0$, then $j_k(f_n) \to j_k(f)$ as $n\to \infty$ for all $k\ge 1$. 
    \label{lemma:jump_time_continuity_cadlag}
\end{myLemma}

\begin{proof}[Proof of \Cref{lemma:jump_time_continuity_cadlag}]
     By the assumption in the lemma, given any $\varepsilon>0$ and $T>0$, there exists a positive integer $N(\varepsilon, T)$ such that
     \begin{align}
         f_n(t) \in B_{\varepsilon}(f(t)), \text{ for all } t\le T, \text{ and for all } n\ge N(\varepsilon, T)\eqcomma
         \label{eq:closeness_requirement}
     \end{align}
     where 
     $
         B_{\varepsilon}(f) \defeq \{g\in E : d(g, f) < \varepsilon\}
     $
    is the open $\varepsilon$-ball around $f$ in $(E, d)$. Suppose $ j_1(f) $ is finite.  Set $T= 2 j_1(f)$, and $ \varepsilon = \delta_{\star} /4 $. We claim that  $ j_1(f_n)=j_1(f)$ for every $ n\geq N({\delta_{\star} / 4}, T)$. If not, there must exist $m \ge N({\delta_{\star} / 4}, T)$ such that $j_1(f) \ne t_1 = j_1(f_m) \in [0, \infty] $. 
    If $t_1 < j_1(f)$, we have 
    \begin{align*}
        \delta_\star \le d(f_m(t_1), f_m(t_1-)) \le d(f_m(t_1), f(t_1)) + d(f(t_1), f_m(t_1-)) \le  d(f_m(t_1), f(t_1)) + \varepsilon\eqstop 
    \end{align*}
    But this implies $d(f_m(t_1), f(t_1)) \ge 3\delta_\star/4 > \varepsilon$, a contradiction to \eqref{eq:closeness_requirement}. On the other hand, if $\tilde{t}_1 = j_1(f) < t_1 \le \infty $, then by a similar argument 
    \begin{align*}
        d(f_m(\tilde{t}_1), f(\tilde{t}_1)) \ge \frac{3\delta_\star}{4} > \varepsilon, 
    \end{align*}
    which contradicts \eqref{eq:closeness_requirement} since $\tilde{t}_1 \le T$. Therefore, $j_1(f_m) = j_1(f)$ for all $m\ge N(\varepsilon, T)$. 
    
    For any $k\ge 1$, if $j_k(f)<\infty$, then by repeating similar arguments with $T= 2 j_k(f)$, we can show that for any $\varepsilon>0$, there exists a positive integer $N(\varepsilon, T)$ such that $j_i(f_m) = j_i(f)$ for all $i=1, 2, \ldots, k$, for all $m \ge N(\varepsilon, T)$. 

    Let $k_\star \ge 0$ be such that $j_k(f)<\infty$ for all $k\le k_\star$ and $j_k(f) = \infty$ for all $k>  k_\star$. We claim that $j_k(f_n) \to \infty$ as $n\to \infty$ for all $k>  k_\star$. Given $K>0$, choose $T$ such that $T > \max\{2K, j_{k_\star}(f)\} $, and $\varepsilon = \delta_\star/4$. Then, following similar arguments as before, $j_i(f_n) = j_i(f)$ for $i=1, 2, \ldots, k_\star$, and $j_{k}> T > K$ for all $k> k_\star$ for all $n \ge N(\varepsilon, T)$, but this proves $\lim_{n\to\infty}j_k(f_n)=j_k(f)=\infty$ for all $k> k_\star$.

\end{proof}

\section{Integral probability metrics}
\label{sec:IPM}
Let $(E, d)$ be a Polish metric space. Although the metrics discussed in this section can be defined on arbitrary measurable spaces \cite{Rachev1991ProbabilityMetrics,Zolotarev1983ProbabilityMetrics,Mueller1997IntegralProbMetrics}, it will suffice for our purposes to consider Polish metric spaces. \cite{Sriperumbudur2012EstimationIPM} provides some statistical methods for their empirical estimation. 

\begin{myDefinition}[Integral probability metric]
    Let $\mathbb{F}$ be a class of real-valued measurable functions on $E$. The probability integral semi-metric $\mathsf{d}_{\mathbb{F}}$ corresponding the class $\mathbb{F}$ is defined as 
    \begin{align}
        \mathsf{d}_{\mathbb{F}}(\mu, \nu) \defeq \sup_{f \in \mathbb{F}} \absolute{\measureIntegral{f}{\mu -\nu}}\eqcomma 
    \end{align}
    for two probability measures $\mu$ and $\nu$ on $E$. 
\end{myDefinition}
Integral probability (semi-)metrics are also sometimes referred to as \emph{probability metrics with a $\zeta$-structure} \citep{Zolotarev1983ProbabilityMetrics,Rachev1991ProbabilityMetrics}. Note that $\mathsf{d}_{\mathbb{F}}(\cdot, \cdot)$ is not necessarily finite.

Let $\mathsf{Lip}_1(E, \setOfReals)$ be  the set of real-valued continuous functions defined on $E$ whose best Lipschitz constant is at most 1, \ie, 
\begin{align*}
    \mathsf{Lip}_1(E, \setOfReals) \defeq \{ f \in C(E, \setOfReals) : \sup_{x,y \in E: x\ne y} \frac{\absolute{f(x) - f(y) }}{ d(x, y) } \le 1\}\eqstop 
\end{align*}
The Lipschitz (semi-)norm of a function $f$ on $E$ with respect to the metric $d$ will be denoted by $\norm{f}_{\mathsf{Lip}(E, d)}$. When the space and the metric are clear from the context, we will simply write $\norm{f}_{\mathsf{Lip}}$ for the Lipschitz (semi-)norm of $f$.

\begin{myExample}[Kantorovich--Rubinstein distance]
    Choosing $\mathbb{F} =  \mathsf{Lip}_1(E, \setOfReals)$ yields the Kantorovich--Rubinstein distance, 
    \begin{align}
    \Kantorovich{\mu, \nu} \defeq  \sup \Big\{ \absolute{ \measureIntegral{f}{\mu} - \measureIntegral{f}{\nu} } : f \in  \mathsf{Lip}_1(E, \setOfReals) \Big\} \eqcomma 
    \label{eq:Kantorovich_Rubinstein}
\end{align}
By the Kantorovich--Rubinstein duality, this equals the Vaser\v{s}te\u{\i}n-$1$ distance $\Wasserstein{\mu, \nu}{1}$. 
\label{example:Kantorovich}
\end{myExample}

\begin{myExample}[Fortet--Mourier distance]
The Fortet--Mourier distance is defined as 
    \begin{align}
        \mathsf{d}_{\mathsf{FM}}(\mu, \nu) \defeq  \mathsf{d}_{\mathbb{F}}(\mu, \nu) \eqcomma 
        \label{eq:Fortet_Mourier_metric}
    \end{align}
    for probability measures $\mu, \nu$, with  $\mathbb{F} = \{f : \norm{f}_{c} \le 1\}$, where 
    \begin{align*}
        \norm{f}_{c} \defeq \sup_{x, y \in E : x\ne y}   \frac{\absolute{f(x) - f(y) }}{ c(x, y) }\eqcomma \text{ with } c(x, y) \defeq d(x,y)\max\{1, d(x, z)^{p-1}, d(y, z)^{p-1}\}\eqcomma 
    \end{align*}
    for $p\ge 1$, and for some $z\in E$. 
    \label{example:Fortet_Mourier}
\end{myExample}

\begin{myExample}[Dudley metric]
    The Dudley  (also called bounded Lipschitz) metric is defined as 
    \begin{align}
        \boundedLipschitz{\mu, \nu} \defeq  \mathsf{d}_{\mathbb{F}}(\mu, \nu) \eqcomma 
        \label{eq:Dudley_metric}
    \end{align}
    for probability measures $\mu, \nu$, with  $\mathbb{F} = \{f \in  C(E, \setOfReals):  \norm{f}_{\infty} + \norm{f}_{\mathsf{Lip}} \le 1 \}$. 
    \label{example:DudleyMetric}
\end{myExample}

\begin{myExample}[Total variation metric]
    The total variation metric between two probability measures $\mu$ and $\nu$ on $E$ is given by 
    \begin{align}
        \TotalVariation{\mu, \nu} \defeq \mathsf{d}_{\mathbb{F}}(\mu, \nu) \eqcomma \label{eq:total_variation_defn}
    \end{align}
    with $\mathbb{F} = \{f : E \to \setOfReals : \norm{f}_{\infty} \le 1\}$, or equivalently $\mathbb{F} = \{f : E \to \setOfReals : \sup f - \inf f \le 2\}$ (see \cite[Theorem 5.4]{Mueller1997IntegralProbMetrics}). The quantity $(\sup f - \inf f)$ is called the span of the function $f$. 
    \label{example:TotalVariation}
\end{myExample}

\section{Limits of stochastic processes}
\label{sec:limits_of_stoch_pro}


 

In order to make this manuscript self-contained, we record a few basic probabilistic facts and results. Let $E$ be a complete separable metric space.

\begin{myDefinition}
        For an element $x \in D([0, T], E)$ where $(E, d)$ is a complete, separable metric space, define the modulus of continuity
\begin{align}\label{eq:modu_cont_defn}
    \modulusOfcontinuity(x,T,\delta) \defeq \sup_{\substack{t_1,t_2\in [0,T]\eqcomma \\ |t_1-t_1|\leq \delta}}d(x(t_1), x(t_2))\eqstop 
\end{align}
We will use the same notation $\modulusOfcontinuity$ for different metric spaces $(E, d)$. However, the appropriate metric $d$ to use will be clear from the context. 
\label{defn:modulus_continuity_C}
\end{myDefinition}

\begin{myDefinition}
    A sequence of stochastic processes $\{Z_n\}_{n\ge 1}$ taking values in a metric space $E$ is said to satisfy the \emph{compact containment condition} if for each $\varepsilon>0$ and $T>0$, there exists a compact set $K_\varepsilon \subseteq E$ such that
    \begin{align*}
        \inf_{n} \probOf{ Z_n(t) \in K_\varepsilon, \forall t\in [0, T]} \ge 1-\varepsilon \eqstop
    \end{align*}
    \label{defn:compact_containment_condition}
\end{myDefinition}

\begin{myDefinition}
        A sequence  of stochastic processes $\{Z_n :n \ge 1\}$ is said to be $C$-tight in $D([0, T], \setOfReals)$ if the sequence is relatively compact and hence tight in the space $D([0, T], \setOfReals)$ as stochastic processes, and limit points of every weakly convergent subsequence lie in the space $C([0, T], \setOfReals)$. 
        \label{def:c-tight-in-d}
    \end{myDefinition}

The following theorem, borrowed from \citet{Dawson1993Measure_valued}, gives a criterion for tightness of a sequence of (finite) measure-valued stochastic processes. 

\begin{myTheorem}
Let $E$ be compact, and $\mathcal{H}$ be a dense subset of $C(E, \setOfReals)$ closed under addition. Then, a sequence $\{\mu_n : n \ge 1\}$ of probability measures on $D([0, \infty), \spaceOfMeasures{E})$ is tight in $\spaceOfOccupationMeasures{1}{D([0, \infty), \spaceOfMeasures{E})}$ if and only if $\{\mu_n \circ f_\phi^{-1} : n \ge 1\}$ is tight in $\spaceOfOccupationMeasures{1}{D([0, \infty), \setOfReals) }$ for each $\phi \in \mathcal{H}$, where $f_\phi : D([0, \infty), \spaceOfMeasures{E} )  \to D([0, \infty), \setOfReals) $ is defined by $f_\phi(x)(t) = \measureIntegral{\phi}{x(t)}$ for each $t \in [0, \infty)$, and $x \in D([0, \infty), \spaceOfMeasures{E})$.
\label{thm:tightness_criterion_measure_valued}
\end{myTheorem}
\begin{proof}[Proof of \Cref{thm:tightness_criterion_measure_valued}]
    See \citet[Theorem 3.7.1 (a)]{Dawson1993Measure_valued}. 
\end{proof}



\section{Propagation of chaos}
\label{sec:prop_chaos_appendix}
As before, we assume $E$ is a separable metric space. We equip the set $E$ with the Borel $\sigma$-algebra $\borel{E}$. The following definition is standard (see, e.g., \cite{Sznitman1991Topics,Hauray2014KacChaos}). 

\begin{myDefinition}
    Let $\{\mu_n : n\ge 1\}$ be a sequence of symmetric probability measures on $(E^n, \borel{E^n})$. Let $\mu$ be a probability measure on $(E, \borel{E})$. We say that $\{\mu_n : n\ge 1\}$ is $\mu$-chaotic if, for each $k\ge 1$, the following limit holds:
    \begin{align}
        \lim_{n\to \infty} \int_{E^n} \left(\prod_{i=1}^k f_i(x_i) \right) \mu_n(\d x_1, \ldots, \d x_n) = \prod_{i=1}^{k} \left( \int_E f_i(x) \mu(\d x)\right)\eqcomma 
        \label{eq:prop_chaos_defn_1}
    \end{align}
    for $f_i \in C_b(E, \setOfReals)$, $i=1,\ldots,k$. 
    \label{def:prop_chaos}
\end{myDefinition}

An equivalent definition of propagation of chaos is as follows. 

\begin{myDefinition}
    Let $\{\mu_n : n\ge 1\}$ be a sequence of symmetric probability measures on $(E^n, \borel{E^n})$. Let $\mu$ be a probability measure on $(E, \borel{E})$. We say that $\{\mu_n : n\ge 1\}$ is $\mu$-chaotic if, for each $k\ge 1$, the following limit holds: 
    \begin{align}
        (\nU_1, \nU_2, \ldots, \nU_k) \ConvInDist(U_1, U_2, \ldots, U_k) 
        \label{eq:prop_chaos_defn_2}
    \end{align}
    as $n\to \infty$, where the  law of the $E^n$-valued random variable $(\nU_1, \nU_2, \ldots, \nU_n)$ is $\mu_n$, and $U_1, U_2, \ldots, U_k$ are \ac{iid} $E$-valued random variables with the probability law $\mu$. 
    \label{def:prop_chaos_2}
\end{myDefinition}

When the convergence in \eqref{eq:prop_chaos_defn_2} holds in a stronger notion (e.g., in $r$-th moment), we define propagation of chaos analogously. The following useful result is borrowed from \citet{Sznitman1991Topics}. 


\begin{myTheorem}
    Let $\{\mu_n : n\ge 1\}$ be a sequence of symmetric probability measures on $(E^n, \borel{E^n})$. Let $\mu$ be a probability measure on $(E, \borel{E})$. Then, $\{\mu_n : n\ge 1\}$ is $\mu$-chaotic if and only if the empirical random measure $$\tilde{\mu}_n \defeq \frac{1}{n}\sum_{i=1}^n \delta_{\nU_i}$$ converges weakly to the deterministic probability measure $\mu$ as $n\to \infty$, where the $E^n$-valued random variable $(\nU_1, \nU_2, \ldots, \nU_n)$ is distributed according to the probability measure $\mu_n$. Moreover, it is also equivalent to the convergence in \eqref{eq:prop_chaos_defn_1} (or equivalently, in \eqref{eq:prop_chaos_defn_2}) for $k=2$. 
    \label{thm:prop_chaos_empirical_measure}
\end{myTheorem}
\begin{proof}[Proof of \Cref{thm:prop_chaos_empirical_measure}]
See \citet[Proposition 2.2]{Sznitman1991Topics}.
\end{proof}

The state of the art on the theory of quantitative propagation of chaos in the context of \acp{IPS} allows the number $k$ to vary with $n$. However, we do not discuss this aspect here, and refer the readers to \cite{Hauray2014KacChaos} for more details. 


\section{Acronyms}
\label{sec:acronyms}
    
\begin{acronym}[OWL-QN]
    \acro{a.s.}{almost surely}
	\acro{BDG}{Burkholder--Davis--Gundy}
	\acro{CDC}{Centers for Disease Control and Prevention}
	\acro{CDF}{Cumulative Distribution Function}
	\acro{CDM}{Compartmental Disease Model}
	\acro{CLT}{Central Limit Theorem}
	\acro{CM}{Configuration Model}
	\acro{CME}{Chemical Master Equation}
	\acro{CRN}{Chemical Reaction Network}
	\acro{CTMC}{Continuous Time Markov Chain}
	\acro{DCT}{Dominated Convergence Theorem}
	\acro{DSA}{Dynamic Survival Analysis}
	\acro{DTMC}{Discrete Time Markov Chain}
	\acro{ER}{Erd\"{o}s--R\'{e}nyi}
	\acro{FCLT}{Functional Central Limit Theorem}
	\acrodefplural{FCLT}[FCLTs]{Functional Central Limit Theorems}
	\acro{FLLN}{Functional Law of Large Numbers}
	\acrodefplural{FLLN}[FLLNs]{Functional Laws of Large Numbers}
	\acro{GP}{Gaussian Process}
	\acrodefplural{GP}[GPs]{Gaussian Processes}
	\acro{HJB}{Hamilton–Jacobi–Bellman}
	\acro{HMC}{Hamiltonian Monte Carlo}
	\acro{iid}{independent and identically distributed}
	\acro{IPS}{Interacting Particle System}
    \acro{IVP}{Initial Value Problem}
	\acro{KL}{Kullback-Leibler}
	\acro{LDP}{Large Deviations Principle}
	\acro{LLN}{Law of Large Numbers}
	\acrodefplural{LLN}[LLNs]{Laws of Large Numbers}
	\acro{LNA}{Linear Noise Approximation}
	\acro{MCMC}{Markov Chain Monte Carlo}
	\acro{MFPT}{Mean First Passage Time}
	\acro{MGF}{Moment Generating Function}
	\acro{MLE}{Maximum Likelihood Estimate}
	\acro{MM}{Michaelis--Menten}
	\acro{ODE}{Ordinary Differential Equation}
	\acro{OU}{Ornstein--Uhlenbeck}
	\acro{PDE}{Partial Differential Equation}
	\acro{PDF}{Probability Density Function}
	\acro{PDP}{Piecewise Deterministic Process}
	\acro{PGF}{Probability Generating Function}
	\acro{PMF}{Probability Mass Function}
    \acro{PPM}{Poisson Point Measure}
    \acro{PRM}{Poisson Random Measure}
	\acro{psd}{positive semi-definite}
	\acro{PT}{Poisson-type}
	\acro{QSSA}{Quasi-Steady State Approximation}
	\acro{rQSSA}{reversible QSSA}
	\acro{SD}{Standard Deviation}
	\acro{SDE}{Stochastic Differential Equation}
	\acro{SEIR}{Susceptible-Exposed-Infected-Recovered}
	\acro{SI}{Susceptible-Infected}
	\acro{SIR}{Susceptible-Infected-Recovered}
	\acro{SIS}{Susceptible-Infected-Susceptible}
	\acro{SLLN}{Strong Law of Large Numbers}
	\acro{sQSSA}{standard QSSA}
	\acro{tQSSA}{total QSSA}
	\acro{TK}{Togashi--Kaneko}
	\acro{WS}{Watts--Strogatz}
	\acro{whp}{with high probability}
\end{acronym}

%
%

\section*{Acknowledgements}
The authors thank Gabriela Gomes (University of Strathclyde), Eben Kenah (The Ohio State University), Thomas House (University of Manchester), Lorenzo Pellis (University of Manchester) for helpful discussions at various stages of this work.


\section*{Funding}
Olga Izyumtseva was supported by the British Academy through grant number RaR{\textbackslash}100741, and in part by British Academy, Cara, Leverhulme Trust through grant LTRSF24{\textbackslash}100014 and LTRSF26{\textbackslash}100038. Kushankur Dutta acknowledges the PhD studentship at the University of Nottingham.


\section*{Roles}
All authors made significant contributions to the work. Following standard conventions in mathematics, the authors are listed alphabetically.





\bibliographystyle{plainnat}
\bibliography{references}

\end{document}